\documentclass[fleqn,10pt]{article}
\usepackage{latexsym, graphicx, epsfig, amsmath, amssymb,amsfonts}
\usepackage{natbib,amsthm,version}
\usepackage{amsbsy,bm,multirow,enumerate}
\usepackage[titletoc,page]{appendix}
\usepackage[mathscr]{eucal}
\usepackage{mathtools}
\usepackage{color}
\usepackage[utf8]{inputenc}
\usepackage[english]{babel}
\usepackage{amsthm}
\usepackage{enumerate}
\usepackage[hidelinks]{hyperref}
\usepackage{url}
\usepackage{subfigure}
\usepackage[lined,boxed,linesnumbered,ruled]{algorithm2e}
\usepackage{float}
\usepackage{setspace}

\usepackage{natbib}
\newcommand{\bs}[1]{\boldsymbol{#1}}
\newcommand{\mbs}[1]{\mathbf{#1}}
\newcommand{\mbb}[1]{\mathbb{#1}}

\DeclareMathOperator*{\argmin}{arg\,min}

\newtheorem{remark}{Remark}[section]

\theoremstyle{definition}

\title{
  A Method for Computing Inverse Parametric PDE Problems
  with Random-Weight Neural Networks
} 
\author{
  Suchuan Dong\thanks{Author of correspondence. Emails: sdong@purdue.edu (S.~Dong),
    wang2335@purdue.edu (Y.~Wang).}, \ \ Yiran Wang  \\
  Center for Computational and Applied Mathematics \\
  Department of Mathematics \\
  Purdue University, USA 
 } 

\date{(October 8, 2022)}

\begin{document}
\maketitle


\begin{abstract}

  We present a method for computing the inverse parameters
  and the solution field
  to inverse parametric partial differential equations (PDE) based on
  randomized neural networks. This extends the local extreme learning machine
  technique originally developed for forward PDEs to inverse problems.
  We develop three algorithms for training the neural network to solve the inverse PDE
  problem. The first algorithm (termed NLLSQ) determines the inverse parameters
  and the trainable network parameters all together by
  the nonlinear least squares method with perturbations (NLLSQ-perturb).
  The second algorithm (termed VarPro-F1) eliminates the inverse parameters
  from the overall problem by variable projection to attain
  a reduced problem about the trainable network parameters only.
  It solves the reduced problem first by the NLLSQ-perturb algorithm for
  the trainable network parameters, and then computes the inverse parameters
  by the linear least squares method. The third algorithm (termed VarPro-F2)
  eliminates the trainable network parameters from the overall
  problem by variable projection to attain a reduced problem
  about the inverse parameters only.
  It solves the reduced problem  for
  the inverse parameters first, and then computes the trainable network parameters 
  afterwards. VarPro-F1 and VarPro-F2
  are reciprocal to each other in some sense.
  The presented method produces accurate results for inverse PDE problems,
  as shown by the numerical examples herein. For noise-free data, the errors of
  the inverse parameters and the solution field decrease exponentially as
  the number of collocation points or the number of trainable network parameters increases,
  and can reach a level close to the machine accuracy.
  For noisy data, the accuracy degrades compared with the case of noise-free data,
  but the method remains quite accurate. The presented method has been compared with
  the physics-informed neural network method.

\end{abstract}


\vspace{0.05cm}




\section{Introduction}
\label{sec:intro}


In this work we focus on the simultaneous determination of the parameters (as
constants or field distributions)
and  the solution field to parametric PDEs
based on artificial neural networks (ANN/NN),
given sparse and noisy measurement data of certain variables.
This type of problems is often referred to
as the inverse PDE problems in the literature~\cite{Karniadakisetal2021}.
Typical examples include the determination of the diffusion
coefficient given certain concentration data
or the computation of the wave speed given sparse measurement of the wave profile.
When the parameter values in the PDE are known, approximation of the PDE solution
is often referred to as the forward PDE problem. We will adopt these notations in this paper.


Closely related to the inverse PDE problems is the data-driven ``discovery'' of 
PDEs (see e.g.~\cite{BongardL2007,SchmidtL2009,BruntonPK2016,RudyBPK2017,Schaeffer2017,RaissiK2018,ZhangL2018,BergN2019,LongLD2019,RudyABK2019,WuX2020,BothCSK2021,XuZZ2021,NorthWS2022}, among others),
in which, given certain
measurement data, the functional form
of the PDE is to be discerned. In order to acquire a parsimonious PDE form,
techniques such as sparsity promotion~\cite{BruntonPK2016,RudyBPK2017} or dimensional
analysis~\cite{ZhangL2018} are often adopted.

As advocated in~\cite{Tartakovskyetal2020,Karniadakisetal2021},
data-driven scientific machine learning problems can be viewed in terms of
the amount of data that is available and the amount of physics that is known.
They are broadly classified into three categories in~\cite{Karniadakisetal2021}:
(i) those with ``lots of physics and small data'' (e.g.~forward PDE problems),
(ii) those with ``some physics and some data'' (e.g.~inverse PDE problems),
and (iii) those with ``no physics and big data'' (e.g.~general PDE discovery).
The authors of~\cite{Karniadakisetal2021} point out that those in the second category
are typically the more interesting and representative in
real applications, where the physics is partially known and sparse measurements
are available. One illustrating example is from multiphase flows,
where the conservation laws (mass/momentum conservations) and
thermodynamic principles (second law of thermodynamics, Galilean invariance)
lead to a thermodynamically-consistent phase field model, but with an incomplete system of
governing equations~\cite{Dong2018,Dong2017}. One has the freedom to choose
the form of the free energy, the wall energy, the form and coefficients
of the constitutive relation, and the
form and coefficient of the interfacial mobility~\cite{Dong2014,Dong2015,YangD2018}.
Different choices will lead to different specific models, which are all thermodynamically consistent.
The different models cannot be distinguished by the thermodynamic principles,
but can be differentiated with experimental measurements.


The development of machine learning techniques for solving
inverse PDE problems   has attracted
a great deal of interest recently, with a variety of contributions from
different researchers.
In~\cite{RaissiK2018} a method for estimating the parameters in nonlinear PDEs
is developed based on Gaussian processes, where the state variable at two consecutive
snapshots are assumed to be known.
The physics informed neural network (PINN) method is introduced in
the influential work~\cite{RaissiPK2019} for solving forward and inverse nonlinear PDEs.
The residuals of the PDE, the boundary and initial conditions,
and the measurement data are encoded into the loss function as soft constraints,
and the neural network is trained by gradient descent or back propagation type algorithms.
The PINN method has influenced significantly subsequent developments
and stimulated applications in many related areas
(see e.g.~\cite{MaoJK2020,RaissiYK2020,JagtapKK2020,MengK2020,Chenetal2020,Schiassietal2020,Caietal2021,LuMMK2021,YangMK2021,JagtapMAK2022,Pateletal2022}, among others).
A hybrid finite element and neural network method is developed in~\cite{BergN2018inv}.
The finite element method (FEM) is used to solve the underlying PDE, which is augmented
by a neural network for representing the PDE coefficient~\cite{BergN2018inv}.
A conservative PINN method is proposed in~\cite{JagtapKK2020} together with
domain decomposition for simulating nonlinear conservation laws, in which the flux
continuity is enforced along the sub-domain interfaces, and interesting results are presented
for a number of forward and inverse problems.
This method is further developed and extended in a subsequent work~\cite{JagtapK2020}
with domain decompositions in both space and time; see a recent study in~\cite{JagtapMAK2022}
of this extended technique for supersonic flows.
Interesting applications are described in~\cite{RaissiYK2020,Caietal2021}, where
the PINN technique is employed to infer the 3D velocity and pressure
fields based on scattered flow visualization data or Schlieren images
from experiments.
In~\cite{DwivediPS2021} a distributed PINN method based on domain decomposition
is presented, and the loss function is optimized by a gradient descent algorithm.
For nonlinear PDEs, the method solves a related linearized equation with certain
variables fixed at their initial values~\cite{DwivediPS2021}.
An auxiliary PINN technique is developed in~\cite{YuanNDH2022} for solving
nonlinear integro-differential equations, in which auxiliary variables are introduced
to represent the anti-derivatives and thus avoiding the integral computation.
We would also like to refer the reader
to e.g.~\cite{Chenetal2020,Tartakovskyetal2020,MathewsFHH2020,LiXHD2020}
(among others) for inverse applications
of neural networks in other related fields.



In the current work we consider the use of randomized neural networks,
also known as extreme learning machines (ELM)~\cite{HuangZS2006} (or
random vector functional link (RVFL) networks~\cite{PaoPS1994}),
for solving inverse PDE problems.
ELM was originally developed for linear classification and regression problems.
It is characterized by two ideas: (i) randomly assigned but fixed (non-trainable)
hidden-layer coefficients, and (ii) trainable linear output-layer coefficients
determined by linear
least squares or by using the
Moore-Penrose inverse~\cite{HuangZS2006}. 
This technique has been extended to scientific computing
in the past few years, for function approximations and for solving
 ordinary and partial differential equations (ODE/PDE);
see e.g.~\cite{YangHL2018,PanghalK2020,DwivediS2020,DongL2021,DongL2021bip,CalabroFS2021,FabianiCRS2021,Schiassietal2021,DongY2022rm}, among others.
The random-weight neural networks are universal function approximators.
As established by the theoretical results of~\cite{IgelnikP1995,HuangCS2006,NeedellNSS2020},
a single-hidden-layer feed-forward neural network (FNN)
having random but fixed (not trained) hidden units
can approximate any continuous function to any desired degree of accuracy,
provided that the number of hidden units is sufficiently large.

In this paper we present a method for computing inverse PDE problems
based on randomized neural networks.
This extends the local extreme learning machine (locELM) technique 
originally developed in~\cite{DongL2021} for forward PDEs
to inverse problems.
Because of the coupling between the unknown PDE parameters
(referred to as the inverse parameters
hereafter) and the solution field, the inverse PDE problem is
fully nonlinear with respect to the unknowns, even though the associated forward PDE
may be linear. We partition the overall domain into sub-domains, and represent
the solution field (and the inverse parameters, if they are field distributions)
by a local FNN  on each sub-domain, imposing
$C^{\mbs k}$ (with appropriate $\mbs k$)
continuity conditions across
the sub-domain boundaries. The weights/biases in the hidden layers of
the local NNs are assigned to random values and fixed (not trainable),
and only the output-layer coefficients are trainable.
The inverse PDE problem is thus reduced to a nonlinear problem about
the inverse parameters and the output-layer coefficients of the solution field,
or if the inverse parameters are field distributions, about the output-layer
coefficients for the inverse parameters and the solution field.

We develop three algorithms for training the neural network to
solve the inverse PDE problem:
\begin{itemize}
\item
  The first algorithm (termed NLLSQ) computes the inverse parameters
  and the trainable parameters of the local NNs all together
  by the nonlinear least squares method~\cite{Bjorck1996}. This extends the
  nonlinear least squares method with perturbations (NLLSQ-perturb)
  from~\cite{DongL2021} (developed for forward nonlinear PDEs)
  to inverse PDE problems.

\item
  The second algorithm (termed VarPro-F1) eliminates the inverse parameters
  from the overall problem based on the variable projection (VarPro)
  strategy~\cite{GolubP1973,GolubP2003}
  to attain a reduced problem about the trainable network parameters
  only. It solves the reduced problem first for
  the trainable parameters of the local NNs
  by the NLLSQ-perturb algorithm, and
  then computes the inverse parameters by the linear least squares method.

\item
  The third algorithm (termed VarPro-F2) eliminates the trainable network parameters
   from the overall inverse problem by
  variable projection to arrive at a reduced problem about the
  inverse parameters only. It solves the reduced problem first for the inverse parameters
  by the NLLSQ-perturb algorithm, and then computes the trainable parameters
  of the local NNs based on the inverse parameters already obtained.
  The VarPro-F2 and VarPro-F1 algorithms both employ the variable projection
  idea and are reciprocal formulations in a sense.
  For inverse problems with an associated forward nonlinear PDE, VarPro-F2
  needs to be combined with a Newton iteration.
  
\end{itemize}

The presented method  produces
accurate solutions to inverse PDE problems, as shown by a number of
numerical examples presented herein. For noise-free data, the errors for the inverse parameters
and the solution field decrease exponentially as the number of
training collocation points or the number of trainable parameters
in the neural network increases. These errors can reach a level close to
the machine accuracy when the simulation parameters become large.
For noisy data, the current method remains quite accurate, although the accuracy
degrades compared with the case of noise-free data.
We observe that, by scaling the measurement-residual vector by a factor, one
can markedly improve the accuracy of the current method for noisy data,
while only slightly degrading the accuracy for noise-free data.
We have compared the current method with the PINN method (see Appendix C).
The current method exhibits an advantage
in terms of the accuracy and the computational cost (network training time).

The method and algorithms developed herein are implemented in
Python based on the Tensorflow (https://www.tensorflow.org/),
Keras (https://keras.io/), and the scipy (https://scipy.org/) libraries.
The numerical simulations are performed on a MAC computer (3.2GHz Intel Core i5 CPU,
24GB memory) in the authors' institution.


The main contribution of this paper lies in the local extreme learning machine
based technique together with the three algorithms
for solving inverse PDE problems.
The exponential convergence behavior  exhibited by the current method
for inverse problems is particularly interesting, and
can be analogized to the observations in~\cite{DongL2021} for forward PDEs.
For inverse problems
such fast convergence seems not available
in the existing techniques (e.g.~PINN based methods).


The rest of this paper is structured as follows. 
In Section~\ref{sec:method} we first discuss the representation of
functions by local randomized neural networks and domain decomposition,
and then present the NLLSQ, VarPro-F1 and VarPro-F2 algorithms
for training the neural network to solve the inverse PDE.
Section~\ref{sec:tests} uses a number of  inverse
parametric PDEs to demonstrate the exponential convergence and
the accuracy of our method, as well as the effects of the noise
and the number of measurement points.
Section~\ref{sec:summary} concludes the discussion with
some closing remarks.
Appendix A summarizes the NLLSQ-perturb algorithm from~\cite{DongL2021}
(with modifications), which forms the basis for the three algorithms
in the current paper for solving inverse PDEs.
Appendix B provides the matrices in the VarPro-F2
algorithm. Appendix C compares the current method with
PINN for several inverse problems
from Section~\ref{sec:tests}.
Appendix D lists the parameter values in the NLLSQ-perturb
algorithm for all the numerical simulations in Section~\ref{sec:tests}.



\section{Algorithms for Inverse PDEs with Randomized Neural Networks}
\label{sec:method}

\subsection{Inverse Parametric PDEs and Local Randomized Neural Networks}


We focus on the inverse problem described by the following parametric PDE,
boundary conditions, and measurement operations on some
domain $\Omega\subset\mathbb{R}^d$ ($d=1,2,3$): 
\begin{subequations}\label{eq_1}\small
  \begin{align}
    &
    \alpha_1\mathcal{L}_1(u) + \alpha_2\mathcal{L}_2(u) + \dots + \alpha_n\mathcal{L}_n(u)
    + \mathcal{F}(u) = f(\mbs x), \quad \mbs x\in\Omega,
    \label{eq_1a} \\
    &
    \mathcal{B}u(\mbs x) = g(\mbs x), \quad \mbs x\in \partial\Omega, \label{eq_1b} \\
    &
    \mathcal{M}u(\bm\xi) = S(\bm\xi), \quad \bm\xi\in\Omega_s\subset\Omega.
  \end{align}
\end{subequations}
In this system, $\mathcal{L}_i$ ($1\leqslant i\leqslant n$)
and $\mathcal{F}$ are differential or algebraic
operators, which can be linear or nonlinear, and $f$ and $g$ are prescribed source terms.
$u(\mbs x)$ is an unknown scalar field, where $\mbs x$ denotes the coordinates.
$\alpha_i$ ($1\leqslant i\leqslant n$) are $n$ unknown constants.
The case with any $\alpha_i$ being an unknown field distribution will be dealt with
later in a remark (Remark~\ref{rem_05}).
We assume that the highest derivative term in~\eqref{eq_1a} is linear with
respect to $u$, while the nonlinear terms with respect to $u$ involve only lower
derivatives (if any).
$\mathcal{B}$ is a linear differential or algebraic operator,
and $\mathcal{B}u$ denotes the boundary
condition(s) on the domain boundary $\partial\Omega$.
$\mathcal{M}$ is a linear algebraic or differential operator
representing the measurement operations.
$\mathcal{M}u(\bm\xi)$ denotes the measurement of $\mathcal{M}u$ at the point $\bm\xi$,
and $S(\bm\xi)$ denotes the measurement data. $\Omega_s$ denotes the set of measurement points.
Given $S(\bm\xi)$, the goal here is to
determine the parameters $\alpha_i$ ($1\leqslant i\leqslant n$)
and the solution field $u(\mbs x)$.
Hereafter we will refer to the parameters $\bm\alpha=(\alpha_1,\dots,\alpha_n)^T$ as
the inverse parameters.
Suppose the inverse parameters are given. The boundary value
problem consisting of the equations~\eqref{eq_1a}--\eqref{eq_1b} will be
referred to as the associated forward PDE problem, with $u(\mbs x)$
as the unknown. We assume that
the formulation is such that the forward PDE problem is well-posed.

\begin{remark}\label{rem_1}
  We assume that the operators $\mathcal{L}_i$ ($1\leqslant i\leqslant n$)
  or $\mathcal{F}$ may contain
  time derivatives (e.g.~$\frac{\partial}{\partial t}$, $\frac{\partial^2}{\partial t^2}$, where
  $t$ denotes time), thus leading to an initial-boundary value problem
  on a spatial-temporal domain $\Omega$.
  In this case, we treat $t$ in the same way as the spatial coordinate $\mbs x$, and
  use the last dimension in $\mbs x=(x_1,x_2,\dots,x_d)$ to denote $t$ (i.e.~$x_d\equiv t$).
  Accordingly, we assume that the equation~\eqref{eq_1b} should include 
  conditions on the appropriate initial boundaries from $\partial\Omega$.
  The point here is that the system~\eqref{eq_1} may refer to time-dependent problems,
  and we will not distinguish this case in subsequent discussions.

\end{remark}

\begin{figure}
  \small
  \centerline{
    \includegraphics[width=2.2in]{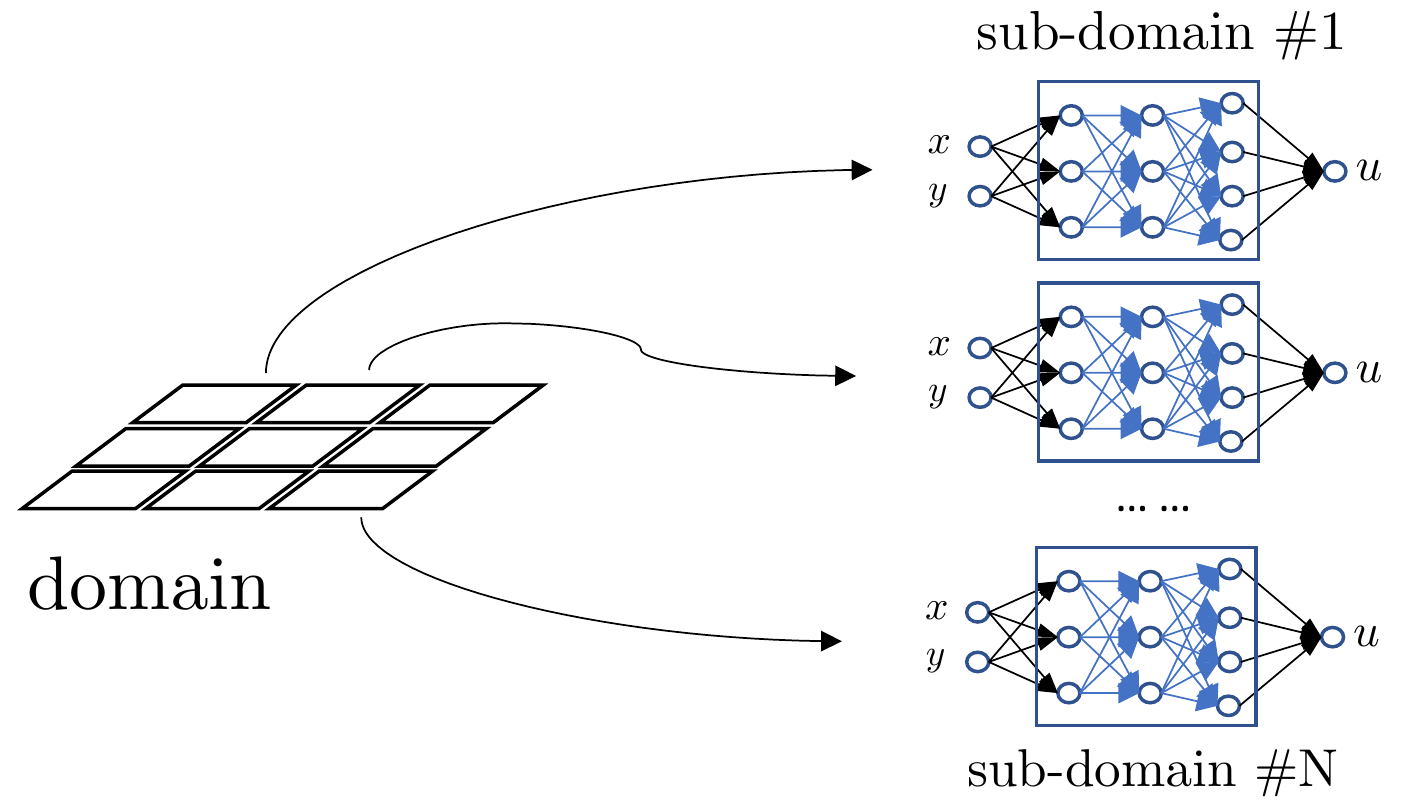}
  }
  \caption{\small Cartoon illustrating domain decomposition and local
    random-weight neural networks.}
  \label{fg_1}
\end{figure}

We devise numerical algorithms to compute a least squares solution
to the system~\eqref{eq_1} based on local randomized  neural networks (or ELM).
%
We decompose the domain $\Omega$ into sub-domains, and represent
$u(\mbs x)$ on each sub-domain by a local ELM in a way
analogous to in~\cite{DongL2021}.
Let
$ 
  \Omega = \Omega_1\cup\Omega_2\cup\dots\cup\Omega_{N},
$ 
  where $\Omega_i$ ($1\leqslant i\leqslant N$) denote $N$
  non-overlapping sub-domains
  (see Figure~\ref{fg_1} for an illustration).
Let
\begin{equation}\label{eq_2}\small
  u(\mbs x) = \left\{
  \begin{array}{ll}
    u_1(\mbs x), & \mbs x\in\Omega_1, \\
    u_2(\mbs x), & \mbs x\in \Omega_2, \\
    \dots \\
    u_N(\mbs x), & \mbs x\in \Omega_{N},
  \end{array}
  \right.
\end{equation}
where $u_i(\mbs x)$ ($1\leqslant i\leqslant N$) denotes the solution field
restricted to the sub-domain $\Omega_i$.
On the interior sub-domain boundaries shared by adjacent sub-domains
we impose $C^{\mbs k}$
continuity conditions on  $u(\mbs x)$,
where $\mbs k=(k_1,\dots,k_d)$ denotes a set of appropriate
non-negative integers related to the order of the PDE~\eqref{eq_1a}. 
If the PDE order (highest derivative) is $m_i$ along the $x_i$ ($1\leqslant i\leqslant d$)
direction, we would in general
impose $C^{m_i-1}$ (i.e.~$k_i=m_i-1$) continuity conditions in this direction
on the shared sub-domain boundaries.

On $\Omega_i$ ($1\leqslant i\leqslant N$) we employ
a local FNN,  whose hidden-layer coefficients
are randomly assigned and fixed, to represent $u_i(\mbs x)$.
More specifically, the local neural network is set as follows. 
The input layer consists of $d$ nodes, representing the input coordinate
$\mbs x=(x_1,x_2,\dots,x_d)\in\Omega_i$. The output layer consists of
a single node, representing $u_i(\mbs x)$.
The network contains $(L-1)$ (with integer $L\geqslant 2$) hidden layers
in between.
Let $\sigma: \mathbb{R}\rightarrow\mathbb{R}$ denote the activation
function for all the hidden nodes.
Hereafter we use the following vector (or list) $\mbs M$ of $(L+1)$ positive integers
to represent the architecture of the local NN,
\begin{equation}\small
  \mbs M=\begin{bmatrix} m_0, m_1,\dots, m_{L-1}, m_L  \end{bmatrix},
  \quad \text{(architectural vector)}
\end{equation}
where $m_0=d$ and $m_L=1$ denote the number of nodes in the input/output
layers respectively, and $m_i$ is
the number of nodes in the $i$-th hidden layer
($1\leqslant i\leqslant L-1$).
We refer to $\mbs M$ as an architectural vector.

We make the following assumptions:
\begin{itemize}
\item
  The output layer should contain (i) no bias, and (ii) no activation function
  (or equivalently, the activation function be $\sigma(x)=x$).

\item
  The weights/biases in all the hidden layers are pre-set to
  uniform random values on $[-R_m,R_m]$,
  where $R_m>0$ is a user-provided constant. The hidden-layer
  coefficients are fixed once they are set.

\item
  The output-layer weights constitute the
  the trainable parameters of the local neural network.

\end{itemize}
We employ the same architecture,
same activation function, and the same $R_m$
for the local neural networks on different sub-domains.

In light of these settings, the logic in the output layer
of the local NNs leads to the following
relation on the sub-domain $\Omega_i$ ($1\leqslant i\leqslant N$),
\begin{equation}\label{eq_4}\small
  u_i(\mbs x) = \sum_{j=1}^M \beta_{ij}\phi_{ij}(\mbs x) = \bm\Phi_i(\mbs x)\bm\beta_i,
\end{equation}
where $M=m_{L-1}$ denotes the width of the last hidden layer of
the local NN, $\phi_{ij}(\mbs x)$ ($1\leqslant j\leqslant M$)
denote the set of output fields of the last hidden layer on $\Omega_i$,
$\beta_{ij}$ ($1\leqslant j\leqslant M$) denote the set of output-layer
coefficients (trainable parameters) on $\Omega_i$,
and $\bm\Phi_i = (\phi_{i1},\phi_{i2},\dots,\phi_{iM})$ and
$\bm\beta_i = (\beta_{i1},\beta_{i2},\dots,\beta_{iM})^T$.
Note that, once the random hidden-layer coefficients
are assigned, $\bm\Phi_i(\mbs x)$ in~\eqref{eq_4}
denotes a set of random (but fixed and known) nonlinear basis functions.
Therefore, with local ELMs 
the output field on each sub-domain is represented by an expansion of
a set of random basis functions as given by~\eqref{eq_4}.

With domain decomposition and local ELMs,
the system~\eqref{eq_1} is symbolically transformed into the following form,
which includes the continuity conditions across shared sub-domain boundaries:
\begin{subequations}\label{eq_5}\small
  \begin{align}
    &
    \alpha_1\mathcal{L}_1(u_i) + \alpha_2\mathcal{L}_2(u_i) + \dots + \alpha_n\mathcal{L}_n(u_i)
    + \mathcal{F}(u_i) = f(\mbs x), \quad \mbs x\in\Omega_i, \ 1\leqslant i\leqslant N;
    \label{eq_5a} \\
    &
    \mathcal{B}u_i(\mbs x) = g(\mbs x), \quad \mbs x\in \partial\Omega\cap\Omega_i,
    \ 1\leqslant i\leqslant N;
    \label{eq_5b} \\
    &
    \mathcal{M}u_i(\bm\xi) = S(\bm\xi), \quad \bm\xi\in\Omega_s\cap\Omega_i,
    \ 1\leqslant i\leqslant N;  \label{eq_5c} \\
    &
    \mathcal{C}u_i(\mbs x) - \mathcal{C}u_j(\mbs x)=0, \quad
    \mbs x\in\partial\Omega_i\cap\partial\Omega_j,
    \ \text{for all adjacent sub-domains}\ (\Omega_i,\Omega_j),
    \ \ 1\leqslant i,j\leqslant N.
    \label{eq_5d}
  \end{align}
\end{subequations}
In this system $u_i(\mbs x)$ is given by \eqref{eq_4}, and
the operator $\mathcal{C}u$
denotes the set of $C^{\mbs k}$ continuity conditions imposed across the shared
sub-domain boundaries on $u$ or its derivatives.
Define the residual of this system as,
\begin{equation}\label{eq_6}\small
  \mbs R(\bm\alpha,\bm\beta,\mbs x,\bm\xi)=\begin{bmatrix}
  \alpha_1\mathcal{L}_1(u_i) + \alpha_2\mathcal{L}_2(u_i)
  + \dots + \alpha_n\mathcal{L}_n(u_i)
  + \mathcal{F}(u_i) - f(\mbs x),
  \ \mbs x\in\Omega_i,\ 1\leqslant i\leqslant N
  \\
  \mathcal{B}u_i(\mbs x) - g(\mbs x),
  \  \mbs x\in\partial\Omega\cap\Omega_i,\ 1\leqslant i\leqslant N
  \\
  \mathcal{M}u_i(\bm\xi) - S(\bm\xi),
  \ \bm\xi\in\Omega_s\cap\Omega_i,\ 1\leqslant i\leqslant N
  \\
  \mathcal{C}u_i(\mbs x) - \mathcal{C}u_j(\mbs x),
  \ \mbs x\in\partial\Omega_i\cap\partial\Omega_j,
  \ \text{for all adjacent}\ (\Omega_i,\Omega_j),\ 1\leqslant i,j\leqslant N
  \end{bmatrix},
\end{equation}
where $\bm\beta$ is the vector of all trainable parameters,
$\bm\beta=(\bm\beta_1^T,\dots,\bm\beta_N^T)^T
=(\beta_{11},\beta_{12},\dots,\beta_{1M},\beta_{21},\dots,\beta_{NM})^T$.

The system~\eqref{eq_5} is what we would solve
numerically by least squares
for the inverse parameters $\bm\alpha$
and the trainable network parameters $\bm \beta$.
After $(\bm\alpha,\bm\beta)$ are determined, the field
solution $u(\mbs x)$ is computed by~\eqref{eq_2} and~\eqref{eq_4}.
In what follows we present three algorithms, one based on
the nonlinear least squares method with perturbations and
the other two based on the variable projection idea,
for determining the $\bm\alpha$ and $\bm\beta$.

\subsection{Nonlinear Least Squares (NLLSQ) Method for Network Training}
\label{sec:nllsq}

\begin{figure}
  \small
  \centerline{
    \includegraphics[width=1.4in]{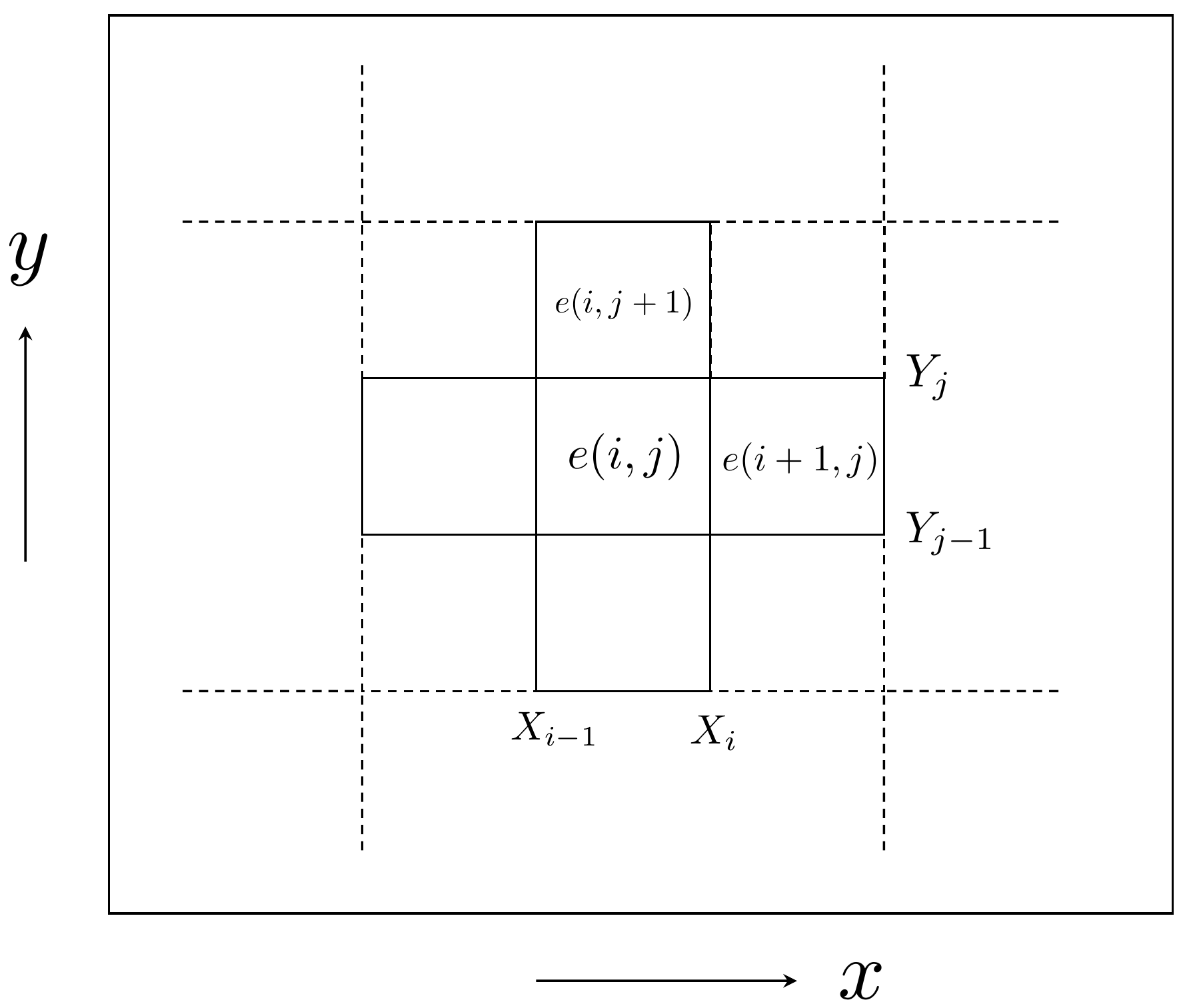}(a) \qquad
    \includegraphics[width=1.4in]{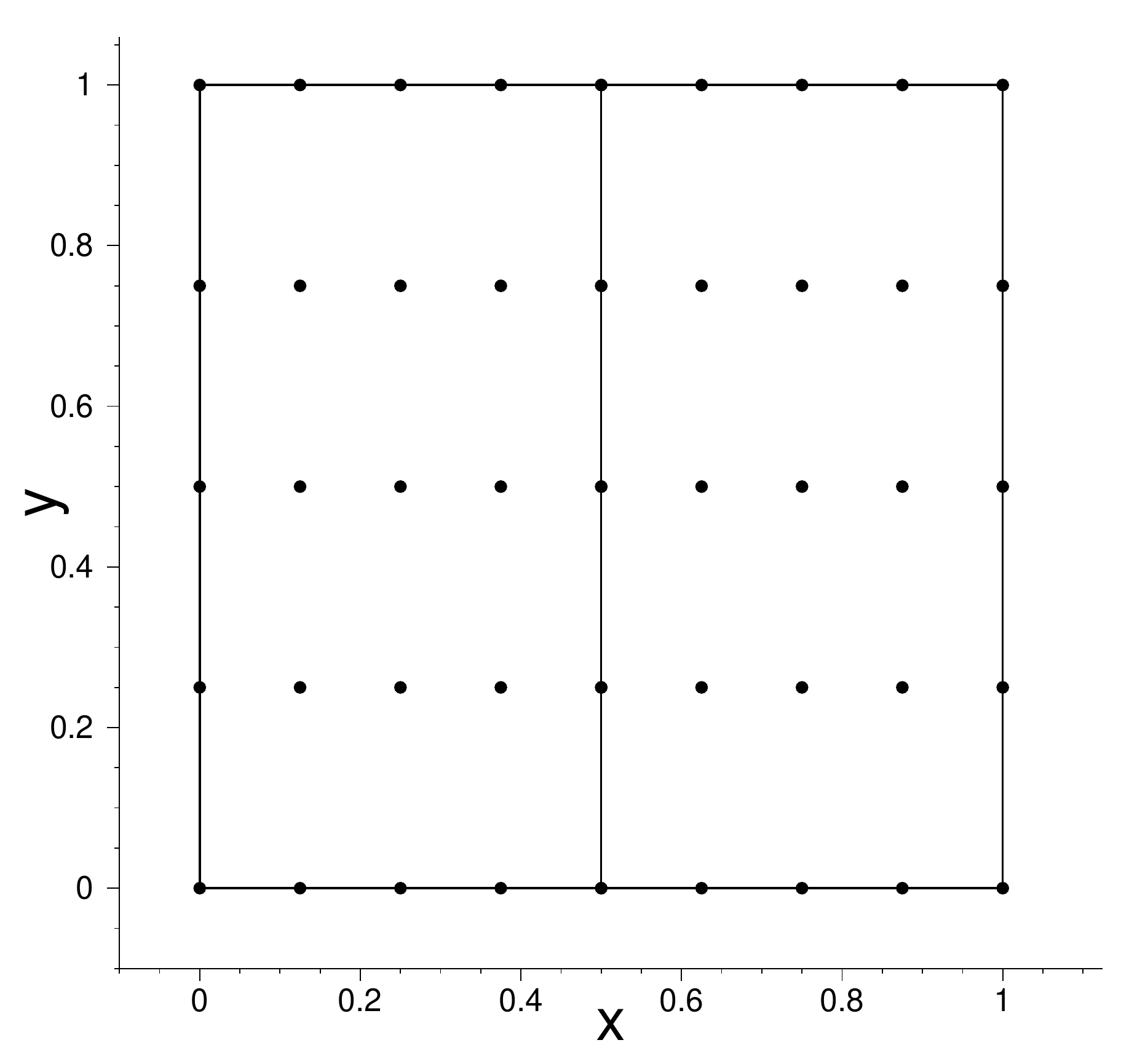}(b)
  }
  \caption{\small Sub-domains and collocation points:
    (a) Sketch of adjacent sub-domains. (b) Sketch of uniform grid
    points as collocation points ($5\times 5$ here) on
    two adjacent sub-domains.
  }
  \label{fg_2}
\end{figure}

We first outline a basic algorithm for computing
$(\bm\alpha,\bm\beta)$ by the nonlinear least squares (NLLSQ) method
with perturbations~\cite{DongL2021}.
It forms the basis for the variable projection algorithms
presented in the next subsection.

For the simplicity of presentation we focus on rectangular domains,
i.e.~$\Omega=[a_1,b_1]\times[a_2,b_2]\times\dots\times[a_d,b_d]$,
where $a_i$ and $b_i$ ($1\leqslant i\leqslant d$) denote the lower/upper
bounds of $\Omega$ in the $x_i$ direction, and
assume that $\Omega$ is partitioned into $N_i$ ($N_i\geqslant 1$)
sub-domains along $x_i$ ($1\leqslant i\leqslant d$).

To make the discussion  more concrete,
we specifically consider a second-order PDE  in two dimensions
($d=2$, $\mbs x=(x_1,x_2)=(x,y)$)
as an example in this and the next subsections. 
In the following discussions
we assume that equation~\eqref{eq_1a} is of second order with
respect to both $x$ and $y$, and we impose $C^1$ continuity conditions
across the sub-domain boundaries in both $x$ and $y$ directions.

Let the vectors $\bm{\mathcal{X}}=(X_0,X_1,\dots,X_{N_1})$ and
$\bm{\mathcal{Y}}=(Y_0,Y_1,\dots,Y_{N_2})$ denote the sub-domain boundary points
along the two directions, respectively, where $(X_0,X_{N_1})=(a_1,b_1)$ and
$(Y_0,Y_{N_2})=(a_2,b_2)$. The total number of sub-domains
is $N=N_1N_2$. We assume that the sub-domain $\Omega_e$ ($1\leqslant e\leqslant N$)
is characterized by the partition indices $(i,j)$ along the $x$
and $y$ directions (see Figure~\ref{fg_2}(a)), with the following relation,
\begin{equation}\label{eq_7}\small
  \begin{split}
    &
    \Omega_e=\Omega_{e(i,j)}=[X_{i-1},X_i]\times[Y_{j-1},Y_j], \quad
    e = e(i,j) = (i-1)N_2 + j, \ \ \text{for}\ 1\leqslant (i,j)\leqslant (N_1,N_2),
  \end{split}
\end{equation}
where ``$1\leqslant (i,j)\leqslant (N_1,N_2)$''
or ``$(1,1)\leqslant (i,j)\leqslant (N_1,N_2)$'' stands for
$1\leqslant i\leqslant N_1$ and $1\leqslant j\leqslant N_2$.
We will use this and similar notations hereafter for
conciseness.

With these settings the boundary conditions in~\eqref{eq_5b}
are reduced to,
\begin{subequations}\label{eq_08}\small
  \begin{align}
    &
    \mathcal{B}u_{e(1,j)}(a_1,y) = g(a_1,y), \quad
    \mathcal{B}u_{e(N_1,j)}(b_1,y) = g(b_1,y),
    \quad \text{for}\ 1\leqslant j\leqslant N_2;
    \\
    &
    \mathcal{B}u_{e(i,1)}(x,a_2) = g(x,a_2), \quad
    \mathcal{B}u_{e(i,N_2)}(x,b_2) = g(x,b_2),
    \quad \text{for}\ 1\leqslant i\leqslant N_1.
  \end{align}
\end{subequations}
Here $u_{e(i,j)}$ denotes $u$
on $\Omega_{e(i,j)}$, and
$e(i,j)$ is given by~\eqref{eq_7}.
The $C^1$ continuity conditions in~\eqref{eq_5d}
reduce to,
\begin{subequations}\label{eq_8}\footnotesize
  \begin{align}
    &
    u_{e(i,j)}(X_i,y) - u_{e(i+1,j)}(X_i,y) = 0,
    \quad \text{for}\ 1\leqslant (i,j)\leqslant (N_1-1,N_2);
    \label{eq_8a} \\
    & \left.\frac{\partial u_{e(i,j)}}{\partial x}\right|_{(X_i,y)}
    - \left.\frac{\partial u_{e(i+1,j)}}{\partial x}\right|_{(X_i,y)} = 0,
    \quad \text{for}\ 1\leqslant (i,j)\leqslant (N_1-1,N_2);
    \label{eq_8b} \\
    &
    u_{e(i,j)}(x,Y_j) - u_{e(i,j+1)}(x,Y_j) = 0,
    \quad \text{for}\ 1\leqslant (i,j)\leqslant (N_1,N_2-1);
    \label{eq_8c} \\
    & \left.\frac{\partial u_{e(i,j)}}{\partial y}\right|_{(x,Y_j)}
    - \left.\frac{\partial u_{e(i,j+1)}}{\partial y}\right|_{(x,Y_j)} = 0,
    \quad \text{for}\ 1\leqslant (i,j)\leqslant (N_1,N_2-1).
    \label{eq_8d}
  \end{align}
\end{subequations}
The equations~\eqref{eq_8a} and~\eqref{eq_8c} are the $C^0$
conditions on the horizontal/vertical sub-domain
boundaries, and the equations~\eqref{eq_8b} and~\eqref{eq_8d} are
the corresponding $C^1$ conditions.

The system to solve now consists of  equations~\eqref{eq_5a},
\eqref{eq_08}, \eqref{eq_5c}, and~\eqref{eq_8}.
This is a continuous system. We next enforce this system on
a set of collocation points and measurement points
to arrive at a discrete system
about the parameters $\bm\alpha$ and $\bm\beta$.


We choose a set of $Q$ ($Q\geqslant 1$) collocation points on
each sub-domain $\Omega_e$ ($1\leqslant e\leqslant N$),
denoted by $\mbs x_p^{e}=(x_p^e,y_p^e)$ ($1\leqslant p\leqslant Q$),
among which $Q_b$ ($1\leqslant Q_b<Q$) points reside on $\partial\Omega_e$.
Let $\mathbb{X}_e$ denote the set of collocation points
on $\Omega_e$, and $\mathbb{X}_e^b=\mathbb{X}_e\cap\partial\Omega_e$
denote the set of collocation points residing on
the sub-domain boundaries. The boundary collocation points on
adjacent sub-domains are required to be compatible.
That is, for any two adjacent sub-domains $(\Omega_{e_1},\Omega_{e_2})$,
those boundary collocation points from $\Omega_{e_1}$
that reside on the shared boundary $\partial\Omega_{e_1}\cap\partial\Omega_{e_2}$
are required to be identical to
those boundary collocation points from $\Omega_{e_2}$ that
reside on the same boundary.

The collocation points can in principle be chosen based on various
distributions (e.g.~random, uniform). In this paper
we focus on using uniform grid points as the collocation points;
see Figure~\ref{fg_2}(b) for an illustration with a $5\times 5$
uniform grid points as the collocation points on two neighboring
sub-domains. Let $Q_1$ and $Q_2$ denote the number of
uniform grid points along $x$ and $y$, with $Q=Q_1Q_2$.
The uniform collocation points on
the sub-domain $\Omega_e=\Omega_{e(m,l)}$ ($1\leqslant (m,l)\leqslant (N_1,N_2)$) are given by
\begin{equation}\label{eq_10}\footnotesize
  \left\{
  \begin{split}
    &
    \mbs x_p^{e} = \mbs x_{p(i,j)}^{e(m,l)}
    = \left(x_{p(i,j)}^{e(m,l)},y_{p(i,j)}^{e(m,l)} \right), \quad
  x_{p(i,j)}^{e(m,l)} = X_{m-1}+(i-1)(X_m-X_{m-1})/(Q_1-1),
  \\
  &
  y_{p(i,j)}^{e(m,l)} = Y_{l-1}+(j-1)(Y_l-Y_{l-1})/(Q_2-1),
  \quad \text{for}\ 1\leqslant (m,l,i,j)\leqslant (N_1,N_2,Q_1,Q_2);
  \\
  &
  p = p(i,j) = (i-1)Q_2 + j,
  \quad \text{for}\ 1\leqslant (p,i,j)\leqslant (Q,Q_1,Q_2).
  \end{split}
  \right.
\end{equation}

We assume that the measurement data 
is given on a set of $Q_s$ ($Q_s\geqslant 1$) random measurement points
(with a uniform distribution)
on each $\Omega_e$ ($1\leqslant e\leqslant N$), denoted
by $\bm\xi_p^e=(\xi_p^e,\eta_p^e)$ ($1\leqslant p\leqslant Q_s$).
We use $\mathbb{Y}_e$ to denote the set of measurement
points on $\Omega_e$ ($1\leqslant e\leqslant N$).

Once the hidden-layer coefficients of local NNs
are randomly assigned and the collocation and measurement points
are chosen, we compute the last hidden-layer field data
$\bm\Phi_e(\mbs x_p^e)$ and their derivatives (up to a certain order),
and the data for
$\mathcal{M}\bm\Phi_e(\bm\xi_p^e)$, by forward evaluations of the
neural network and by auto-differtiations.
We then store these data for subsequent use.
In light of~\eqref{eq_4},
for any given 
$\bm\beta=(\bm\beta_1^T,\dots,\bm\beta_N^T)^T$, we have
\begin{equation}\label{eq_011}\small
  \begin{split}
    &
  u_e(\mbs x_p^e) = \bm\Phi_e(\mbs x_p^e)\bm\beta_e, \
  \mathcal{D}u_e(\mbs x_p^e) = \mathcal{D}\bm\Phi_e(\mbs x_p^e)\bm\beta_e, \
  \mathcal{M}u_e(\bm\xi_q^e) = \mathcal{M}\bm\Phi_e(\bm\xi_q^e)\bm\beta_e, \
  1\leqslant (e,p,q)\leqslant (N,Q,Q_s),
  \end{split}
\end{equation}
where $\mathcal{D}$ is a linear differential operator and $\mathcal{M}$
is the measurement operator.

\begin{remark}\label{rem_02}
  To compute $\bm\Phi_e(\mbs x_p^e)$, $\mathcal{D}\bm\Phi_e(\mbs x_p^e)$
  and $\mathcal{M}\bm\Phi_e(\bm\xi_p^e)$, in the implementation we create
  a Keras sub-model, referred to as the last-hidden-layer-model,
  to the local NN for each sub-domain.
  The input nodes to
  this sub-model are identical to those of the original local NN,
  and the output nodes of this sub-model consist of those nodes in the last hidden
  layer of the original local NN.
  We compute $\bm\Phi_e(\mbs x_p^e)$ ($1\leqslant p\leqslant Q$) and
  $\bm\Phi_e(\bm\xi_p^e)$ ($1\leqslant p\leqslant Q_s$)
  by a forward evaluation of the last-hidden-layer-model for $\Omega_e$
  on the input data (collocation points, or measurement points).
  We compute the derivatives of $\bm\Phi_e$ on $\mbs x_p^e$ 
  or on $\bm\xi_p^e$ 
  by a forward-mode auto-differentiation of the last-hidden-layer-model,
  implemented by the ``ForwardAccumulator'' in the Tensorflow library.
  The forward-mode auto-differentiation is crucial to the performance
  of the ELM method (see~\cite{DongY2022rm}).
  
\end{remark}

To derive the discrete system we enforce~\eqref{eq_5a} on
all the collocation points in
$\mathbb{X}_e$ ($1\leqslant e\leqslant N$),
enforce~\eqref{eq_08} on all
the boundary collocation points in $\mathbb{X}_e^b\cap\partial\Omega$
for $1\leqslant e\leqslant N$,
enforce~\eqref{eq_5c} on all the measurement points
in $\mathbb{Y}_e$ ($1\leqslant e\leqslant N$),
and enforce~\eqref{eq_8} on
those collocation points from $\mathbb{X}_e^b$ ($1\leqslant e\leqslant N$)
that reside on the shared boundaries of adjacent sub-domains.

The discrete system corresponding to~\eqref{eq_5a} enforced on
the collocation points is,
\begin{equation}\label{eq_11}\footnotesize
  \alpha_1\mathcal{L}_1\left(u_{e}(\mbs x_p^{e})\right)
  +\dots+ \alpha_n\mathcal{L}_n\left(u_{e}(\mbs x_p^{e})\right)
  +\mathcal{F}\left(u_{e}(\mbs x_p^{e})\right)
  - f\left(\mbs x_p^{e}\right)=0,
  \ \text{for}\ \mbs x_p^e\in\mathbb{X}_e,
  \ 1\leqslant (e,p)\leqslant (N,Q).
\end{equation}
The discrete system corresponding to~\eqref{eq_08} on
the boundary collocation points is given by,
\begin{subequations}\label{eq_12}
  \begin{align}\footnotesize
    &
    \mathcal{B}u_{e(1,l)}(a_1,y_{p(1,j)}^{e(1,l)})
    -g(a_1,y_{p(1,j)}^{e(1,l)})=0,
    \quad \text{for}\ 1\leqslant (l,j)\leqslant (N_2,Q_2);
    \label{eq_12a} \\
    &
    \mathcal{B}u_{e(N_1,l)}(b_1,y_{p(Q_1,j)}^{e(N_1,l)})
    - g(b_1,y_{p(Q_1,j)}^{e(N_1,l)})=0,
    \quad \text{for}\ 1\leqslant (l,j)\leqslant (N_2,Q_2); \label{eq_12b}
    \\
    &
    \mathcal{B}u_{e(m,1)}(x_{p(i,1)}^{e(m,1)},a_2)
    - g(x_{p(i,1)}^{e(m,1)},a_2)=0,
    \quad \text{for}\ 1\leqslant (m,i)\leqslant (N_1,Q_1); \label{eq_12c}
    \\
    &
    \mathcal{B}u_{e(m,N_2)}(x_{p(i,Q_2)}^{e(m,N_2)},b_2)
    - g(x_{p(i,Q_2)}^{e(m,N_2)},b_2)=0,
    \quad \text{for}\ 1\leqslant (m,i)\leqslant (N_1,Q_1). \label{eq_12d}
  \end{align}
\end{subequations}
Here the functions $e(\cdot,\cdot)$ and $p(\cdot,\cdot)$ are
defined in~\eqref{eq_7} and~\eqref{eq_10}, respectively.
The discrete system corresponding to~\eqref{eq_5c} enforced on
the measurement points is given by
\begin{equation}\label{eq_13}\footnotesize
  \mathcal{M}u_e(\bm\xi_p^e) - S(\bm\xi_p^e)=0,
  \quad \text{for}\ \bm\xi_p^e\in\mathbb{Y}_e,
  \ 1\leqslant (e,p)\leqslant (N,Q_s).
\end{equation}
The discrete system corresponding to~\eqref{eq_8} enforced on
the interior sub-domain boundary points is,
\begin{subequations}\label{eq_14}
  \begin{align}\footnotesize
    &
    u_{e(m,l)}(X_m,y_{p(Q_1,j)}^{e(m,l)}) - u_{e(m+1,l)}(X_m,y_{p(1,j)}^{e(m+1,l)})=0,
    \quad \text{for}\ 1\leqslant (m,l,j)\leqslant (N_1-1,N_2,Q_2); \label{eq_14a}
    \\
    &
    \left.\frac{\partial u_{e(m,l)}}{\partial x}\right|_{(X_m,y_{p(Q_1,j)}^{e(m,l)})}
    - \left.\frac{\partial u_{e(m+1,l)}}{\partial x}\right|_{(X_m,y_{p(1,j)}^{e(m+1,l)})}=0,
    \quad \text{for}\ 1\leqslant (m,l,j)\leqslant (N_1-1,N_2,Q_2); \label{eq_14b}
    \\
    &
    u_{e(m,l)}(x_{p(i,Q_2)}^{e(m,l)},Y_l) - u_{e(m,l+1)}(x_{p(i,1)}^{e(m,l+1)},Y_l)=0,
    \quad \text{for}\ 1\leqslant (m,l,i)\leqslant (N_1,N_2-1,Q_1); \label{eq_14c}
    \\
    &
    \left.\frac{\partial u_{e(m,l)}}{\partial y}\right|_{(x_{p(i,Q_2)}^{e(m,l)},Y_l)}
    - \left.\frac{\partial u_{e(m,l+1)}}{\partial y}\right|_{(x_{p(i,1)}^{e(m,l+1)},Y_l)}=0,
    \quad \text{for}\ 1\leqslant (m,l,i)\leqslant (N_1,N_2-1,Q_1).
    \label{eq_14d}
  \end{align}
\end{subequations}
In the above equations $\mbs x_p^e$, $x_{p(i,j)}^{e(m,l)}$ and $y_{p(i,j)}^{e(m,l)}$
are defined in~\eqref{eq_10}, and $u_e(\mbs x)$ is given by~\eqref{eq_4}
and~\eqref{eq_011}.

The  equations~\eqref{eq_11}--\eqref{eq_14d}
form the system we would solve to determine
the inverse parameters
$\bm\alpha=(\alpha_1,\dots,\alpha_n)^T$ and
the trainable network parameters  $\bm\beta=(\beta_{11},\dots,\beta_{NM})^T$.
This is a system of nonlinear algebraic equations about ($\bm\alpha$, $\bm\beta$).
Note that the functions $\bm\Phi_e(\mbs x)$ ($1\leqslant e\leqslant N$)
and their derivatives evaluated on the collocation/measurement points,
which are involved in the operators such as $\mathcal{L}_i(u_e)$,
$\mathcal{F}(u_e)$, $\mathcal{B}u_e$, $\mathcal{M}u_e$, and $\mathcal{C}u_e$,
are computed by evaluations of the neural network and
auto-differentiations (see Remark~\ref{rem_02}).
This system consists of $N_c$ equations and a total of $N_a$ unknowns,
where
\begin{equation}\small
N_c = N(Q+Q_s+2Q_1+2Q_2), \quad
N_a = N_L+n = NM+n,
\end{equation}
and $N_L=NM$ is the total number of trainable parameters
 in the neural network.

We seek a least squares solution to this system,
and solve this system for $(\bm\alpha,\bm\beta)$ by the
nonlinear least squares (NLLSQ) method~\cite{Bjorck1996,DongL2021}.
In our implementation we take advantage of the quality
implementations of the nonlinear least squares method
in the scientific libraries, specficially the ``least\_squares()''
routine from the scipy.optimize package in Python for the current work.
This library routine implements the Gauss-Newton method~\cite{Bjorck1996}
together with a trust region algorithm~\cite{BranchCL1999,ByrdSS1988}.


Since the nonlinear least squares
method is a local optimization algorithm, it can be trapped to
a local-minimum solution that is unacceptable. It is therefore
crucial to combine the nonlinear least squares
method with some perturbation strategy when solving the nonlinear
least squares problem, in order to prevent the method from
being trapped to the worst local-minimum solutions.
In this paper we adopt the strategy for
the initial guess perturbation and sub-iteration procedure
developed in~\cite{DongL2021}, with some modifications, and combine it with
the nonlinear least squares method for solving the current system
arising from the inverse PDE problem.
We refer to the combined algorithm as the nonlinear least
squares method with perturbations (NLLSQ-perturb).
The NLLSQ-perturb algorithm is listed 
in the Appendix A of this paper (as Algorithm~\ref{alg_1}),
which contains explanations of
the various input parameters to the algorithm.

The NLLSQ-perturb algorithm (Algorithm~\ref{alg_1}) requires
two routines, one for computing the residual vector and
the other for computing the Jacobian matrix for an arbitrary given
approximation to the solution.
When the system~\eqref{eq_5} is enforced on the collocation points,
the residual function in~\eqref{eq_6} is reduced to the vector,
\begin{equation}\label{eq_16}\footnotesize
  \mbs R(\bm\alpha,\bm\beta) = \begin{bmatrix}
    \mbs R^{\text{pde}}(\bm\alpha,\bm\beta) \\
    \mbs R^{\text{bc}}(\bm\beta) \\
    \mbs R^{\text{mea}}(\bm\beta) \\
    \mbs R^{\text{ck}}(\bm\beta)
  \end{bmatrix}_{N_c\times 1}.
\end{equation}
In this expression,
\begin{equation}\label{eq_17}\footnotesize
  \left\{
  \begin{split}
    &
  \mbs R^{\text{pde}} = \begin{bmatrix}
    \vdots \\ R^{\text{pde}}_{ep} \\ \vdots
  \end{bmatrix}_{NQ\times 1}; \
  \mbs R^{\text{mea}} = \begin{bmatrix}
    \vdots \\ R^{\text{mea}}_{ep} \\ \vdots
  \end{bmatrix}_{NQ_s\times 1}; \
  \mbs R^{\text{bc}} = \begin{bmatrix}
    \mbs R^{\text{bc1}} \\ \mbs R^{\text{bc2}} \\ \mbs R^{\text{bc3}} \\ \mbs R^{\text{bc4}}
  \end{bmatrix}; \
  \mbs R^{\text{ck}} = \begin{bmatrix}
    \mbs R^{\text{ck1}} \\ \mbs R^{\text{ck2}} \\ \mbs R^{\text{ck3}} \\ \mbs R^{\text{ck4}}
  \end{bmatrix}; \
  \mbs R^{\text{bc1}}=\begin{bmatrix}\vdots\\ R_{lj}^{bc1}\\ \vdots  \end{bmatrix}_{N_2Q_2\times 1}, \\
  &
  \mbs R^{\text{bc2}}=\begin{bmatrix}\vdots\\ R_{lj}^{bc2}\\ \vdots  \end{bmatrix}_{N_2Q_2\times 1}, \
  \mbs R^{\text{bc3}}=\begin{bmatrix}\vdots\\ R_{mi}^{bc3}\\ \vdots  \end{bmatrix}_{N_1Q_1\times 1}, \
  \mbs R^{\text{bc4}}=\begin{bmatrix}\vdots\\ R_{mi}^{bc4}\\ \vdots  \end{bmatrix}_{N_1Q_1\times 1}; \
  \mbs R^{\text{ck1}}=\begin{bmatrix}\vdots\\ R_{mlj}^{ck1}\\ \vdots  \end{bmatrix}_{(N-N_2)Q_2\times 1},
  \\
  &
  \mbs R^{\text{ck2}}=\begin{bmatrix}\vdots\\ R_{mlj}^{ck2}\\ \vdots  \end{bmatrix}_{(N-N_2)Q_2\times 1},\
  \mbs R^{\text{ck3}}=\begin{bmatrix}\vdots\\ R_{mli}^{ck3}\\ \vdots  \end{bmatrix}_{(N-N_1)Q_1\times 1},\
  \mbs R^{\text{ck4}}=\begin{bmatrix}\vdots\\ R_{mli}^{ck4}\\ \vdots  \end{bmatrix}_{(N-N_1)Q_1\times 1}.
  \end{split}
  \right.
\end{equation}
In the above expressions, $R^{\text{pde}}_{ep}$ is the left hand side (LHS)
of~\eqref{eq_11},
and $R^{\text{mea}}_{ep}$ is the LHS of~\eqref{eq_13}.
$R_{lj}^{bc1}$, $R_{lj}^{bc2}$, $R_{mi}^{bc3}$ and $R_{mi}^{bc4}$ are the LHSs
of~\eqref{eq_12a}--\eqref{eq_12d}, respectively.
$R_{mlj}^{ck1}$, $R_{mlj}^{ck2}$, $R_{mli}^{ck3}$ and $R_{mli}^{ck4}$ are the LHSs
of~\eqref{eq_14a}--\eqref{eq_14d}, respectively.

We therefore compute the residual vector $\mbs R(\bm\alpha,\bm\beta)$ as follows.
Given arbitrary $(\bm\alpha,\bm\beta)$, we compute $u_e(\mbs x_p^e)$
($\mbs x_p^e\in\mathbb{X}_e$) for $1\leqslant e\leqslant N$,
and their derivatives by~\eqref{eq_011}.
Then we compute the LHSs of the equations~\eqref{eq_11},
\eqref{eq_12a}--\eqref{eq_12d}, \eqref{eq_13}, and~\eqref{eq_14a}--\eqref{eq_14d},
and assemble them to form the vectors
$\mbs R^{\text{pde}}$, $\mbs R^{\text{bc}}$, $\mbs R^{\text{mea}}$ and $\mbs R^{\text{ck}}$.
The residual vector $\mbs R(\bm\alpha,\bm\beta)$
is finally assembled according to~\eqref{eq_16}.
The procedure for computing $\mbs R(\bm\alpha,\bm\beta)$
is summarized in Algorithm~\ref{alg_2}.

\begin{algorithm}[tb]\small
  \DontPrintSemicolon
  \SetKwInOut{Input}{input}\SetKwInOut{Output}{output}

  \Input{ vector $\bm\theta=(\bm\alpha,\bm\beta)$;
    $\bm\Phi_e(x_p^e)$ and derivatives $(1\leqslant(e,p)\leqslant(N,Q))$;
    $\mathcal{M}\bm\Phi_e(\bm\xi_p^e)$ $(1\leqslant(e,p)\leqslant(N,Q_s))$.
  }
  \Output{residual vector $\mbs R(\bm\theta)$}
  \BlankLine
  \eIf{$\bm\theta = \bm\theta_s$}{
    retrieve $u_e(\mbs x_p^e)$ $(1\leqslant (e,p)\leqslant (N,Q))$ and their derivatives,
    and $\mathcal{M}u_e(\bm\xi_p^e)$  $(1\leqslant (e,p)\leqslant (N,Q_s))$\;
  }{
  compute $u_e(\mbs x_p^e)$ $(1\leqslant (e,p)\leqslant (N,Q))$,  
  their derivatives (up to a necessary order), and
  $\mathcal{M}u_e(\bm\xi_p^e)$ $(1\leqslant (e,p)\leqslant (N,Q_s))$  by~\eqref{eq_011}\;
  set $\bm\theta_s=\bm\theta$, and save $u_e(\mbs x_p^e)$ $(1\leqslant (e,p)\leqslant (N,Q))$,
  their derivatives,
    and $\mathcal{M}u_e(\bm\xi_p^e)$  $(1\leqslant (e,p)\leqslant (N,Q_s))$\;
  }
  \BlankLine
  compute $\mbs R^{\text{pde}}(\bm\theta)$, $\mbs R^{\text{bc}}(\bm\theta)$,
  $\mbs R^{\text{mea}}(\bm\theta)$, $\mbs R^{\text{ck}}(\bm\theta)$ by the LHSs
  of~\eqref{eq_11}--\eqref{eq_14d}, respectively\;
  form $\mbs R(\bm\theta)$ according to~\eqref{eq_16}\;
  
  \caption{Computing the residual $\mbs R(\bm\alpha,\bm\beta)$ for NLLSQ algorithm}
  \label{alg_2}
\end{algorithm}

\begin{remark}\label{rem_2}
  On line $4$ of Algorithm~\ref{alg_2}, the ``necessary order''
  refers to the order of all the derivative terms of $u_e$ involved in the
  system consisting of~\eqref{eq_11}--\eqref{eq_14d}.
  For example, if $\frac{\partial^2u_e}{\partial y^2}$ and $\frac{\partial u_e}{\partial x}$
  are involved in this system, one would need to compute these derivatives based on
  $\frac{\partial^2\bm\Phi_e}{\partial y^2}$ and $\frac{\partial\bm\Phi_e}{\partial x}$
  on line $4$ of this algorithm.

\end{remark}


The Jacobian matrix is given by
\begin{equation}\label{eq_18}\footnotesize
  \frac{\partial\mbs R}{\partial(\bm\alpha,\bm\beta)}=\begin{bmatrix}
  \frac{\partial\mbs R}{\partial\bm\alpha} &
  \frac{\partial\mbs R}{\partial\bm\beta}
  \end{bmatrix} =\begin{bmatrix}
  \frac{\partial\mbs R^{\text{pde}}}{\partial\bm\alpha} &
  \frac{\partial\mbs R^{\text{pde}}}{\partial\bm\beta} \\
  \mbs 0 & \frac{\partial\mbs R^{\text{bc}}}{\partial\bm\beta} \\
  \mbs 0 & \frac{\partial\mbs R^{\text{mea}}}{\partial\bm\beta} \\
  \mbs 0 & \frac{\partial\mbs R^{\text{ck}}}{\partial\bm\beta}
  \end{bmatrix}_{N_c\times N_L}.
\end{equation}
In this expression,
\begin{equation}\label{eq_19}\footnotesize
  \left\{
  \begin{split}
    &
  \frac{\partial\mbs R^{\text{pde}}}{\partial\bm\alpha}=\begin{bmatrix}
    \frac{\partial R_{ep}^{\text{pde}}}{\partial\alpha_i} 
    \end{bmatrix}_{NQ\times n} = \begin{bmatrix}
      \mathcal{L}_1(u_e(\mbs x_p^e)) & \dots & \mathcal{L}_n(u_e(\mbs x_p^e))
  \end{bmatrix}_{NQ\times n}, \
    \frac{\partial\mbs R^{\text{pde}}}{\partial\bm\beta}=\begin{bmatrix}
    \frac{\partial R_{ep}^{\text{pde}}}{\partial\beta_{ij}} 
    \end{bmatrix}_{NQ\times NM},
    \\
    &
    \frac{\partial\mbs R^{\text{mea}}}{\partial\bm\beta}=\begin{bmatrix}
    \frac{\partial R_{ep}^{\text{mea}}}{\partial\beta_{ij}} 
    \end{bmatrix}_{NQ_s\times NM},
    \frac{\partial\mbs R^{\text{bc}}}{\partial\bm\beta}=\begin{bmatrix}
    \frac{\partial\mbs R^{\text{bc1}} }{\partial\bm\beta} \\
    \frac{\partial\mbs R^{\text{bc2}} }{\partial\bm\beta} \\
    \frac{\partial\mbs R^{\text{bc3}} }{\partial\bm\beta} \\
    \frac{\partial\mbs R^{\text{bc4}} }{\partial\bm\beta} 
    \end{bmatrix}, \
    \frac{\partial\mbs R^{\text{ck}}}{\partial\bm\beta}=\begin{bmatrix}
    \frac{\partial\mbs R^{\text{ck1}} }{\partial\bm\beta} \\
    \frac{\partial\mbs R^{\text{ck2}} }{\partial\bm\beta} \\
    \frac{\partial\mbs R^{\text{ck3}} }{\partial\bm\beta} \\
    \frac{\partial\mbs R^{\text{ck4}} }{\partial\bm\beta} 
    \end{bmatrix}, \
    \frac{\partial\mbs R^{\text{bc1}} }{\partial\bm\beta} = \begin{bmatrix}
      \frac{\partial R_{lj}^{\text{bc1}} }{\partial\beta_{ik}} 
    \end{bmatrix}_{N_2Q_2\times NM},
    \\
    &
    \frac{\partial\mbs R^{\text{bc2}} }{\partial\bm\beta} = \begin{bmatrix}
      \frac{\partial R_{lj}^{\text{bc2}} }{\partial\beta_{ik}} 
    \end{bmatrix}_{N_2Q_2\times NM}, \
    \frac{\partial\mbs R^{\text{bc3}} }{\partial\bm\beta} = \begin{bmatrix}
      \frac{\partial R_{mi}^{\text{bc3}} }{\partial\beta_{lk}} 
    \end{bmatrix}_{N_1Q_1\times NM},
    \frac{\partial\mbs R^{\text{bc4}} }{\partial\bm\beta} = \begin{bmatrix}
      \frac{\partial R_{mi}^{\text{bc4}} }{\partial\beta_{lk}} 
    \end{bmatrix}_{N_1Q_1\times NM},
    \\
    &
    \frac{\partial\mbs R^{\text{ck1}} }{\partial\bm\beta} = \begin{bmatrix}
      \frac{\partial R_{mlj}^{\text{ck1}} }{\partial\beta_{iq}} 
    \end{bmatrix}_{(N-N_2)Q_2\times NM}, \
    \frac{\partial\mbs R^{\text{ck2}} }{\partial\bm\beta} = \begin{bmatrix}
      \frac{\partial R_{mlj}^{\text{ck2}} }{\partial\beta_{iq}} 
    \end{bmatrix}_{(N-N_2)Q_2\times NM}, \
    \frac{\partial\mbs R^{\text{ck3}} }{\partial\bm\beta} = \begin{bmatrix}
      \frac{\partial R_{mli}^{\text{ck3}} }{\partial\beta_{jq}} 
    \end{bmatrix}_{(N-N_1)Q_1\times NM},
    \\
    &
    \frac{\partial\mbs R^{\text{ck4}} }{\partial\bm\beta} = \begin{bmatrix}
      \frac{\partial R_{mli}^{\text{ck4}} }{\partial\beta_{jq}} 
    \end{bmatrix}_{(N-N_1)Q_1\times NM}.
  \end{split}
  \right.
\end{equation}
In the matrix $\frac{\partial\mbs R^{\text{pde}}}{\partial\bm\beta}$ the only
non-zero terms are
\begin{multline}\label{eq_20}\footnotesize
  \frac{\partial R_{ep}^{\text{pde}}}{\partial\beta_{ej}}
  = \alpha_1\mathcal{L}_1'(u_e(\mbs x_p^e))\phi_{ej}(\mbs x_p^e) + \dots
  + \alpha_n\mathcal{L}_n'(u_e(\mbs x_p^e))\phi_{ej}(\mbs x_p^e)
  + \mathcal{F}'(u_e(\mbs x_p^e))\phi_{ej}(\mbs x_p^e), \\
  \text{for}\ 1\leqslant (e,p,j)\leqslant (N,Q,M),
\end{multline}
where $\mathcal{L}_i'(u)$ ($1\leqslant i\leqslant n$) denote
the derivatives of $\mathcal{L}_i(u)$ with respect to $u$,
and $\mathcal{F}'(u)$ denotes the derivative of $\mathcal{F}(u)$
with respect to $u$.
In the matrix $\frac{\partial\mbs R^{\text{mea}}}{\partial\bm\beta}$ the only
non-zero terms are
\begin{equation}\label{eq_21}\footnotesize
  \frac{\partial R_{ep}^{\text{mea}}}{\partial\beta_{ej}} = \mathcal{M}\phi_{ej}(\bm\xi_p^e),
  \quad \text{for}\ 1\leqslant (e,p,j)\leqslant (N,Q_s,M).
\end{equation}
In the matrices $\frac{\partial\mbs R^{\text{bc1}}}{\partial\bm\beta}$,
$\frac{\partial\mbs R^{\text{bc2}}}{\partial\bm\beta}$,
$\frac{\partial\mbs R^{\text{bc3}}}{\partial\bm\beta}$ and
$\frac{\partial\mbs R^{\text{bc4}}}{\partial\bm\beta}$ the only non-zero terms are,
\begin{equation}\label{eq_22}
  \footnotesize
  \left\{
  \begin{split}
    &
    \frac{\partial R_{lj}^{\text{bc1}}}{\partial\beta_{lq}}
    = \mathcal{B}\phi_{eq}(a_1,y_{p}^{e}), \
    \text{where}\ e = e(1,l),\ p = p(1,j), \quad
    \text{for}\ 1\leqslant (l,j,q)\leqslant (N_2,Q_2,M);
    \\
    &
    \frac{\partial R_{lj}^{\text{bc2}}}{\partial\beta_{lq}}
    = \mathcal{B}\phi_{eq}(b_1,y_{p}^{e}), \
    \text{where}\ e = e(N_1,l),\ p = p(Q_1,j), \quad
    \text{for}\ 1\leqslant (l,j,q)\leqslant (N_2,Q_2,M);
    \\
    &
     \frac{\partial R_{mi}^{\text{bc3}}}{\partial\beta_{mq}}
    = \mathcal{B}\phi_{eq}(x_{p}^{e},a_2), \
    \text{where}\ e = e(m,1),\ p = p(i,1), \quad
    \text{for}\ 1\leqslant (m,i,q)\leqslant (N_1,Q_1,M);
    \\
    &
    \frac{\partial R_{mi}^{\text{bc4}}}{\partial\beta_{mq}}
    = \mathcal{B}\phi_{eq}(x_{p}^{e},b_2), \
    \text{where}\ e = e(m,N_2),\ p = p(i,Q_2), \quad
    \text{for}\ 1\leqslant (m,i,q)\leqslant (N_1,Q_1,M).
  \end{split}
  \right.
\end{equation}
In the matrices $\frac{\partial\mbs R^{\text{ck1}}}{\partial\bm\beta}$,
$\frac{\partial\mbs R^{\text{ck2}}}{\partial\bm\beta}$,
$\frac{\partial\mbs R^{\text{ck3}}}{\partial\bm\beta}$ and
$\frac{\partial\mbs R^{\text{ck4}}}{\partial\bm\beta}$ the only non-zero terms are,
\begin{equation}\label{eq_23}
  \footnotesize
  \left\{
  \begin{split}
    &\frac{\partial R_{mlj}^{\text{ck1}}}{\partial\beta_{e_1q}}
    = \phi_{e_1q}(X_m,y_{p_1}^{e_1}), \quad
    \frac{\partial R_{mlj}^{\text{ck1}}}{\partial\beta_{e_2q}}
    = -\phi_{e_2q}(X_m,y_{p_2}^{e_2}), \quad
    \text{where}\ e_1=e(m,l), \ e_2 = e(m+1,l), \\
    &\qquad\qquad\qquad\qquad
    p_1=p(Q_1,j), \ p_2=p(1,j), \
    \text{for}\ 1\leqslant (m,l,j,q)\leqslant (N_1-1,N_2,Q_2,M);
    \\
    &\frac{\partial R_{mlj}^{\text{ck2}}}{\partial\beta_{e_1q}}
    = \left.\frac{\partial\phi_{e_1q}}{\partial x}\right|_{(X_m,y_{p_1}^{e_1})}, \quad
    \frac{\partial R_{mlj}^{\text{ck2}}}{\partial\beta_{e_2q}}
    = -\left.\frac{\partial\phi_{e_2q}}{\partial x}\right|_{(X_m,y_{p_2}^{e_2})}, \
    \text{where}\ e_1=e(m,l), \ e_2 = e(m+1,l), \\
    &\qquad\qquad\qquad\qquad
    p_1=p(Q_1,j), \ p_2=p(1,j), \
    \text{for}\ 1\leqslant (m,l,j,q)\leqslant (N_1-1,N_2,Q_2,M);
    \\
    &\frac{\partial R_{mli}^{\text{ck3}}}{\partial\beta_{e_1q}}
    = \phi_{e_1q}(x_{p_1}^{e_1},Y_l), \quad
    \frac{\partial R_{mli}^{\text{ck3}}}{\partial\beta_{e_2q}}
    = -\phi_{e_2q}(x_{p_2}^{e_2},Y_l), \quad
    \text{where}\ e_1=e(m,l), \ e_2 = e(m,l+1), \\
    &\qquad\qquad\qquad\qquad
    p_1=p(i,Q_2), \ p_2=p(i,1), \
    \text{for}\ 1\leqslant (m,l,i,q)\leqslant (N_1,N_2-1,Q_1,M);
    \\
    &\frac{\partial R_{mli}^{\text{ck4}}}{\partial\beta_{e_1q}}
    = \left.\frac{\partial\phi_{e_1q}}{\partial y}\right|_{(x_{p_1}^{e_1},Y_l)}, \quad
    \frac{\partial R_{mli}^{\text{ck4}}}{\partial\beta_{e_2q}}
    = -\left.\frac{\partial\phi_{e_2q}}{\partial y}\right|_{(x_{p_2}^{e_2},Y_l)}, \
    \text{where}\ e_1=e(m,l), \ e_2 = e(m,l+1), \\
    &\qquad\qquad\qquad\qquad
    p_1=p(i,Q_2), \ p_2=p(i,1), \
    \text{for}\ 1\leqslant (m,l,i,q)\leqslant (N_1,N_2-1,Q_1,M).
  \end{split}
  \right.
\end{equation}

Therefore the Jacobian matrix can be computed as follows.
Given arbitrary $(\bm\alpha,\bm\beta)$, we compute
$u_e(\mbs x_p^e)$ ($1\leqslant (e,p)\leqslant(N,Q)$), their
derivatives, and $\mathcal{M}u_e(\bm\xi_p^e)$ ($1\leqslant (e,p)\leqslant(N,Q_s)$)
based on $\bm\beta$ and the pre-computed $\bm\Phi_e(\mbs x_p^e)$, their derivatives,
and the $\mathcal{M}\bm\Phi_e(\bm\xi_p^e)$ data.
Then we compute the Jacobian and related matrices by
the equations~\eqref{eq_18}--\eqref{eq_23}.
Algorithm~\ref{alg_3} summarizes the routine for computing
the Jacobian matrix.

\begin{algorithm}[tb]\small
  \DontPrintSemicolon
  \SetKwInOut{Input}{input}\SetKwInOut{Output}{output}

  \Input{ vector $\bm\theta=(\bm\alpha,\bm\beta)$;
    $\bm\Phi_e(x_p^e)$ and derivatives $(1\leqslant(e,p)\leqslant(N,Q))$;
    $\mathcal{M}\bm\Phi_e(\bm\xi_p^e)$ $(1\leqslant(e,p)\leqslant(N,Q_s))$.
  }
  \Output{Jacobian matrix $\frac{\partial\mbs R}{\partial\bm\theta}$}
  \BlankLine
  \eIf{$\bm\theta = \bm\theta_s$}{
    retrieve $u_e(\mbs x_p^e)$ $(1\leqslant (e,p)\leqslant (N,Q))$ and their derivatives,
    and $\mathcal{M}u_e(\bm\xi_p^e)$  $(1\leqslant (e,p)\leqslant (N,Q_s))$\;
  }{
  compute $u_e(\mbs x_p^e)$ $(1\leqslant (e,p)\leqslant (N,Q))$,  
  their derivatives (up to a necessary order), and
  $\mathcal{M}u_e(\bm\xi_p^e)$ $(1\leqslant (e,p)\leqslant (N,Q_s))$  by~\eqref{eq_011}\;
  set $\bm\theta_s=\bm\theta$, and save $u_e(\mbs x_p^e)$ $(1\leqslant (e,p)\leqslant (N,Q))$,
  their derivatives,
    and $\mathcal{M}u_e(\bm\xi_p^e)$  $(1\leqslant (e,p)\leqslant (N,Q_s))$\;
  }
  \BlankLine
  compute $\frac{\partial\mbs R^{\text{pde}}}{\partial\bm\alpha}$,
  $\frac{\partial\mbs R^{\text{pde}}}{\partial\bm\beta}$,
  $\frac{\partial\mbs R^{\text{bc}}}{\partial\bm\beta}$,
  $\frac{\partial\mbs R^{\text{mea}}}{\partial\bm\beta}$,
  $\frac{\partial\mbs R^{\text{ck}}}{\partial\bm\beta}$ by~\eqref{eq_19}--\eqref{eq_23}\;
  form $\frac{\partial\mbs R}{\partial\bm\theta}=\frac{\partial\mbs R}{\partial(\bm\alpha,\bm\beta)}$
  by~\eqref{eq_18}\;
  
  \caption{Computing the Jacobian matrix
    $\frac{\partial\mbs R}{\partial(\bm\alpha,\bm\beta)}$ for NLLSQ algorithm}
  \label{alg_3}
\end{algorithm}

\begin{remark}\label{rem_3}
  In Algorithms~\ref{alg_2} and~\ref{alg_3} we have stored the data for
  $u$, its derivatives, and $\mathcal{M}u$ on the collocation/measurement points
  corresponding to the $\bm\theta=(\bm\alpha,\bm\beta)$ value last computed
  (denoted by $\bm\theta_s$); see lines $1$ to $6$ in both algorithms.
  This saves computations, because in the nonlinear least squares iterations
  Algorithm~\ref{alg_2} is typically invoked first to compute the residual
  corresponding to some $(\bm\alpha,\bm\beta)$, and then Algorithm~\ref{alg_3}
  is invoked to compute the Jacobian for the same $(\bm\alpha,\bm\beta)$.

\end{remark}


\begin{remark}\label{rem_4}

    In this work the hidden-layer coefficients are assigned to uniform random values
    generated on the interval $[-R_m,R_m]$, where $R_m>0$ is a constant.
    The $R_m$ value influences the accuracy of the simulation results of inverse PDE
    problems, similar
    to what has been observed in forward problems (see~\cite{DongL2021,DongY2022rm}).
    In this paper we compute a near-optimal $R_m$
    using the method from~\cite{DongY2022rm} based on the differential
    evolution algorithm, and employ this value (or a value nearby)
    in numerical simulations of inverse PDEs.

\end{remark}

\begin{remark}\label{rem_04}
  For noisy measurement data $S(\bm\xi)$, we observe that
  scaling the residual vector associated with the measurement ($\mbs R^{\text{mea}}$) by
  a constant factor can improve the accuracy of the
  results (more robust to noise). Let $\lambda_{mea}>0$ denote a prescribed
  constant. We scale the equation~\eqref{eq_13} by $\lambda_{mea}$,
  \begin{equation}\label{eq_a25}\footnotesize
    \lambda_{mea}\mathcal{M}u_e(\bm\xi_p^e) - \lambda_{mea}S(\bm\xi_p^e)=0,
    \quad \text{for}\ \bm\xi_p^e\in\mathbb{Y}_e,
    \ 1\leqslant (e,p)\leqslant (N,Q_s).
  \end{equation}
  Then in the presented method we replace equation~\eqref{eq_13} by the scaled equation~\eqref{eq_a25},
  with corresponding changes to the computation of the residual vector and
  the Jacobian matrix. The scaling factor $\lambda_{mea}$ will cause some change to the
  least squares solution to ($\bm\alpha,\bm\beta$). 
  When the data $S(\bm\xi)$ is noisy, numerical experiments indicate that
  employing a constant $0<\lambda_{mea}<1$ can in general improve the accuracy of
  the computed $\bm\alpha$ and $u(\mbs x)$ markedly, compared with the case without
  scaling (i.e.~$\lambda_{mea}=1$). Note that employing the scaled equation~\eqref{eq_a25}
  is equivalent to using a scaled term $\frac12\lambda_{mea}^2\|\mbs R^{\text{mea}}\|^2$
  in the underlying loss function for the nonlinear least squares method.
  
\end{remark}

\begin{remark}\label{rem_05}

    The method developed here can be applied to inverse PDEs in which the inverse parameters
    may be an unknown field distribution. Consider for example,
    \begin{equation}\label{eq_025}\footnotesize
      \gamma(\mbs x)\mathcal{L}(u) + \mathcal{F}(u) = f(\mbs x),
    \end{equation}
    where the coefficient $\gamma(\mbs x)$ is an unknown field and $u(\mbs x)$
    is the unknown solution to the forward problem.
    In this case we can expand $\gamma(\mbs x)$ in terms of a set of basis functions
    and transform \eqref{eq_025} into a form similar to~\eqref{eq_1a}, in which
    the expansion coefficients of $\gamma(\mbs x)$ become the inverse parameters.
    Therefore the inverse problem can be computed using the method presented above.
    In this work we employ the same bases in the expansion for $u(\mbs x)$ (see~\eqref{eq_4})
    and for $\gamma(\mbs x)$. This translates into two nodes in the output layer
    of the neural network architecture, one representing $u(\mbs x)$ and the other representing
    $\gamma(\mbs x)$. When more inverse coefficient fields are involved, one can correspondingly
    increase the number of nodes in the output layer of the neural network.
    We will present a numerical example for an inverse PDE similar to~\eqref{eq_025}
    in Section~\ref{sec:tests}.

\end{remark}

\subsection{Variable Projection Algorithms for Network Training}
\label{sec:varpro}

This subsection outlines two algorithms for computing $(\bm\alpha,\bm\beta)$,
both based on the variable projection (VarPro)
idea~\cite{GolubP1973,GolubP2003,DongY2022}
but with different formulations.
In the first formulation (VarPro-F1), the inverse parameters ($\bm\alpha$) are
eliminated from the problem to attain a reduced problem about
$\bm\beta$ only. The reduced problem is solved by
the nonlinear least squares method first for $\bm\beta$,
and then $\bm\alpha$ is computed by
the linear least squares method.
In the second formulation (VarPro-F2),
the field solution (equivalently, the $\bm\beta$ parameters)
is eliminated from the problem to attain a reduced problem about
$\bm\alpha$ only. The reduced problem
is solved first by the nonlinear least squares method for $\bm\alpha$, and then
 $\bm\beta$ is computed based on the $\bm\alpha$ already obtained.
The problem settings and notations here
follow those of Section~\ref{sec:nllsq}.

\subsubsection{Formulation \#1 (VarPro-F1): Eliminating the Inverse Parameters}
\label{sec:vp_1}

We start with the discrete system consisting of
equations~\eqref{eq_11}--\eqref{eq_14d}.
We  re-arrange this system symbolically into a matrix equation
about the parameters $\bm\alpha=(\alpha_1,\dots,\alpha_n)^T$,
\begin{equation}\label{eq_26}\small
  \mbs H(\bm\beta)\bm\alpha = \mbs b(\bm\beta),
\end{equation}
where
\begin{equation}\label{eq_27}
  \footnotesize
  \left\{
  \begin{split}
    &
  \mbs H(\bm\beta)=\begin{bmatrix}
  \mbs H^{\text{pde}}(\bm\beta) \\ \mbs 0 \\ \mbs 0 \\ \mbs 0
  \end{bmatrix}_{N_c\times n},
  \mbs b(\bm\beta)=\begin{bmatrix}
  \mbs b^{\text{pde}}(\bm\beta) \\
  -\mbs R^{\text{bc}}(\bm\beta) \\
  -\mbs R^{\text{mea}}(\bm\beta) \\
  -\mbs R^{\text{ck}}(\bm\beta)
  \end{bmatrix}_{N_c\times 1}, 
  \mbs H^{\text{pde}}(\bm\beta)=\begin{bmatrix}
  \vdots &  & \vdots \\
  \mathcal{L}_1\left(u_e(\mbs x_p^e)\right) & \cdots & \mathcal{L}_n\left(u_e(\mbs x_p^e) \right)\\
  \vdots &  & \vdots
  \end{bmatrix}_{NQ\times n}, \\
  &
  \mbs b^{\text{pde}}(\bm\beta)=\begin{bmatrix}b_{ep}^{\text{pde}}  \end{bmatrix}_{NQ\times 1}
  = \begin{bmatrix}
  \vdots \\
  f\left(\mbs x_p^e \right) - \mathcal{F}\left(u_e(\mbs x_p^e) \right) \\
  \vdots 
  \end{bmatrix}_{NQ\times 1}.
  \end{split}
  \right.
\end{equation}
In these expressions, $\mbs R^{\text{bc}}$, $\mbs R^{\text{mea}}$
and $\mbs R^{\text{ck}}$ are defined in~\eqref{eq_17}.

For any given $\bm\beta$, the least squares solution to~\eqref{eq_26}
with the minimum norm is given by
\begin{equation}\label{eq_28}\small
  \bm\alpha = \mbs H^+(\bm\beta)\mbs b(\bm\beta),
\end{equation}
where $\mbs H^+(\bm\beta)$ denotes the Moore-Penrose inverse
of $\mbs H(\bm\beta)$.
Substituting this expression into~\eqref{eq_26} gives rise to
a reduced system about $\bm\beta$ only.
The residual of this reduced system  (see also~\eqref{eq_16})
is given by
\begin{equation}\label{eq_29}\small
  \mbs r(\bm\beta) = \mbs R(\bm\alpha,\bm\beta)
  = \mbs H(\bm\beta)\bm\alpha - \mbs b(\bm\beta)
  = \mbs H(\bm\beta)\mbs H^+(\bm\beta)\mbs b(\bm\beta) - \mbs b(\bm\beta).
\end{equation}
We determine the optimum $\bm\beta^*$ by minimizing the Euclidean norm of this residual,
\begin{equation}\label{eq_30}\small
  \bm\beta^* = \argmin_{\bm\beta}\frac12\|\mbs r(\bm\beta) \|^2
  = \argmin_{\bm\beta}\frac12\|\mbs H(\bm\beta)\mbs H^+(\bm\beta)\mbs b(\bm\beta)
  - \mbs b(\bm\beta)  \|^2
\end{equation}
where $\|\cdot\|$ denotes the Euclidean norm.
With $\bm\beta$ determined by~\eqref{eq_30}, we solve the system~\eqref{eq_26} for $\bm\alpha$
by the linear least squares method with the minimum-norm solution (or by directly
using~\eqref{eq_28}).

Equation~\eqref{eq_30} represents a nonlinear least squares problem about $\bm\beta$.
We solve this problem by the NLLSQ-perturb algorithm (Algorithm~\ref{alg_1} in
Appendix A).
As noted previously, two routines are required for this algorithm, one for computing
the reduced residual $\mbs r(\bm\beta)$ and the other
for computing the Jacobian matrix of the reduced problem,
$\frac{\partial\mbs r}{\partial\bm\beta}$, for any given $\bm\beta$.

\begin{algorithm}[tb]\small
  \DontPrintSemicolon
  \SetKwInOut{Input}{input}\SetKwInOut{Output}{output}

  \Input{ $\bm\beta$;
    $\bm\Phi_e(x_p^e)$ and derivatives $(1\leqslant(e,p)\leqslant(N,Q))$;
    $\mathcal{M}\bm\Phi_e(\bm\xi_p^e)$ $(1\leqslant(e,p)\leqslant(N,Q_s))$.
  }
  \Output{reduced residual $\mbs r(\bm\beta)$}
  \BlankLine
  \eIf{$\bm\beta = \bm\beta_s$}{
    retrieve $\mbs H(\bm\beta_s)$, $\mbs b(\bm\beta_s)$, $\bm\alpha^{LS}$\;
    set $\mbs H(\bm\beta)= \mbs H(\bm\beta_s)$ and $\mbs b(\bm\beta)=\mbs b(\bm\beta_s)$\;
  }{
  compute $u_e(\mbs x_p^e)$ $(1\leqslant (e,p)\leqslant (N,Q))$,  
  their derivatives (up to a necessary order), and
  $\mathcal{M}u_e(\bm\xi_p^e)$ $(1\leqslant (e,p)\leqslant (N,Q_s))$  by~\eqref{eq_011}\;
  compute $\mbs H(\bm\beta)$ and $\mbs b(\bm\beta)$ by~\eqref{eq_27} and~\eqref{eq_17}\;
  solve equation~\eqref{eq_26} for $\bm\alpha$ by the linear least squares method,
  and let $\bm\alpha^{LS}=\bm\alpha$\;
  set $\bm\beta_s=\bm\beta$, and save $\mbs H(\bm\beta)$, $\mbs b(\bm\beta)$, $\bm\alpha^{LS}$\;
  }
  \BlankLine
  compute $\mbs r(\bm\beta)$ by equation~\eqref{eq_31}\;
  
  \caption{Computing reduced residual $\mbs r(\bm\beta)$
    for VarPro-F1.
  }
  \label{alg_4}
\end{algorithm}

We compute the reduced residual as follows. For any given $\bm\beta$,
we solve equation~\eqref{eq_26} for $\bm\alpha$ (with minimum norm) by
the linear least squares method. Let $\bm\alpha^{LS}$ denote this solution.
Then the residual is given by
\begin{equation}\label{eq_31}
  \mbs r(\bm\beta) = \mbs H(\bm\beta)\bm\alpha^{LS} - \mbs b(\bm\beta).
\end{equation}
Algorithm~\ref{alg_4} summarizes the procedure for
computing the reduced residual.


To compute the Jacobian of the reduced residual, we note
the following formula owing to~\cite{GolubP1973},
\begin{equation}\label{eq_32}\footnotesize
  \begin{split}
  \frac{\partial}{\partial\bm\theta}\left[\mbs H(\bm\theta)\mbs H^+(\bm\theta) \right]
  &= \left[\mbs I - \mbs H(\bm\theta)\mbs H^+(\bm\theta) \right]
  \frac{\partial\mbs H}{\partial\bm\theta}\mbs H^+(\bm\theta)
  + \left[\mbs H^T(\bm\theta) \right]^+\frac{\partial\mbs H^T}{\partial\bm\theta}
  \left[\mbs I - \mbs H(\bm\theta)\mbs H^+(\bm\theta) \right] \\
  &\approx \left[\mbs I - \mbs H(\bm\theta)\mbs H^+(\bm\theta) \right]
  \frac{\partial\mbs H}{\partial\bm\theta}\mbs H^+(\bm\theta),
  \end{split}
\end{equation}
where $\mbs I$ is the identity matrix and
on the second line we have kept only the first term in the formula
as an approximation to the LHS, thanks to the suggestion of~\cite{Kaufman1975}.
In light of~\eqref{eq_29} and~\eqref{eq_32}, we have
\begin{equation}\label{eq_33}\footnotesize
  \begin{split}
    \frac{\partial\mbs r}{\partial\bm\beta} &=
    \left(\frac{\partial}{\partial\bm\beta}\left[\mbs H(\bm\beta)\mbs H^+(\bm\beta) \right]\right)
    \mbs b(\bm\beta)
    + \left[\mbs H(\bm\beta)\mbs H^+(\bm\beta) \right]\frac{\partial\mbs b}{\partial\bm\beta}
    - \frac{\partial\mbs b}{\partial\bm\beta} \\
    &\approx \frac{\partial\mbs H}{\partial\bm\beta}\mbs H^+(\bm\beta)\mbs b(\bm\beta)
    - \mbs H(\bm\beta)\mbs H^+(\bm\beta)
    \frac{\partial\mbs H}{\partial\bm\beta}\mbs H^+(\bm\beta)\mbs b(\bm\beta)
    + \mbs H(\bm\beta)\mbs H^+(\bm\beta) \frac{\partial\mbs b}{\partial\bm\beta}
    - \frac{\partial\mbs b}{\partial\bm\beta} \\
    &= \left(\frac{\partial\mbs H}{\partial\bm\beta}\mbs H^+(\bm\beta)\mbs b(\bm\beta)
    - \frac{\partial\mbs b}{\partial\bm\beta}\right)
    - \mbs H(\bm\beta)\mbs H^+(\bm\beta)
    \left(\frac{\partial\mbs H}{\partial\bm\beta}\mbs H^+(\bm\beta)\mbs b(\bm\beta)
    - \frac{\partial\mbs b}{\partial\bm\beta}\right) \\
    &= \mbs J_1(\bm\beta) - \mbs J_2(\bm\beta),
  \end{split}
\end{equation}
where
\begin{equation}\label{eq_34}\footnotesize
  \begin{split}
    &
    \mbs J_1(\bm\beta) = \mbs J_0(\bm\beta) - \frac{\partial\mbs b}{\partial\bm\beta},\quad 
    \mbs J_2(\bm\beta) = \mbs H(\bm\beta)\mbs H^+(\bm\beta)\mbs J_1(\bm\beta), \quad
    \mbs J_0(\bm\beta) = \frac{\partial\mbs H}{\partial\bm\beta}\mbs H^+(\bm\beta)\mbs b(\bm\beta).
  \end{split}
\end{equation}

Therefore, we need a procedure for computing
$\mbs J_0(\bm\beta)$, $\frac{\partial\mbs b}{\partial\bm\beta}$ and $\mbs J_2(\bm\beta)$.
$\mbs J_0(\bm\beta)$ can be computed as follows,
\begin{equation}\label{eq_35}\footnotesize
  \begin{split}
  \mbs J_0(\bm\beta) &= \frac{\partial\mbs H}{\partial\bm\beta}\mbs H^+(\bm\beta)\mbs b(\bm\beta)
  = \frac{\partial\mbs H}{\partial\bm\beta}\bm\alpha^{LS}
  = \frac{\partial\left[\mbs H(\bm\beta)\bm\alpha^{LS} \right]}{\partial\bm\beta}
  =\begin{bmatrix}
  \frac{\partial\mbs R^{\text{pdeI}}}{\partial\bm\beta} \\
  \mbs 0 \\ \mbs 0 \\ \mbs 0
  \end{bmatrix}_{N_c\times NM}.
  \end{split}
\end{equation}
In this equation, $\bm\alpha^{LS}=(\alpha_1^{LS},\dots,\alpha_n^{LS})^T$ is the minimum-norm
solution to~\eqref{eq_26} computed by the linear least squares method, and
\begin{subequations}\label{eq_36}
\begin{align}\footnotesize
  &\mbs R^{\text{pdeI}}(\bm\beta)=\begin{bmatrix}R_{ep}^{\text{pdeI}}  \end{bmatrix}_{NQ\times 1}
  =\begin{bmatrix}
  \vdots \\
  \alpha_1^{LS}\mathcal{L}_1(u_e(\mbs x_p^e)) + \dots +
  \alpha_n^{LS}\mathcal{L}_n(u_e(\mbs x_p^e)) \\
  \vdots
  \end{bmatrix}_{NQ\times 1}, \\
  &
  \frac{\partial\mbs R^{\text{pdeI}}}{\partial\bm\beta}=\begin{bmatrix}
  \frac{\partial R_{ep}^{\text{pdeI}}}{\partial\beta_{ij}} 
  \end{bmatrix}_{NQ\times NM}.
\end{align}
\end{subequations}
In the matrix $\frac{\partial\mbs R^{\text{pdeI}}}{\partial\bm\beta}$ the only non-zero
terms are,
\begin{equation}\label{eq_37}\small
  \frac{\partial R_{ep}^{\text{pdeI}}}{\partial\beta_{ej}}
  = \alpha_1^{LS}\mathcal{L}_1'(u_e(\mbs x_p^e))\phi_{ej}(\mbs x_p^e) + \dots
  + \alpha_n^{LS}\mathcal{L}_n'(u_e(\mbs x_p^e))\phi_{ej}(\mbs x_p^e),
  \quad \text{for}\ 1\leqslant (e,p,j) \leqslant (N,Q,M).
\end{equation}
It is important to note that, when computing $\frac{\partial\mbs R^{\text{pdeI}}}{\partial\bm\beta}$,
we treat $\bm\alpha^{LS}$ as a constant vector independent of $\bm\beta$.

$\frac{\partial\mbs b}{\partial\bm\beta}$ is computed as follows,
\begin{equation}\label{eq_38}\footnotesize
  \frac{\partial\mbs b}{\partial\bm\beta} = \begin{bmatrix}
    \frac{\partial\mbs b^{\text{pde}}}{\partial\bm\beta} \\
    - \frac{\partial\mbs R^{\text{bc}}}{\partial\bm\beta} \\
    - \frac{\partial\mbs R^{\text{mea}}}{\partial\bm\beta} \\
    - \frac{\partial\mbs R^{\text{ck}}}{\partial\bm\beta} 
  \end{bmatrix},
  \quad
  \frac{\partial\mbs b^{\text{pde}}}{\partial\bm\beta} = \begin{bmatrix}
    \frac{\partial b_{ep}^{\text{pde}}}{\partial\beta_{ij}} 
  \end{bmatrix}_{NQ\times NM},
\end{equation}
where $\frac{\partial\mbs R^{\text{bc}}}{\partial\bm\beta}$,
$\frac{\partial\mbs R^{\text{mea}}}{\partial\bm\beta}$ and
$\frac{\partial\mbs R^{\text{ck}}}{\partial\bm\beta}$ are given
in~\eqref{eq_19} and \eqref{eq_21}--\eqref{eq_23}.
The only non-zero terms in  $\frac{\partial\mbs b^{\text{pde}}}{\partial\bm\beta}$
are,
\begin{equation}\label{eq_39}\small
  \frac{\partial b_{ep}^{\text{pde}}}{\partial\beta_{ej}}
  = -\mathcal{F}'\left(u_e(\mbs x_p^e) \right)\phi_{ej}(\mbs x_p^e),
  \quad \text{for}\ 1\leqslant(e,p,j)\leqslant (N,Q,M).
\end{equation}
With $J_0(\bm\beta)$ and $\frac{\partial\mbs b}{\partial\bm\beta}$ determined,
we can compute $\mbs J_1(\bm\beta)$ by~\eqref{eq_34}.

\begin{algorithm}[tb]\small
  \DontPrintSemicolon
  \SetKwInOut{Input}{input}\SetKwInOut{Output}{output}

  \Input{ $\bm\beta$;
    $\bm\Phi_e(x_p^e)$ and derivatives $(1\leqslant(e,p)\leqslant(N,Q))$;
    $\mathcal{M}\bm\Phi_e(\bm\xi_p^e)$ $(1\leqslant(e,p)\leqslant(N,Q_s))$.
  }
  \Output{Jacobian matrix $\frac{\partial\mbs r}{\partial\bm\beta}$ }
  \BlankLine
  \eIf{$\bm\beta = \bm\beta_s$}{
    retrieve $\mbs H(\bm\beta_s)$, $\mbs b(\bm\beta_s)$, $\bm\alpha^{LS}$\;
    set $\mbs H(\bm\beta)= \mbs H(\bm\beta_s)$ and $\mbs b(\bm\beta)=\mbs b(\bm\beta_s)$\;
  }{
  compute $u_e(\mbs x_p^e)$ $(1\leqslant (e,p)\leqslant (N,Q))$,  
  their derivatives (up to a necessary order), and
  $\mathcal{M}u_e(\bm\xi_p^e)$ $(1\leqslant (e,p)\leqslant (N,Q_s))$  by~\eqref{eq_011}\;
  compute $\mbs H(\bm\beta)$ and $\mbs b(\bm\beta)$ by~\eqref{eq_27} and~\eqref{eq_17}\;
  solve equation~\eqref{eq_26} for $\bm\alpha$ by the linear least squares method,
  and let $\bm\alpha^{LS}=\bm\alpha$\;
  set $\bm\beta_s=\bm\beta$, and save $\mbs H(\bm\beta)$, $\mbs b(\bm\beta)$, $\bm\alpha^{LS}$\;
  }
  \BlankLine
  compute $\mbs J_0(\bm\beta)$ by equations~\eqref{eq_35}--\eqref{eq_37}\;
  compute $\frac{\partial\mbs b}{\partial\bm\beta}$ by~\eqref{eq_38}, \eqref{eq_39},
  \eqref{eq_19}, and~\eqref{eq_21}--\eqref{eq_23}\;
  compute $\mbs J_1(\bm\beta)$ by~\eqref{eq_34}\;
  compute $\mbs J_2(\bm\beta)$ by~\eqref{eq_40}--\eqref{eq_41}\;
  compute $\frac{\partial\mbs r}{\partial\bm\beta}$ by~\eqref{eq_33}\;
  \caption{Computing Jacobian matrix $\frac{\partial\mbs r}{\partial\bm\beta}$
    for VarPro-F1.
  }
  \label{alg_5}
\end{algorithm}

In light of~\eqref{eq_34}, we compute $\mbs J_2(\bm\beta)$ by the following equations,
\begin{align}\small
  &
  \mbs H(\bm\beta)\mbs K = \mbs J_1(\bm\beta), \label{eq_40} \\
  &
  \mbs J_2(\bm\beta) = \mbs H(\bm\beta)\mbs K. \label{eq_41}
\end{align}
We first solve equation~\eqref{eq_40} for the $n\times NM$ matrix $\mbs K$
by the linear least squares method, and then compute $\mbs J_2(\bm\beta)$
by equation~\eqref{eq_41} with a matrix multiplication.

Therefore, given an arbitrary $\bm\beta$, we compute $\mbs J_0(\bm\beta)$
by~\eqref{eq_35}--\eqref{eq_37}, $\frac{\partial\mbs b}{\partial\bm\beta}$
by~\eqref{eq_38} and~\eqref{eq_39}, and $\mbs J_1(\bm\beta)$ by~\eqref{eq_34}.
Then we compute $\mbs J_2(\bm\beta)$ by~\eqref{eq_40}--\eqref{eq_41}.
The (approximate) Jacobian matrix of the reduced problem
is then given by~\eqref{eq_33}.
The procedure for computing the Jacobian matrix is summarized in
the Algorithm~\ref{alg_5}.

The overall VarPro-F1 algorithm for solving the inverse problem consists of
two steps: (i) Invoke the
NLLSQ-perturb algorithm (Algorithm~\ref{alg_1} in Appendix A) to compute $\bm\beta$
from the reduced problem~\eqref{eq_30}, with the routines
given in Algorithms~\ref{alg_4} and~\ref{alg_5} as input.
(ii) Solve~\eqref{eq_26} for $\bm\alpha$ by the linear least squares method.

\begin{remark}
  In the VarPro-F1 algorithm, one only needs to solve
  linear systems by the linear least squares
  method. The Moore-Penrose inverse of the coefficient matrix is not explicitly computed.
  In our implementation we employ the linear least squares routine
  scipy.linalg.lstsq() from the scipy package
  in Python, which in turn uses the linear least squares implementation
  in the LAPACK library.
\end{remark}

\subsubsection{Formulation \#2 (VarPro-F2): Eliminating the Field Function}
\label{sec:vp_2}

We next present an alternative formulation (VarPro-F2) of variable projection,
which is reciprocal to the VarPro-F1 algorithm of Section \ref{sec:vp_1}.
In this formulation,
we eliminate the field function $u$ (or the parameters $\bm\beta$)
from the problem to attain a reduced problem about  $\bm\alpha$ only.
We then solve the reduced problem first for $\bm\alpha$, and compute
the parameters $\bm\beta$ afterwards.

This formulation applies to cases in which
the operators $\mathcal{L}_i$ ($1\leqslant i\leqslant n$)
and $\mathcal{F}$ are all linear with respect to $u$. We first present the
algorithm with regard to this case below. 
Then we outline an extension in a remark (Remark~\ref{rem_6})
by combining this algorithm with a Newton iteration
to deal with cases in which these operators
are nonlinear with respect to $u$.

Let us now assume that $\mathcal{L}_i$ ($1\leqslant i\leqslant n$)
and $\mathcal{F}$ are all linear operators, and we again start with
the discrete system consisting of the equations~\eqref{eq_11}--\eqref{eq_14d}.
We re-arrange this system into a matrix equation about the
trainable network parameters 
$\bm\beta=(\bm\beta_1^T,\dots,\bm\beta_N^T)^T=(\beta_{11},\dots,\beta_{NM})^T$,
\begin{equation}\label{eq_42}\small
  \mbs H(\bm\alpha)\bm\beta = \mbs b,
\end{equation}
where
\begin{equation}\label{eq_43}\footnotesize
  \begin{split}
    &
    \mbs H(\bm\alpha)=\begin{bmatrix}
    \mbs H^{\text{pde}}(\bm\alpha) \\ \mbs H^{\text{bc}} \\ \mbs H^{\text{mea}} \\ \mbs H^{\text{ck}}
    \end{bmatrix}_{N_c\times NM}, \
    \mbs b=\begin{bmatrix}\mbs b^{\text{pde}}\\ \mbs b^{\text{bc}}\\ \mbs b^{\text{mea}}\\ \mbs 0
    \end{bmatrix}_{N_c\times 1}, \
    \mbs H^{\text{bc}} = \begin{bmatrix}\mbs H^{\text{bc1}}\\ \mbs H^{\text{bc2}}\\
      \mbs H^{\text{bc3}}\\ \mbs H^{\text{bc4}} \end{bmatrix}, \
    \mbs H^{\text{ck}} = \begin{bmatrix}\mbs H^{\text{ck1}}\\ \mbs H^{\text{ck2}}\\
      \mbs H^{\text{ck3}}\\ \mbs H^{\text{ck4}} \end{bmatrix}, \
    \mbs b^{\text{bc}}=\begin{bmatrix}\mbs b^{\text{bc1}}\\ \mbs b^{\text{bc2}}\\ \mbs b^{\text{bc3}}
    \\ \mbs b^{\text{bc4}}
    \end{bmatrix},
  \end{split}
\end{equation}
and the specific forms for these matrices 
are provided in the Appendix B.

For any given $\bm\alpha$ the least squares solution (with minimum norm)
to the system~\eqref{eq_42} is,
\begin{equation}\label{eq_44}\small
  \bm\beta = \mbs H^+(\bm\alpha)\mbs b.
\end{equation}
Substitution of this expression into~\eqref{eq_42} results in a reduced system
about $\bm\alpha$ only, with a residual given by
\begin{equation}\small
  \mbs r(\bm\alpha) = \mbs H(\bm\alpha) \mbs H^+(\bm\alpha)\mbs b - \mbs b.
\end{equation}
We determine the optimum $\bm\alpha^*$ by minimizing the Euclidean norm
of this residual,
\begin{equation}\label{eq_46}\small
  \bm\alpha^* = \argmin_{\bm\alpha}\frac12\|\mbs r(\bm\alpha) \|^2
  = \argmin_{\bm\alpha}\frac12\|\mbs H(\bm\alpha) \mbs H^+(\bm\alpha)\mbs b - \mbs b  \|^2.
\end{equation}
After $\bm\alpha$ is obtained, we compute $\bm\beta$ by solving
the system~\eqref{eq_42} with the linear least squares method.

The problem~\eqref{eq_46} is a nonlinear least squares problem about $\bm\alpha$. We
employ the NLLSQ-perturb algorithm (Algorithm~\ref{alg_1})
to solve this problem. In light of~\eqref{eq_32},
we can obtain the Jacobian matrix for this problem,
\begin{equation}\label{eq_47}\footnotesize
  \frac{\partial\mbs r}{\partial\bm\alpha} \approx \mbs J_0(\bm\alpha) - \mbs J_1(\bm\alpha),
  \quad \mbs J_0(\bm\alpha) = \frac{\partial\mbs H}{\partial\bm\alpha}\mbs H^+(\bm\alpha)\mbs b,
  \quad \mbs J_1(\bm\alpha) = \mbs H(\bm\alpha) \mbs H^+(\bm\alpha)\mbs J_0(\bm\alpha).
\end{equation}
$\mbs J_0(\bm\alpha)$ can be computed as follows. For any given $\bm\alpha$,
let $\bm\beta^{LS}=\mbs H^+(\bm\alpha)\mbs b=((\bm\beta_1^{LS})^T,\dots,(\bm\beta_N^{LS})^T)^T$
denote a {\em constant} vector. Then
\begin{equation}\label{eq_48}\footnotesize
  \left\{
  \begin{split}
    &
  \mbs J_0(\bm\alpha) =
  \frac{\partial\left[\mbs H(\bm\alpha)\bm\beta^{LS} \right]}{\partial\bm\alpha}
  =\frac{\partial}{\partial\bm\alpha}\begin{bmatrix}
    \mbs H^{\text{pde}}(\bm\alpha)\bm\beta^{LS}\\
    \mbs H^{\text{bc}}\bm\beta^{LS} \\
    \mbs H^{\text{mea}}\bm\beta^{LS} \\
    \mbs H^{\text{ck}}\bm\beta^{LS}
  \end{bmatrix}_{N_c\times 1}
  = \begin{bmatrix}
    \frac{\partial\left[\mbs H^{\text{pde}}(\bm\alpha)\bm\beta^{LS}\right]}{\partial\bm\alpha}\\
    \mbs 0 \\
    \mbs 0 \\
    \mbs 0
  \end{bmatrix}_{N_c\times n}, \\
  &
  \frac{\partial}{\partial\bm\alpha} \left[\mbs H^{\text{pde}}(\bm\alpha)\bm\beta^{LS}\right]
  = \frac{\partial}{\partial\bm\alpha}\begin{bmatrix}
    \vdots \\
    \alpha_1\mathcal{L}_1u^{LS}_e(\mbs x_p^e) + \dots +
    \alpha_n\mathcal{L}_nu^{LS}_e(\mbs x_p^e) + \mathcal{F}u^{LS}_e(\mbs x_p^e) \\
    \vdots
  \end{bmatrix}_{NQ\times 1} \\
  &\qquad\qquad\qquad\qquad\ 
  =\begin{bmatrix}
  \vdots &  & \vdots \\
  \mathcal{L}_1u^{LS}_e(\mbs x_p^e) & \cdots & \mathcal{L}_nu^{LS}_e(\mbs x_p^e) \\
  \vdots &  & \vdots
  \end{bmatrix}_{NQ\times n},
  \end{split}
  \right.
\end{equation}
where $u_e^{LS}(\mbs x)=\bm\Phi_e(\mbs x)\bm\beta_e^{LS}$ for $1\leqslant e\leqslant N$.
We compute $\mbs J_1(\bm\alpha)$ by the following two equations,
\begin{subequations}
\begin{align}\footnotesize
  &
  \mbs H(\bm\alpha)\mbs K = \mbs J_0(\bm\alpha) \label{eq_49a} \\
  & \mbs J_1(\bm\alpha) = \mbs H(\bm\alpha) \mbs K. \label{eq_49b}
\end{align}
\end{subequations}
We first solve~\eqref{eq_49a} for the $n\times n$ matrix $\mbs K$
by the linear least squares method, and then compute $\mbs J_1(\bm\alpha)$
by~\eqref{eq_49b} with a matrix multiplication.

\begin{algorithm}[tb]\small
  \DontPrintSemicolon
  \SetKwInOut{Input}{input}\SetKwInOut{Output}{output}

  \Input{ $\bm\alpha$;
    $\bm\Phi_e(x_p^e)$ and derivatives $(1\leqslant(e,p)\leqslant(N,Q))$;
    $\mathcal{M}\bm\Phi_e(\bm\xi_p^e)$ $(1\leqslant(e,p)\leqslant(N,Q_s))$.
  }
  \Output{reduced residual $\mbs r(\bm\alpha)$}
  \BlankLine
  \eIf{$\bm\alpha = \bm\alpha_s$}{
    retrieve $\mbs H(\bm\alpha_s)$, $\mbs b$, $\bm\beta^{LS}$\;
    set $\mbs H(\bm\alpha)= \mbs H(\bm\alpha_s)$\;
  }{
    compute $\mbs H(\bm\alpha)$ and $\mbs b$ by~\eqref{eq_43} and the
    equations~\eqref{eq_ab1}--\eqref{eq_ab4} in Appendix B\;
  solve equation~\eqref{eq_42} for $\bm\beta$ by the linear least squares method,
  and let $\bm\beta^{LS}=\bm\beta$\;
  set $\bm\alpha_s=\bm\alpha$, and save $\mbs H(\bm\alpha)$, $\mbs b$, and $\bm\beta^{LS}$\;
  }
  \BlankLine
  compute $\mbs r(\bm\alpha)$ by $\mbs r(\bm\alpha)=\mbs H(\bm\alpha)\bm\beta^{LS}-\mbs b$\;
  
  \caption{Computing reduced residual $\mbs r(\bm\alpha)$
    for VarPro-F2.
  }
  \label{alg_6}
\end{algorithm}

The procedures for computing the residual $\mbs r(\bm\alpha)$
and the Jacobian matrix $\frac{\partial\mbs r}{\partial\bm\alpha}$
for the reduced problem~\eqref{eq_46} are summarized in the Algorithms~\ref{alg_6}
and~\ref{alg_7}.

\begin{algorithm}[tb]\small
  \DontPrintSemicolon
  \SetKwInOut{Input}{input}\SetKwInOut{Output}{output}

  \Input{ $\bm\alpha$;
    $\bm\Phi_e(x_p^e)$ and derivatives $(1\leqslant(e,p)\leqslant(N,Q))$;
    $\mathcal{M}\bm\Phi_e(\bm\xi_p^e)$ $(1\leqslant(e,p)\leqslant(N,Q_s))$.
  }
  \Output{Jacobian matrix $\frac{\partial\mbs r}{\partial\bm\alpha}$ }
  \BlankLine
  \eIf{$\bm\alpha = \bm\alpha_s$}{
    retrieve $\mbs H(\bm\alpha_s)$, $\mbs b$, $\bm\beta^{LS}$\;
    set $\mbs H(\bm\alpha)= \mbs H(\bm\alpha_s)$\;
  }{
    compute $\mbs H(\bm\alpha)$ and $\mbs b$ by~\eqref{eq_43} and the
    equations~\eqref{eq_ab1}--\eqref{eq_ab4} in Appendix B\;
  solve equation~\eqref{eq_42} for $\bm\beta$ by the linear least squares method,
  and let $\bm\beta^{LS}=\bm\beta$\;
  set $\bm\alpha_s=\bm\alpha$, and save $\mbs H(\bm\alpha)$, $\mbs b$, and $\bm\beta^{LS}$\;
  }
  \BlankLine
  compute $u_e^{LS}(\mbs x_p^e)$ $(1\leqslant (e,p)\leqslant (N,Q))$ and 
  their derivatives by~\eqref{eq_011} based on $\bm\beta^{LS}$\;
  compute $\mbs J_0(\bm\alpha)$ by~\eqref{eq_48}\;
  compute $\mbs J_1(\bm\alpha)$ by~\eqref{eq_49a}--\eqref{eq_49b}\;
  compute $\frac{\partial\mbs r}{\partial\bm\alpha}$ by~\eqref{eq_47}\;  
  
  \caption{Computing Jacobian matrix $\frac{\partial\mbs r}{\partial\bm\alpha}$
    for VarPro-F2.
  }
  \label{alg_7}
\end{algorithm}

The overall VarPro-F2 algorithm consists of two steps:
(i) Invoke the NLLSQ-perturb algorithm (Algorithm~\ref{alg_1} in Appendix A)
to solve the problem~\eqref{eq_46} for $\bm\alpha$, with
the routines in Algorithms~\ref{alg_6} and~\ref{alg_7} as input
arguments. (ii) Solve equation~\eqref{eq_42} for $\bm\beta$ by
the linear least squares method.

\begin{remark}\label{rem_6}
  Let us now discuss an extension of the above algorithm
  to deal with the case in which some (or all) of the operators of
  $\mathcal{L}_i$ ($1\leqslant i\leqslant n$) and $\mathcal{F}$ are
  nonlinear with respect to $u$. In this case, we can first use a Newton
  iteration to linearize the nonlinear operators, and then solve the linearized
  system by the VarPro-F2 algorithm as discussed above. Upon convergence
  of the Newton iteration, the solution for $(\bm\alpha,\bm\beta)$ to
  the original system will be attained.
  To make the discussion more concrete and without loss of generality, let us
  assume that $\mathcal{L}_1$ and $\mathcal{F}$ are nonlinear while the other
  operators are linear. Let $u_e^{k}(\mbs x)$ ($1\leqslant e\leqslant N$)
  denote the approximation of $u_e(\mbs x)$ at the $k$-th Newton step.
  Equation~\eqref{eq_11} is nonlinear with respect to $u$, and its
  linearized form is given by,
  \begin{equation}\label{eq_50}\footnotesize
    \begin{split}
      &
    \alpha_1\mathcal{L}_1'(u_e^{k}(\mbs x_p^e))u_e^{k+1}(\mbs x_p^e)
    +\alpha_2\mathcal{L}_2u_e^{k+1}(\mbs x_p^e) + \dots
    + \alpha_n\mathcal{L}_nu_e^{k+1}(\mbs x_p^e)
    + \mathcal{F}'(u_e^k(\mbs x_p^e))u_e^{k+1}(\mbs x_p^e) \\
    &\qquad
    -\left[ f(\mbs x_p^e)
      - \alpha_1\mathcal{L}_1(u_e^k(\mbs x_p^e))
      + \alpha_1\mathcal{L}_1'(u_e^k(\mbs x_p^e))u_e^k(\mbs x_p^e)
      - \mathcal{F}(u_e^k(\bs x_p^e)) + \mathcal{F}'(u_e^k(\mbs x_p^e))u_e^k(\mbs x_p^e)
      \right] = 0, \\
    &\qquad \text{for}\ 1\leqslant (e,p)\leqslant (N,Q).
    \end{split}
  \end{equation}
  Notice that this equation is linear with respect to $u_e^{k+1}$.
  The equations~\eqref{eq_12}--\eqref{eq_14} are linear with respect to $u_e$, and we
  enforce them on the $(k+1)$-th Newton step (i.e.~replacing $u_e$ by $u_e^{k+1}$ in these
  equations).
  The system consisting of~\eqref{eq_50} and the equations~\eqref{eq_12}--\eqref{eq_14}
  (written in terms of $u_e^{k+1}$) are linear with respect to the updated
  approximation field $u_e^{k+1}$. With the expansion
  $u_e^{k+1}(\mbs x)=\bm\Phi_e(\mbs x)\bm\beta_e^{k+1}$, we can solve this system
  for $(\bm\alpha,\bm\beta^{k+1})$ by the VarPro-F2 algorithm
  as discussed above. Upon convergence of the Newton iteration, the solution
  to $(\bm\alpha,\bm\beta)$ is given by the converged result, and
  the neural network coefficients contains the representation for
  the field solution $u(\mbs x)$ to the original nonlinear system.
  For inverse nonlinear PDEs with respect to $u$,
  the combination of the Newton iteration and the VarPro-F2 algorithm in general
  works quite well. We have also observed from numerical experiments that
  for certain problems it appears to be somewhat less robust than
  the VarPro-F1 and NLLSQ methods,
  leading to less accurate results than VarPro-F1 and NLLSQ.

\end{remark}


\section{Numerical Examples}
\label{sec:tests}

In this section we test the presented method and
algorithms using several inverse PDE problems
 in two dimensions (2D) or in one spatial dimension (1D)
plus time. The Gaussian activation function, $\sigma(x)=e^{-x^2}$,
is employed in all the neural networks. We fix the seed value at $25$ in
the random number generator for all the test problems, so that
the reported results here are exactly reproducible.
Note that $\lambda_{mea}$ denotes the scaling coefficient
for the measurement residual (see Remark~\ref{rem_04}), with $\lambda_{mea}=1$
corresponding to the case of no scaling.
We refer the reader to the Appendix C for a comparison between
the current method and the PINN method with several of these test problems.

\subsection{Parametric Poisson Equation}
\label{poisson}

Consider the domain $(x,y)\in\Omega=[0,1.4]\times[0,1.4]$, and the inverse problem,
\begin{subequations}\label{eq_51}\small
  \begin{align}
    &
    \frac{\partial^2u}{\partial x^2} + \alpha\frac{\partial^2u}{\partial y^2}
    = f(x,y), \label{eq_51a} \\
    &
    u(0,y) = g_1(y), \quad
    u(1.4,y) = g_2(y), \quad
    u(x,0) = g_3(x), \quad
    u(x,1.4) = g_4(x), \label{eq_51b} \\
    &
    u(\xi_i,\eta_i) = S(\xi_i,\eta_i), \quad
    (\xi_i,\eta_i)\in\mbb Y, \quad
    1\leqslant i\leqslant NQ_s, \label{eq_51c}
  \end{align}
\end{subequations}
where $f$ and $g_i$ ($1\leqslant i\leqslant 4$) denote a source term
and the boundary data respectively,
$\mbb Y\subset\Omega$ denotes the set of random measurement points,
$\alpha$ and $u(x,y)$ are the unknowns to be solved for, $N$ denotes the number of sub-domains,
and $Q_s$ is the number of measurement points per sub-domain.
We use the following manufactured solution to this problem,
\begin{equation}\label{eq_52}\small
  \begin{split}
    &
    \alpha_{ex} = 1, \quad
    u_{ex}(x,y) = \sin(\pi x^2)\sin(\pi y^2).
  \end{split}
\end{equation}
The source term and the boundary data are chosen such that the expressions in~\eqref{eq_52}
satisfy~\eqref{eq_51a}--\eqref{eq_51b}.
The measurement data are taken to be
\begin{equation}\label{eq_53}\small
  S(\xi_i,\eta_i) = u_{ex}(\xi_i,\eta_i)(1+\epsilon\zeta_i),
  \quad 1\leqslant i\leqslant NQ_s,
\end{equation}
where $\zeta_i$ denotes a uniform random number from $[-1,1]$ representing
the noise and the constant $\epsilon\geqslant 0$
denotes the relative level of the noise.

\begin{figure}
  \small
  \centerline{
    \includegraphics[width=2.0in]{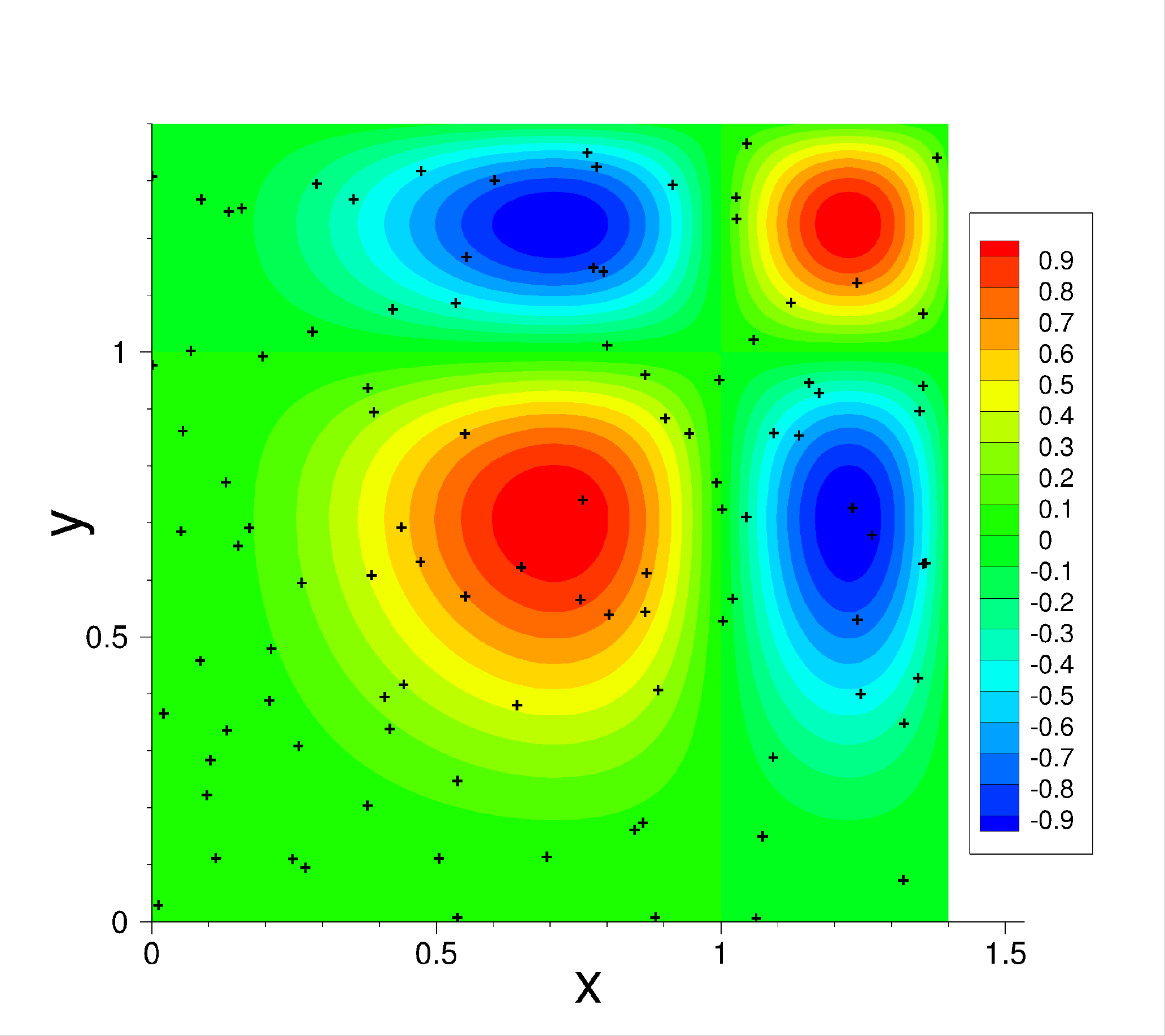}(a)
    \includegraphics[width=2.0in]{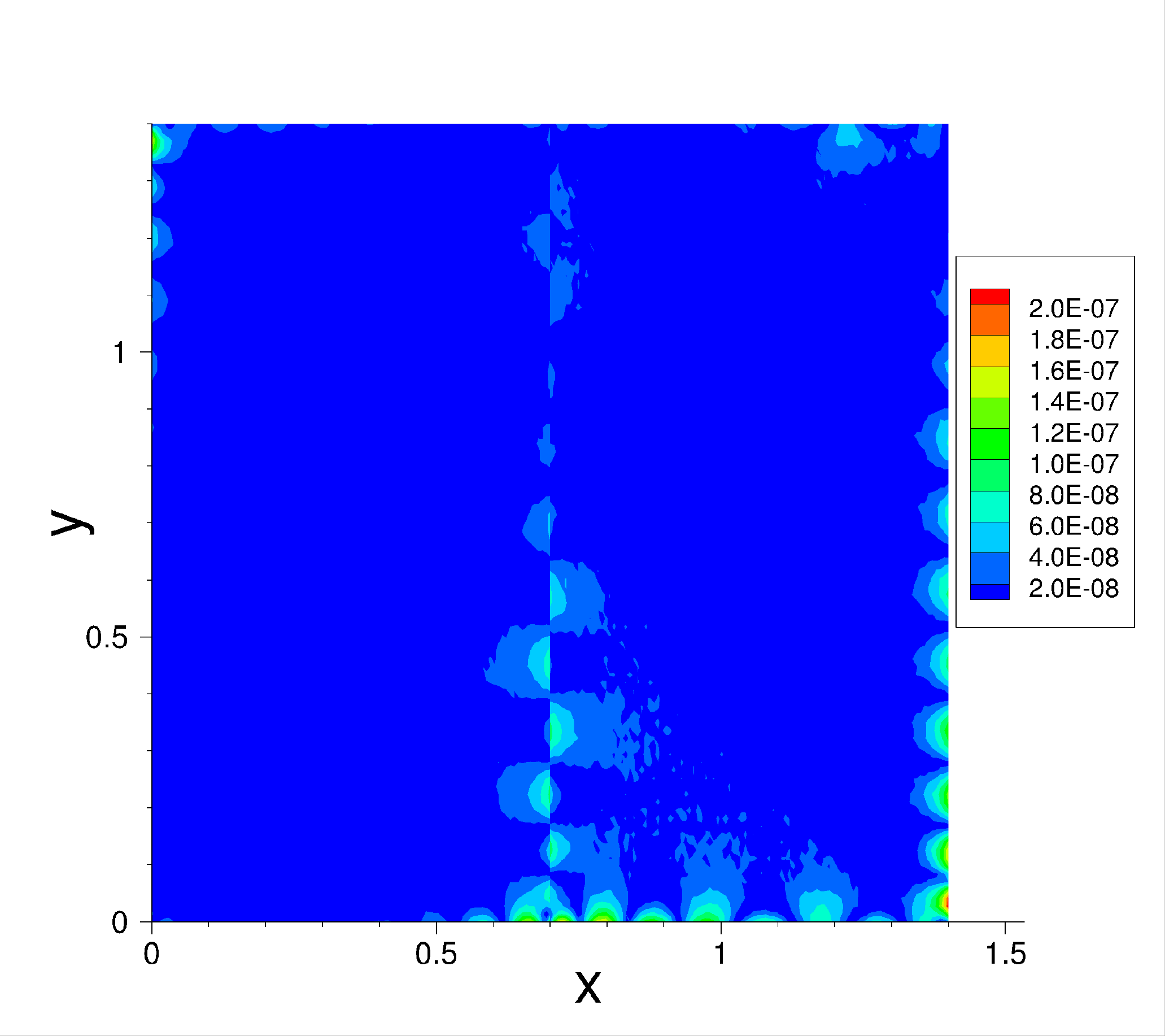}(b)
  }
  \caption{\small Inverse Poisson problem: distributions of (a) the NLLSQ solution and (b) its point-wise
    absolute error, with the random measurement points shown in (a)
    as ``+'' symbols.
    Two uniform sub-domains (along $x$), local NN $[2,400,1]$, $Q=25\times 25$,
    $Q_s=50$, $R_m=2.0$, $\lambda_{mea}$=1, $\epsilon$=0 (no noise in measurement data).
  }
  \label{fg_3}
\end{figure}

Henceforth $Q$ denotes the number of uniform collocation
points per  sub-domain, $Q_s$ denotes the number of random measurement points per sub-domain,
$\epsilon$ denotes the noise level, and $M$ denotes the number of trainable parameters
 of each local NN.
 $R_m$ denotes a constant, and the hidden-layer coefficients are
 assigned to  uniform random values generated on $[-R_m,R_m]$.
The $R_m$ values employed in the tests are obtained by the method
from~\cite{DongY2022rm}, as noted in Remark~\ref{rem_4}.
After the NN is trained, it is evaluated on another set of $Q_{eval}=101\times 101$
uniform grid points (evaluation points) on each sub-domain to obtain $u$, which is compared with~\eqref{eq_52}
to compute the errors. The relative errors of $\alpha$ ($e_{\alpha}$)
and $u$ ($l^{\infty}$ and $l^2$ norms) are defined as,
\begin{equation}\footnotesize\label{eq_a55}
  e_{\alpha}=\frac{\left|\alpha-\alpha_{ex} \right|}{\left|\alpha_{ex} \right|}, \
  l^{\infty}\text{-u} = \frac{\max\left\{\left|u(\mbs x_i)-u_{ex}(\mbs x_i) \right|\right\}_{i=1}^{NQ_{eval}}}{\sqrt{\frac{1}{NQ_{eval}}\sum_{i=1}^{NQ_{eval}}u_{ex}^2(\mbs x_i)}}, \
  l^{2}\text{-u} = \frac{\sqrt{\frac{1}{NQ_{eval}}\sum_{i=1}^{NQ_{eval}}[u(\mbs x_i)-u_{ex}(\mbs x_i)]^2}}{\sqrt{\frac{1}{NQ_{eval}}\sum_{i=1}^{NQ_{eval}}u_{ex}^2(\mbs x_i)}},
\end{equation}
where $N$ is the number of sub-domains and $\mbs x_i$ denotes the evaluation points.

Figure~\ref{fg_3} illustrates $u(x,y)$ and its point-wise absolute
error obtained by the NLLSQ algorithm with $2$ sub-domains.
The caption lists the main simulation parameters.
In particular, the random measurement
points ($100$ total) are shown in Figure~\ref{fg_3}(a), and the there is no noise in
the measurement data. The NLLSQ solution for $u$ is quite accurate,
with a maximum error on the order of $10^{-7}$ in the domain. 
The relative (or absolute) error of the computed $\alpha$ is $9.03\times 10^{-9}$.


\begin{table}[tb]
  \small\centering
  \begin{tabular}{llll}
    \hline
    $Q$ & $\alpha$ (NLLSQ) & $\alpha$ (VarPro-F1) & $\alpha$ (VarPro-F2)  \\ \hline
    5$\times$5 & 1.076466245043E+0 & 9.982719409724E-1 & 0.000000000000E+0 \\
    10$\times$10 & 9.999867935849E-1 & 9.999965494049E-1 & -3.188390321381E-5 \\
    15$\times$15 & 1.000000029498E+0 & 9.999999954822E-1 & 9.999999998978E-1 \\
    20$\times$20 & 9.999999999701E-1 & 9.999999999592E-1 & 9.999999999536E-1 \\
    25$\times$25 & 9.999999987249E-1 & 1.000000000817E+0 & 1.000000001279E+0 \\
    30$\times$30 & 1.000000002811E+0 & 1.000000000906E+0 & 1.000000000002E+0 \\
    35$\times$35 & 1.000000001708E+0 & 1.000000000670E+0 & 1.000000000237E+0 \\
    40$\times$40 & 1.000000001552E+0 & 1.000000000717E+0 & 1.000000000183E+0
    \\ \hline
  \end{tabular}
  \caption{\small Inverse Poisson problem: computed $\alpha$ by
    the NLLSQ, VarPro-F1 and VarPro-F2 algorithms versus $Q$ (number of collocation
    points).  Single sub-domain, NN $[2, 600, 1]$,
    $Q_s=100$, $\lambda_{mea}$=1, $\epsilon$=0; $R_m=3.0$ with NLLSQ, $R_m=2.8$ with
    VarPro-F1, and $R_m=2.0$ with VarPro-F2.
  }
  \label{tab_1}
\end{table}

\begin{table}[tb]
  \small\centering
  \begin{tabular}{l|ccc|ccc|ccc}
    \hline
     & & NLLSQ & & & VarPro-F1 & & & VarPro-F2 & \\ \cline{2-10}
    $Q$ & $e_\alpha$ & $l^{\infty}$-u & $l^2$-u & $e_{\alpha}$ & $l^{\infty}$-u & $l^2$-u 
    & $e_{\alpha}$ & $l^{\infty}$-u & $l^2$-u  \\ \hline
    5$\times$5 & 7.65E-2 & 1.98E-1 & 2.59E-2 & 1.73E-3 & 5.70E-1 & 3.00E-2 & 1.00E+0 & 8.23E+1 & 6.65E+0  \\
    10$\times$10 & 1.32E-5 & 9.01E-3 & 7.88E-4 & 3.45E-6 & 1.14E-2 & 6.32E-4 & 1.00E+0 & 1.06E+3 & 5.08E+1 \\
    15$\times$15 & 2.95E-8 & 5.15E-5 & 3.50E-6 & 4.52E-9 & 3.50E-5 & 2.15E-6 & 1.02E-10 & 4.19E-6 & 2.40E-7 \\
    20$\times$20 & 2.99E-11 & 2.98E-7 & 1.56E-8 & 4.08E-11 & 1.93E-7 & 1.30E-8 & 4.64E-11 & 2.06E-7 & 1.05E-8 \\
    25$\times$25 & 1.28E-9 & 7.06E-8 & 8.17E-9 & 8.17E-10 & 4.46E-8 & 4.00E-9 & 1.28E-9 & 7.37E-8 & 5.96E-9 \\
    30$\times$30 & 2.81E-9 & 9.13E-8 & 7.13E-9 & 9.06E-10 & 6.75E-8 & 3.61E-9 & 2.44E-12 & 7.58E-8 & 7.14E-9 \\
    35$\times$35 & 1.71E-9 & 1.53E-7 & 1.16E-8 & 6.70E-10 & 1.02E-7 & 5.62E-9 & 2.37E-10 & 1.75E-7 & 1.55E-8 \\
    40$\times$40 & 1.55E-9 & 2.09E-7 & 1.69E-8 & 7.17E-10 & 1.64E-7 & 9.29E-9 & 1.83E-10 & 1.16E-7 & 1.19E-8
    \\ \hline
  \end{tabular}
  \caption{\small Inverse Poisson problem: relative errors of $\alpha$ and $u$ versus $Q$ computed
    by the NLLSQ, VarPro-F1 and VarPro-F2 algorithms.
    Single sub-domain, NN $[2, 600, 1]$,
    $Q_s=100$, $\lambda_{mea}$=1, $\epsilon$=0; $R_m=3.0$ with NLLSQ, $R_m=2.8$ with VarPro-F1,
    and $R_m=2.0$ with VarPro-F2.
  }
  \label{tab_2}
\end{table}

The convergence of computation results with respect to $Q$ (number of collocation points)
is illustrated by Tables~\ref{tab_1} and~\ref{tab_2}.
Table~\ref{tab_1} lists the computed $\alpha$ values versus $Q$
by the NLLSQ, VarPro-F1 and VarPro-F2 methods.
Table~\ref{tab_2} lists the relative errors of $\alpha$  and $u$
with respect to $Q$ from the three methods.
The main parameters values for these tests are provided in the table captions.
The $\alpha$ and the $u$ errors generally decrease exponentially with
increasing $Q$, until $Q$ reaches a certain level. The errors generally stagnate
as $Q$ further increases beyond that point.


\begin{table}[tb]
  \small\centering
  \begin{tabular}{l|ccc|ccc|ccc}
    \hline
     & & NLLSQ & & & VarPro-F1 & & & VarPro-F2 & \\ \cline{2-10}
    $M$ & $e_\alpha$ & $l^{\infty}$-u & $l^2$-u & $e_{\alpha}$ & $l^{\infty}$-u & $l^2$-u 
    & $e_{\alpha}$ & $l^{\infty}$-u & $l^2$-u  \\ \hline
    100 & 1.05E+4 & 2.87E+0 & 1.02E+0 & 1.61E+6 & 2.46E+0 & 9.98E-1 & 5.66E-1 & 2.20E+0 & 5.08E-1 \\
    200 & 1.41E-2 & 1.99E-1 & 3.23E-2 & 3.19E-4 & 1.35E-1 & 1.89E-2 & 3.91E-4 & 3.27E-2 & 4.90E-3 \\
    300 & 2.77E-5 & 2.83E-3 & 4.00E-4 & 8.68E-6 & 1.74E-3 & 1.96E-4 & 6.93E-7 & 3.27E-4 & 2.39E-5 \\
    400 & 2.82E-7 & 5.68E-5 & 4.36E-6 & 6.97E-9 & 2.37E-5 & 2.09E-6 & 9.53E-9 & 3.39E-6 & 2.21E-7 \\
    500 & 1.20E-8 & 1.81E-6 & 1.68E-7 & 5.68E-8 & 7.28E-7 & 1.07E-7 & 1.00E-8 & 4.03E-7 & 2.26E-8 \\
    600 & 1.28E-9 & 7.06E-8 & 8.17E-9 & 8.17E-10 & 4.46E-8 & 4.00E-9 & 1.28E-9 & 7.37E-8 & 5.96E-9
    \\ \hline
  \end{tabular}
  \caption{\small Inverse Poisson problem: $\alpha$ and $u$ relative errors versus $M$ (number of training parameters)
    obtained by the NLLSQ, VarPro-F1 and VarPro-F2 algorithms.
    Single sub-domain, NN $[2, M, 1]$,
    $Q=25\times25$, $Q_s=100$, $\lambda_{mea}$=1, $\epsilon$=0;
    $R_m=3.0$ with NLLSQ, $R_m=2.8$ with VarPro-F1, and $R_m=2.0$ with VarPro-F2.
  }
  \label{tab_3}
\end{table}

The convergence of the NLLSQ, VarPro-F1 and VarPro-F2 algorithms with respect to
the number of trainable parameters $M$ is illustrated by Table~\ref{tab_3}.
A single sub-domain and a single hidden layer in the neural network are
employed in the simulations, where the number of hidden nodes ($M$) is varied.
The caption lists the crucial parameter values.
It is evident that the errors for $\alpha$ and $u$
decrease exponentially with increasing number of
training parameters.


\begin{table}[tb]
  \small\centering
  \begin{tabular}{l|ccc|ccc|ccc}
    \hline
     & & NLLSQ & & & VarPro-F1 & & & VarPro-F2 & \\ \cline{2-10}
    $Q_s$ & $e_\alpha$ & $l^{\infty}$-u & $l^2$-u & $e_{\alpha}$ & $l^{\infty}$-u & $l^2$-u 
    & $e_{\alpha}$ & $l^{\infty}$-u & $l^2$-u  \\ \hline
    1 & 1.02E+0 & 1.01E+0 & 3.27E-1 & 1.63E+0 & 2.09E+1 & 7.27E+0 & 6.61E-4 & 1.19E-3 & 3.59E-4 \\
    2 & 5.01E-7 & 3.67E-6 & 3.80E-7 & 1.70E+0 & 3.10E+1 & 1.15E+1 & 1.67E-8 & 5.60E-7 & 3.49E-8 \\
    3 & 5.06E-8 & 3.67E-6 & 2.64E-7 & 9.64E-8 & 9.11E-7 & 9.59E-8 & 4.17E-9 & 5.46E-7 & 3.26E-8 \\
    5 & 3.06E-8 & 3.67E-6 & 2.62E-7 & 2.81E-8 & 8.88E-7 & 8.26E-8 & 7.30E-9 & 5.33E-7 & 3.25E-8 \\
    10 & 1.39E-8 & 3.67E-6 & 2.61E-7 & 1.34E-8 & 9.05E-7 & 8.12E-8 & 4.72E-8 & 5.39E-7 & 4.16E-8 \\
    20 & 5.14E-8 & 3.67E-6 & 2.63E-7 & 1.19E-8 & 9.04E-7 & 8.10E-8 & 2.62E-9 & 5.47E-7 & 3.46E-8 \\
    50 & 1.07E-8 & 3.67E-6 & 2.62E-7 & 6.63E-9 & 8.56E-7 & 8.02E-8 & 4.81E-10 & 5.33E-7 & 3.23E-8 \\
    100 & 3.26E-8 & 3.72E-6 & 2.62E-7 & 3.15E-9 & 9.29E-7 & 7.95E-8 & 1.17E-8 & 5.48E-7 & 3.27E-8
    \\ \hline
  \end{tabular}
  \caption{\small Inverse Poisson problem: $\alpha$ and $u$ relative errors versus $Q_s$ (number of measurement points)
    for the NLLSQ, VarPro-F1 and VarPro-F2 algorithms.
    Single sub-domain, neural network $[2, 500, 1]$,
    $Q=30\times 30$;
    $R_m=3.0$ with NLLSQ, $R_m=2.8$ with VarPro-F1, and $R_m=2.0$ with VarPro-F2;
    $\lambda_{mea}$=1, $\epsilon$=0.
  }
  \label{tab_4}
\end{table}

Table~\ref{tab_4} illustrates the effect of the number of random measurement points ($Q_s$)
on the $\alpha$ and $u$ errors computed by the NLLSQ, VarPro-F1 and VarProf-F2 algorithms.
When $Q_s$ is very small, the computed $\alpha$ and $u$
are inaccurate or less accurate. On the other hand, when $Q_s$ reaches a
certain value ($Q_s=3$ for this problem)
and beyond, the three algorithms produce highly accurate results.
This seems to be a common characteristic of these algorithms for all the test problems
considered in this work.


\begin{table}[tb]
  \small\centering
  \begin{tabular}{ll | ll | ll}
    \hline
    $\epsilon$ & computed-$\alpha$ & $\epsilon$ & computed-$\alpha$ &
    $\epsilon$ & computed-$\alpha$ \\ \hline
    0.0 & 9.99999993208E-1 & 0.01 & 9.9875752E-1 & 0.1 & 9.8779390E-1 \\
    0.001 & 9.9987537E-1 & 0.03 & 9.9630066E-1 & 0.2 & 9.7602056E-1  \\
    0.002 & 9.9975066E-1 & 0.05 & 9.9383633E-1 & 0.5 & 9.4329282E-1  \\
    0.005 & 9.9937764E-1 & 0.07 & 9.9139103E-1 & 0.7 & 9.2316247E-1  \\
    0.007 & 9.9912874E-1 & 0.09 & 9.8897497E-1 & 1.0 & 8.9557261E-1 \\
    \hline
  \end{tabular}
  \caption{\small Inverse Poisson problem: $\alpha$ obtained by
    the NLLSQ algorithm corresponding to several noise levels ($\epsilon$).
    Single sub-domain, NN $[2, 500, 1]$,
    $Q=25\times 25$, $Q_s=50$, $R_m=3.0$, $\lambda_{mea}$=1.
  }
  \label{tab_5}
\end{table}

\begin{table}[tb]
  \small\centering
  \begin{tabular}{l|ccc|ccc|ccc}
    \hline
     & & NLLSQ & & & VarPro-F1 & & & VarPro-F2 & \\ \cline{2-10}
    $\epsilon$ & $e_\alpha$ & $l^{\infty}$-u & $l^2$-u & $e_{\alpha}$ & $l^{\infty}$-u & $l^2$-u 
    & $e_{\alpha}$ & $l^{\infty}$-u & $l^2$-u  \\ \hline
    0.0 & 6.79E-9 & 1.81E-6 & 1.93E-7 & 5.93E-8 & 7.59E-7 & 1.38E-7 & 4.20E-10 & 4.50E-7 & 2.31E-8 \\
    0.001 & 1.25E-4 & 2.79E-4 & 8.04E-5 & 1.33E-4 & 2.82E-4 & 8.50E-5 & 1.23E-4 & 2.81E-4 & 7.85E-5 \\
    0.005 & 6.22E-4 & 1.39E-3 & 4.01E-4 & 6.73E-4 & 1.42E-3 & 4.29E-4 & 6.10E-4 & 1.41E-3 & 3.92E-4 \\
    0.01 & 1.24E-3 & 2.79E-3 & 8.02E-4 & 1.35E-3 & 2.85E-3 & 8.61E-4 & 1.22E-3 & 2.81E-3 & 7.82E-4 \\
    0.05 & 6.16E-3 & 1.39E-2 & 4.00E-3 & 6.63E-3 & 1.42E-2 & 4.26E-3 & 6.49E-3 & 1.41E-2 & 4.07E-3 \\
    0.1 & 1.22E-2 & 2.79E-2 & 7.98E-3 & 1.33E-2 & 2.86E-2 & 8.58E-3 & 1.19E-2 & 2.81E-2 & 7.75E-3 \\
    0.5 & 5.67E-2 & 1.42E-1 & 3.90E-2 & 6.08E-2 & 1.44E-1 & 4.17E-2 & 5.52E-2 & 1.43E-1 & 3.78E-2 \\
    1.0 & 1.04E-1 & 2.88E-1 & 7.63E-2 & 1.11E-1 & 2.93E-1 & 8.14E-2 & 1.08E-1 & 2.90E-1 & 7.62E-2
    \\ \hline
  \end{tabular}
  \caption{\small Inverse Poisson problem: $\alpha$ and $u$ relative errors versus $\epsilon$ (noise level)
    computed by the NLLSQ,
    VarPro-F1 and VarPro-F2 algorithms.
    Single sub-domain, NN $[2, 500, 1]$,
    $Q=25\times 25$, $Q_s=50$, $\lambda_{mea}$=1;
    $R_m=3.0$ with NLLSQ, $R_m=2.8$ with VarPro-F1, and $R_m=2.0$ with VarPro-F2.
  }
  \label{tab_6}
\end{table}

\begin{table}[tb]
  \small\centering
  \begin{tabular}{l|ccc|ccc|ccc}
    \hline
     & & $\lambda_{mea}$=0.5 & & & $\lambda_{mea}$=0.25 & & & $\lambda_{mea}$=0.1 & \\ \cline{2-10}
    $\epsilon$ & $e_\alpha$ & $l^{\infty}$-u & $l^2$-u & $e_{\alpha}$ & $l^{\infty}$-u & $l^2$-u 
    & $e_{\alpha}$ & $l^{\infty}$-u & $l^2$-u  \\ \hline
    0.0 & 1.69E-8 & 1.81E-6 & 2.27E-7 & 2.25E-8 & 1.82E-6 & 2.46E-7 & 1.34E-8 & 1.81E-6 & 2.49E-7 \\
    0.001 & 1.04E-4 & 1.80E-4 & 5.58E-5 & 9.79E-5 & 1.75E-4 & 5.23E-5 & 9.58E-5 & 1.73E-4 & 5.19E-5 \\
    0.005 & 5.22E-4 & 9.01E-4 & 2.79E-4 & 4.89E-4 & 8.72E-4 & 2.61E-4 & 4.79E-4 & 8.61E-4 & 2.59E-4 \\
    0.01 & 1.04E-3 & 1.80E-3 & 5.58E-4 & 9.76E-4 & 1.74E-3 & 5.22E-4 & 9.56E-4 & 1.72E-3 & 5.18E-4 \\
    0.05 & 5.16E-3 & 8.94E-3 & 2.77E-3 & 4.83E-3 & 8.65E-3 & 2.59E-3 & 4.73E-3 & 8.54E-3 & 2.57E-3 \\
    0.1 & 1.02E-2 & 1.78E-2 & 5.49E-3 & 9.54E-3 & 1.71E-2 & 5.13E-3 & 9.33E-3 & 1.69E-2 & 5.08E-3 \\
    0.5 & 4.66E-2 & 8.35E-2 & 2.59E-2 & 4.33E-2 & 8.00E-2 & 2.39E-2 & 4.23E-2 & 7.88E-2 & 2.36E-2 \\
    1.0 & 8.45E-2 & 1.56E-1 & 4.83E-2 & 7.78E-2 & 1.48E-1 & 4.40E-2 & 7.57E-2 & 1.45E-1 & 4.33E-2
    \\ \hline
  \end{tabular}
  \caption{\small Inverse Poisson problem: $\alpha$ and $u$ relative errors versus $\epsilon$ and
    $\lambda_{mea}$ (scaling coefficient of measurement residual) computed by the NLLSQ algorithm.
    Single sub-domain, NN $[2, 500, 1]$,
    $Q=25\times 25$, $Q_s=50$,
    $R_m=3.0$. These can be compared with the NLLSQ data in Table~\ref{tab_6}, which
    is with $\lambda_{mea}=1$.
  }
  \label{tab_06}
\end{table}

In the foregoing tests no noise is considered in the measurement data ($\epsilon=0$).
Tables~\ref{tab_5}, \ref{tab_6} and~\ref{tab_06} demonstrate the effect of noisy
measurement data on the computation results.
Table~\ref{tab_5} shows the computed $\alpha$ values by the NLLSQ algorithm
corresponding to different noise levels, ranging from $\epsilon=0$ ($0\%$)
to $\epsilon=1.0$ ($100\%$). Table~\ref{tab_6} lists the $\alpha$ errors and
the $u$ errors corresponding to several noise levels obtained by
the NLLSQ, VarPro-F1 and VarPro-F2 algorithms.
Table~\ref{tab_06} provides the $\alpha$ and $u$ relative errors corresponding to
different $\epsilon$ and several $\lambda_{mea}$ (scaling factor of measurement residual) values
with the NLLSQ algorithm.
The presence of noise degrades the simulation accuracy.
But the current method and these algorithms appear to be quite robust.
For example, with $10\%$ ($\epsilon=0.1$) noise in the measurement data
the relative error of  $\alpha$ is around $1\%$ for the three methods.
With $100\%$ ($\epsilon=1.0$) noise in the data, the computed $\alpha$
exhibits a relative error around $10\%$ with these algorithms.
For noisy data, scaling the measurement residual by $\lambda_{mea}$ can improve
the accuracy of computation results and make the method more robust (see Table~\ref{tab_06}), compared
with the case of no scaling.
A smaller $\lambda_{mea}$ in general leads to a better accuracy.

\subsection{Parametric Advection Equation}
\label{advec}

Consider the spatial-temporal domain, $(x,t)\in\Omega=[0,3]\times[0,1]$, and the
following inverse problem,
\begin{subequations}
  \begin{align}\small
    &
    \frac{\partial u}{\partial t} - c\frac{\partial u}{\partial x}=0, \\
    &
    u(0,t) = u(3,t), \quad
    u(x,0) = 10\sinh\left[\frac{1}{10}\sin\frac{2\pi}{3}\left(x - \frac52\right) \right], \\
    &
    u(\xi_i,\eta_i) = S(\xi_i,\eta_i), \quad (\xi_i,\eta_i)\in\mbb Y\subset\Omega,
    \ 1\leqslant i\leqslant NQ_s,
  \end{align}
\end{subequations}
where $\mbb Y$ denotes the set of measurement points in $\Omega$. The wave speed $c$ and the field
$u(x,t)$ are the unknowns to be determined in this problem.
We employ the following exact solution to this problem in the tests,
\begin{equation}\label{eq_55}\small
  \begin{split}
    &
    c_{ex} = 3, \quad
    u_{ex}(x,t) = 10\sinh\left[\frac{1}{10}\sin\frac{2\pi}{3}\left(x + 3t - \frac52\right) \right].
  \end{split}
\end{equation}
We employ random measurement points in $\Omega$,
and the measurement data are given by~\eqref{eq_53},
in which $u_{ex}$ is given by~\eqref{eq_55}.
The notations adopted below (e.g.~$Q$, $M$, $N$, $Q_s$, $R_m$, $\epsilon$)
are the same as in Section~\ref{poisson}.
The $l^{\infty}$ and $l^2$ norms of the $u$ relative error reported below
are computed on a set of $Q_{eval}=101\times 101$ uniform grid points in each sub-domain
after the network is trained.

\begin{figure}
  \small
  \centerline{
    \includegraphics[width=2.7in]{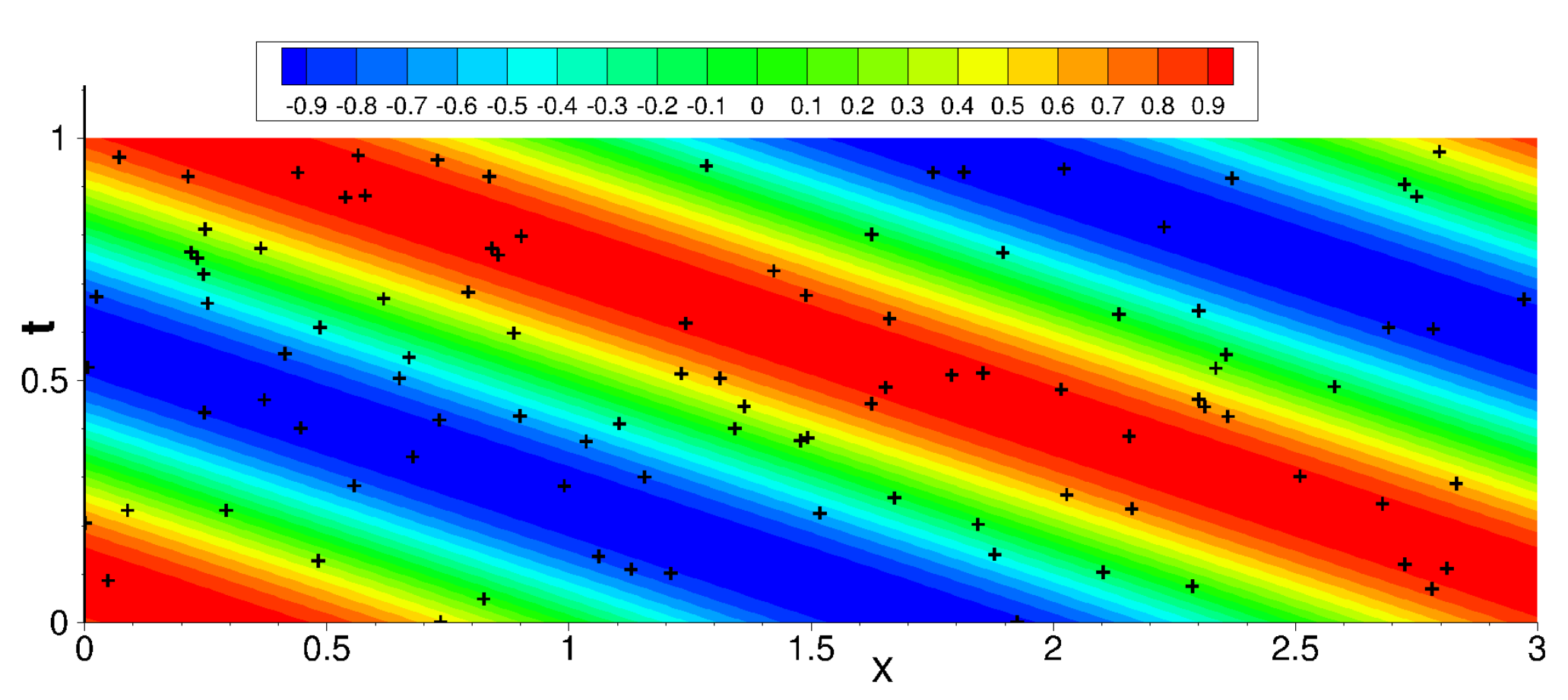}(a)
    \includegraphics[width=2.7in]{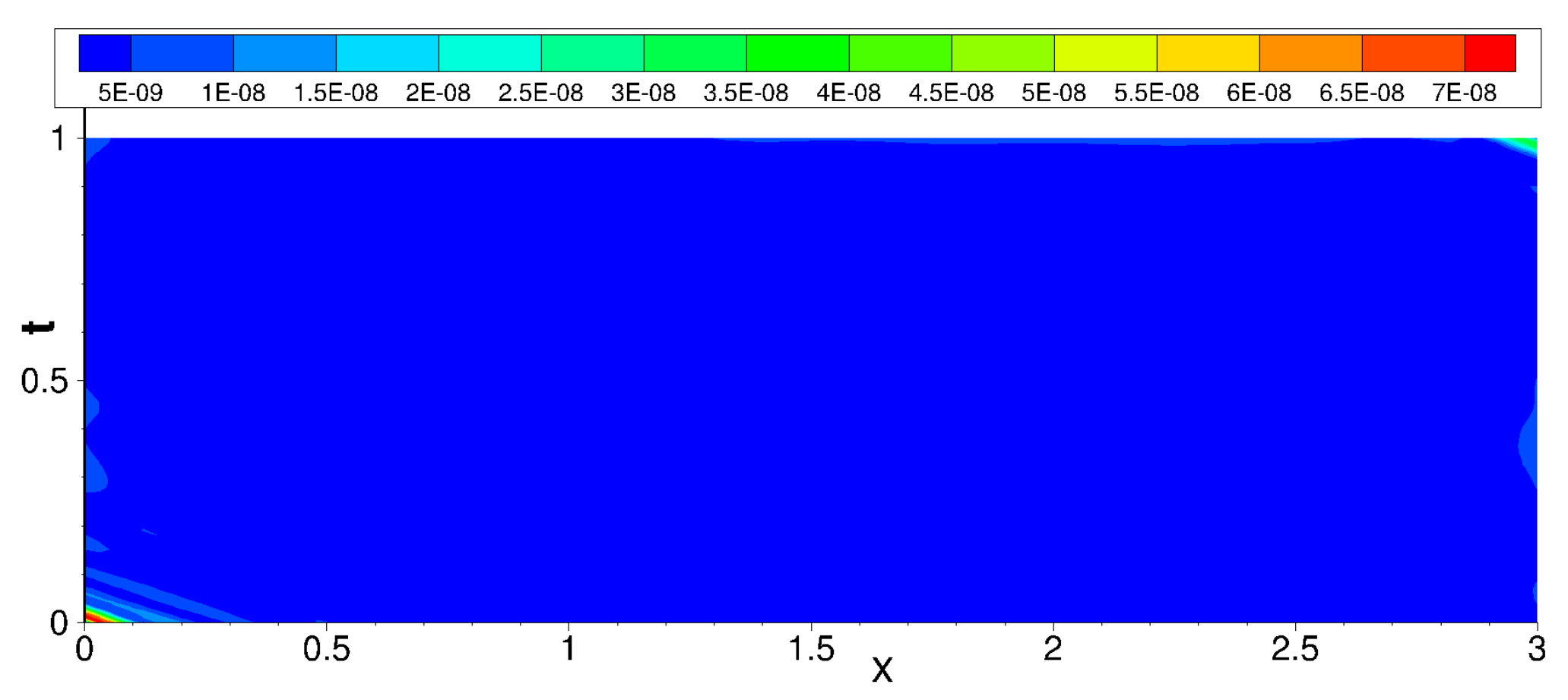}(b)
  }
  \caption{\small Inverse advection problem: distributions of (a) the NLLSQ solution for $u(x,t)$
    and (b) its point-wise
    absolute error,  with the random measurement points shown in (a) as ``+'' symbols.
    Single sub-domain, NN $[2,400,1]$,
    $Q=25\times 25$ (collocation points),  $Q_s=100$ (measurement points),
    $R_m=2.5$, $\lambda_{mea}$=1, $\epsilon=0$ (no noise in measurement).
  }
  \label{fg_4}
\end{figure}

Figure~\ref{fg_4} illustrates the distributions of the NLLSQ solution for $u(x,t)$
and its point-wise absolute error in $\Omega$.
The crucial simulation parameters are listed in the figure caption.
The solution is highly accurate, with a maximum error on the level $10^{-8}$ in the domain.
The computed wave speed $c$ has a relative error $2.82\times 10^{-10}$ for this case.

\begin{table}[tb]
  \small\centering
  \begin{tabular}{llll}
    \hline
    $Q$ & $c$ (NLLSQ) & $c$ (VarPro-F1) & $c$ (VarPro-F2) \\
    5$\times$5 & 3.000074167561E+0 & 2.999935510214E+0 & 6.785575335360E-1  \\
    10$\times$10 & 2.999998340831E+0 & 3.000000635012E+0 & 6.785578125741E-1  \\
    15$\times$15 & 2.999999999982E+0 & 2.999999999967E+0 & -7.284017530389E-2  \\
    20$\times$20 & 3.000000000029E+0 & 3.000000000041E+0 & 3.000000000378E+0  \\
    25$\times$25 & 3.000000000845E+0 & 3.000000000869E+0 & 3.000000000025E+0  \\
    30$\times$30 & 3.000000000534E+0 & 3.000000000542E+0 & 3.000000001047E+0  \\
    35$\times$35 & 3.000000000596E+0 & 3.000000000596E+0 & 3.000000001295E+0  \\
    40$\times$40 & 3.000000000771E+0 & 3.000000000770E+0 & 3.000000001534E+0  \\
    \hline
  \end{tabular}
  \caption{\small Inverse advection problem: the computed $c$ versus $Q$
    obtained by the NLLSQ, VarPro-F1 and VarPro-F2 algorithms.
    Single sub-domain, NN $[2,400,1]$, $Q_s=100$;
    $R_m=2.5$ with NLLSQ and VarPro-F1, and $R_m=2.0$ with VarPro-F2;
    $\lambda_{mea}$=1, $\epsilon=0$.
  }
  \label{tab_7}
\end{table}


\begin{table}[tb]
  \small\centering
  \begin{tabular}{l|ccc|ccc|ccc}
    \hline
     & & NLLSQ & & & VarPro-F1 & & & VarPro-F2 & \\ \cline{2-10}
    $Q$ & $e_c$ & $l^{\infty}$-u & $l^2$-u & $e_c$ & $l^{\infty}$-u & $l^2$-u 
    & $e_c$ & $l^{\infty}$-u & $l^2$-u  \\ \hline
    5$\times$5 & 2.47E-5 & 4.56E-2 & 4.60E-3 & 2.15E-5 & 2.03E-2 & 2.19E-3 & 7.74E-1 & 1.51E+2 & 1.46E+1 \\
    10$\times$10 & 5.53E-7 & 5.25E-4 & 5.70E-5 & 2.12E-7 & 3.16E-4 & 2.56E-5 & 7.74E-1 & 1.66E+3 & 1.11E+2 \\
    15$\times$15 & 5.91E-12 & 8.55E-6 & 3.50E-7 & 1.12E-11 & 8.76E-6 & 3.57E-7 & 1.02E+0 & 1.83E+4 & 5.68E+2 \\
    20$\times$20 & 9.74E-12 & 2.28E-7 & 8.31E-9 & 1.37E-11 & 2.31E-7 & 8.14E-9 & 1.26E-10 & 4.38E-7 & 1.42E-8 \\
    25$\times$25 & 2.82E-10 & 1.08E-7 & 3.50E-9 & 2.90E-10 & 1.09E-7 & 3.49E-9 & 8.47E-12 & 1.67E-7 & 5.47E-9 \\
    30$\times$30 & 1.78E-10 & 8.29E-8 & 3.23E-9 & 1.81E-10 & 8.29E-8 & 3.24E-9 & 3.49E-10 & 9.89E-8 & 5.13E-9 \\
    35$\times$35 & 1.99E-10 & 7.30E-8 & 3.95E-9 & 1.99E-10 & 7.29E-8 & 3.95E-9 & 4.32E-10 & 6.27E-8 & 5.39E-9 \\
    40$\times$40 & 2.57E-10 & 6.58E-8 & 4.67E-9 & 2.57E-10 & 6.57E-8 & 4.67E-9 & 5.11E-10 & 4.87E-8 & 5.77E-9 \\ 
    \hline
  \end{tabular}
  \caption{\small Inverse advection problem: $c$ and $u$ relative errors versus $Q$
    obtained by the NLLSQ, VarPro-F1 and VarPro-F2 algorithms.
    Single sub-domain, NN $[2, 400, 1]$,
    $Q_s=100$, $\lambda_{mea}$=1, $\epsilon$=0; $R_m=2.5$ with NLLSQ and VarPro-F1, and $R_m=2.0$ with VarPro-F2.
    $e_c$ denotes the relative error of $c$.
  }
  \label{tab_8}
\end{table}

\begin{table}[tb]
  \small\centering
  \begin{tabular}{l|ccc|ccc|ccc}
    \hline
     & & NLLSQ & & & VarPro-F1 & & & VarPro-F2 & \\ \cline{2-10}
    $M$ & $e_c$ & $l^{\infty}$-u & $l^2$-u & $e_c$ & $l^{\infty}$-u & $l^2$-u 
    & $e_c$ & $l^{\infty}$-u & $l^2$-u  \\ \hline
    50 & 1.97E-2 & 1.69E-1 & 5.53E-2 & 1.97E-2 & 1.69E-1 & 5.53E-2 & 5.56E-3 & 9.33E-2 & 2.66E-2 \\
    100 & 8.79E-4 & 3.71E-2 & 2.18E-2 & 8.79E-4 & 3.71E-2 & 2.18E-2 & 8.57E-6 & 4.25E-3 & 1.60E-3 \\
    200 & 1.66E-6 & 1.10E-4 & 2.03E-5 & 1.66E-6 & 1.10E-4 & 2.03E-5 & 1.80E-6 & 8.04E-5 & 2.12E-5 \\
    300 & 9.84E-9 & 4.12E-6 & 4.96E-7 & 9.84E-9 & 4.12E-6 & 4.96E-7 & 6.65E-9 & 1.21E-6 & 8.44E-8 \\
    400 & 8.05E-11 & 1.11E-7 & 3.96E-9 & 6.54E-11 & 1.11E-7 & 3.98E-9 & 5.02E-10 & 1.62E-7 & 6.07E-9 \\
    \hline
  \end{tabular}
  \caption{\small Inverse advection problem: $c$ and $u$ relative errors versus $M$ (number of training parameters)
    obtained with the NLLSQ, VarPro-F1 and VarPro-F2 algorithms.
    Single sub-domain, NN $[2, M, 1]$,
    $Q=25\times25$, $Q_s=50$, $\lambda_{mea}$=1, $\epsilon$=0;
    $R_m=2.5$ with NLLSQ and VarPro-F1, and $R_m=2.0$ with VarPro-F2.
  }
  \label{tab_9}
\end{table}

The convergence behaviors of the computed $c$ and $u$ with respect to the collocation points ($Q$)
and to the training parameters ($M$) are illustrated in Tables~\ref{tab_7} to~\ref{tab_9} (without noise).
Table~\ref{tab_7} and Table~\ref{tab_8} list the computed $c$ values, and the relative errors of $c$ errors and $u$,
for several sets of uniform collocation points obtained by the NLLSQ, VarPro-F1 and
VarPro-F2 algorithms.
Table~\ref{tab_9} shows the $c$ errors and the $u$ errors for several sets of training parameters
with the three algorithms.
One can observe the general exponential convergence of the $c$ and $u$ errors with respect to
$Q$ and to $M$.
Tables~\ref{tab_7} and~\ref{tab_8} indicate that the convergence of VarPro-F2 with respect to $Q$
is not quite regular. If the set of collocation points is too small ($Q=15\times 15$ and below),
the computed VarPro-F2 results are not accurate.


\begin{table}[tb]
  \small\centering
  \begin{tabular}{ll | ll | ll}
    \hline
    $\epsilon$ & computed $c$ & $\epsilon$ & computed $c$ &
    $\epsilon$ & computed $c$ \\ \hline
    0.0 & 3.000000000534E+0  & 0.01 & 2.9997368E+0  & 0.1 & 2.9971795E+0  \\
    0.001 & 2.9999739E+0  & 0.03 & 2.9991975E+0  & 0.2 & 2.9937846E+0  \\
    0.002 & 2.9999477E+0  & 0.05 & 2.9986396E+0  & 0.5 & 2.9779459E+0  \\
    0.005 & 2.9998688E+0  & 0.07 & 2.9980677E+0  & 0.7 & 2.9570522E+0  \\
    0.007 & 2.9998158E+0  & 0.09 & 2.9974839E+0  & 1.0 & 2.8441808E+0 \\
    \hline
  \end{tabular}
  \caption{\small Inverse advection problem:  $c$  computed by
    the NLLSQ algorithm corresponding to several noise levels $\epsilon$.
    Single sub-domain, NN $[2, 400, 1]$,
    $Q=30\times 30$, $Q_s=100$, $R_m=2.5$, $\lambda_{mea}$=1. 
  }
  \label{tab_10}
\end{table}

\begin{table}[tb]
  \small\centering
  \begin{tabular}{l|ccc|ccc|ccc}
    \hline
     & & NLLSQ & & & VarPro-F1 & & & VarPro-F2 & \\ \cline{2-10}
    $\epsilon$ & $e_c$ & $l^{\infty}$-u & $l^2$-u & $e_c$ & $l^{\infty}$-u & $l^2$-u 
    & $e_c$ & $l^{\infty}$-u & $l^2$-u  \\ \hline
    0.0 & 1.78E-10 & 8.29E-8 & 3.23E-9 & 1.81E-10 & 8.29E-8 & 3.24E-9 & 3.49E-10 & 9.89E-8 & 5.13E-9 \\
    0.001 & 8.72E-6 & 4.91E-4 & 1.77E-4 & 8.73E-6 & 4.91E-4 & 1.77E-4 & 4.67E-6 & 6.91E-4 & 1.76E-4 \\
    0.005 & 4.37E-5 & 2.46E-3 & 8.85E-4 & 4.41E-5 & 2.46E-3 & 8.84E-4 & 2.26E-5 & 3.49E-3 & 8.80E-4 \\
    0.01 & 8.77E-5 & 4.91E-3 & 1.77E-3 & 8.81E-5 & 4.91E-3 & 1.77E-3 & 4.61E-5 & 7.00E-3 & 1.76E-3 \\
    0.05 & 4.53E-4 & 2.46E-2 & 8.84E-3 & 4.55E-4 & 2.45E-2 & 8.84E-3 & 2.49E-4 & 3.50E-2 & 8.79E-3 \\
    0.1 & 9.40E-4 & 4.92E-2 & 1.77E-2 & 9.54E-4 & 4.92E-2 & 1.77E-2 & 5.53E-4 & 6.95E-2 & 1.76E-2 \\
    0.5 & 7.35E-3 & 2.47E-1 & 8.90E-2 & 7.40E-3 & 2.47E-1 & 8.90E-2 & 5.42E-3 & 3.59E-1 & 8.86E-2 \\
    1.0 & 5.19E-2 & 5.04E-1 & 2.06E-1 & 5.19E-2 & 5.04E-1 & 2.06E-1 & 3.78E-2 & 8.03E-1 & 1.95E-1
    \\ \hline
  \end{tabular}
  \caption{\small Inverse advection problem: $c$ and $u$ relative errors versus $\epsilon$
    obtained with the NLLSQ, VarPro-F1 and VarPro-F2 algorithms.
    Single sub-domain, NN $[2, 400, 1]$,
    $Q=30\times 30$, $Q_s=100$, $\lambda_{mea}$=1;
    $R_m=2.5$ with NLLSQ and VarPro-F1, and $R_m=2.0$ with VarPro-F2.
  }
  \label{tab_11}
\end{table}

\begin{table}[tb]
  \small\centering
  \begin{tabular}{l|ccc|ccc}
    \hline
     & & $\lambda_{mea}$=0.25 & & & $\lambda_{mea}$=0.1 &  \\ \cline{2-7}
    $\epsilon$ & $e_c$ & $l^{\infty}$-u & $l^2$-u & $e_c$ & $l^{\infty}$-u & $l^2$-u  \\ \hline
    0.0 & 2.31E-10 & 8.53E-8 & 4.31E-9 & 2.32E-10 & 8.51E-8 & 4.66E-9 \\
    0.001 & 2.71E-6 & 6.68E-5 & 3.19E-5 & 2.60E-6 & 2.84E-5 & 1.10E-5 \\
    0.005 & 1.36E-5 & 3.34E-4 & 1.59E-4 & 1.30E-5 & 1.42E-4 & 5.49E-5 \\
    0.01 & 2.71E-5 & 6.68E-4 & 3.19E-4 & 2.61E-5 & 2.85E-4 & 1.10E-4 \\
    0.05 & 1.37E-4 & 3.34E-3 & 1.59E-3 & 1.31E-4 & 1.43E-3 & 5.52E-4 \\
    0.1 & 2.75E-4 & 6.69E-3 & 3.19E-3 & 2.65E-4 & 2.88E-3 & 1.11E-3 \\
    0.5 & 1.46E-3 & 3.35E-2 & 1.60E-2 & 1.44E-3 & 1.54E-2 & 5.89E-3 \\
    1.0 & 3.25E-3 & 6.74E-2 & 3.21E-2 & 3.25E-3 & 3.38E-2 & 1.29E-2
    \\ \hline
  \end{tabular}
  \caption{\small Inverse advection problem: $c$ and $u$ relative errors versus $\epsilon$ and
    $\lambda_{mea}$ (scaling coefficient of measurement residual)
    obtained with the NLLSQ algorithm.
    Single sub-domain, NN $[2, 400, 1]$,
    $Q=30\times 30$, $Q_s=100$,
    $R_m=2.5$. These can be compared with the NLLSQ data in Table~\ref{tab_11} for $\lambda_{mea}$=1.
  }
  \label{tab_011}
\end{table}

The effects of noisy measurement data on the computation accuracy
are illustrated by Tables~\ref{tab_10} through~\ref{tab_011}.
Table~\ref{tab_10} lists the computed $c$ by the NLLSQ algorithm
corresponding to several noise levels $\epsilon$ in the measurement data.
Table~\ref{tab_11} shows the $c$ and $u$ relative errors corresponding
to different noise levels, computed by the NLLSQ, VarPro-F1 and VarPro-F2
algorithms. Table~\ref{tab_011} shows the relative errors for $c$ and $u$ corresponding to
different noise levels and several $\lambda_{mea}$ values, illustrating
the effect of scaling the measurement residual (see Remark~\ref{rem_04}).
The computation results are observed to be quite robust to
the noise in the measurement. For example, with $10\%$ noise ($\epsilon=0.1$)
in the measurement, the relative errors of $c$ computed by these
methods are generally on the level of $0.1\%$ (see Table~\ref{tab_11}).
Scaling the measurement residual by $\lambda_{mea}<1$ can markedly improve
the simulation accuracy in the presence of noise, but slightly degrades
the accuracy for the noise-free data; see Table~\ref{tab_011}.

\subsection{Parametric Nonlinear Helmholtz Equation}
\label{nonl_helm}

Consider the 2D domain, $(x,y)\in\Omega=[0,1.4]\times[0,1.4]$, and the
inverse problem on $\Omega$,
\begin{subequations}
  \begin{align}\small
    &
    \frac{\partial^2u}{\partial x^2} + \frac{\partial^2u}{\partial y^2}
    - \alpha_1 u + \alpha_2 \cos(2u) = f(x,y), \\
    &
    u(0,y) = g_1(y), \quad
    u(1.4,y) = g_2(y), \quad
    u(x,0) = g_3(x), \quad
    u(x,1.4) = g_4(x), \\
    &
    u(\xi_i,\eta_i) = S(\xi_i,\eta_i), \quad (\xi_i,\eta_i)\in\mbb Y\subset\Omega,
    \ 1\leqslant i\leqslant NQ_s,
  \end{align}
\end{subequations}
where $f$ and $g_i$ ($1\leqslant i\leqslant 4$) are prescribed source term and boundary data,
$\mbb Y$ denotes the set of random measurement points in $\Omega$,
and the parameters $(\alpha_1,\alpha_2)$ and the field $u(x,y)$ are the unknowns to be determined.
We consider the following manufactured solution to this problem in the tests,
\begin{equation}\label{eq_57}\small
  \begin{split}
    &
    \alpha_1^{ex} = 100, \quad \alpha_2^{ex}=5, \quad
    u_{ex}(x,y) = \cos(\pi x^2)\cos(\pi y^2).
  \end{split}
\end{equation}
The measurement data $S(\xi_i,\eta_i)$ ($1\leqslant i\leqslant NQ_s$) are given by~\eqref{eq_53},
in which $u_{ex}$ is given by~\eqref{eq_57}.
The $u$ errors are computed on a set of $101\times 101$ uniform grid points
in each sub-domain after the neural network is trained.
The notations below follow those of the previous sub-sections.

\begin{figure}
  \small
  \centerline{
    \includegraphics[height=1.8in]{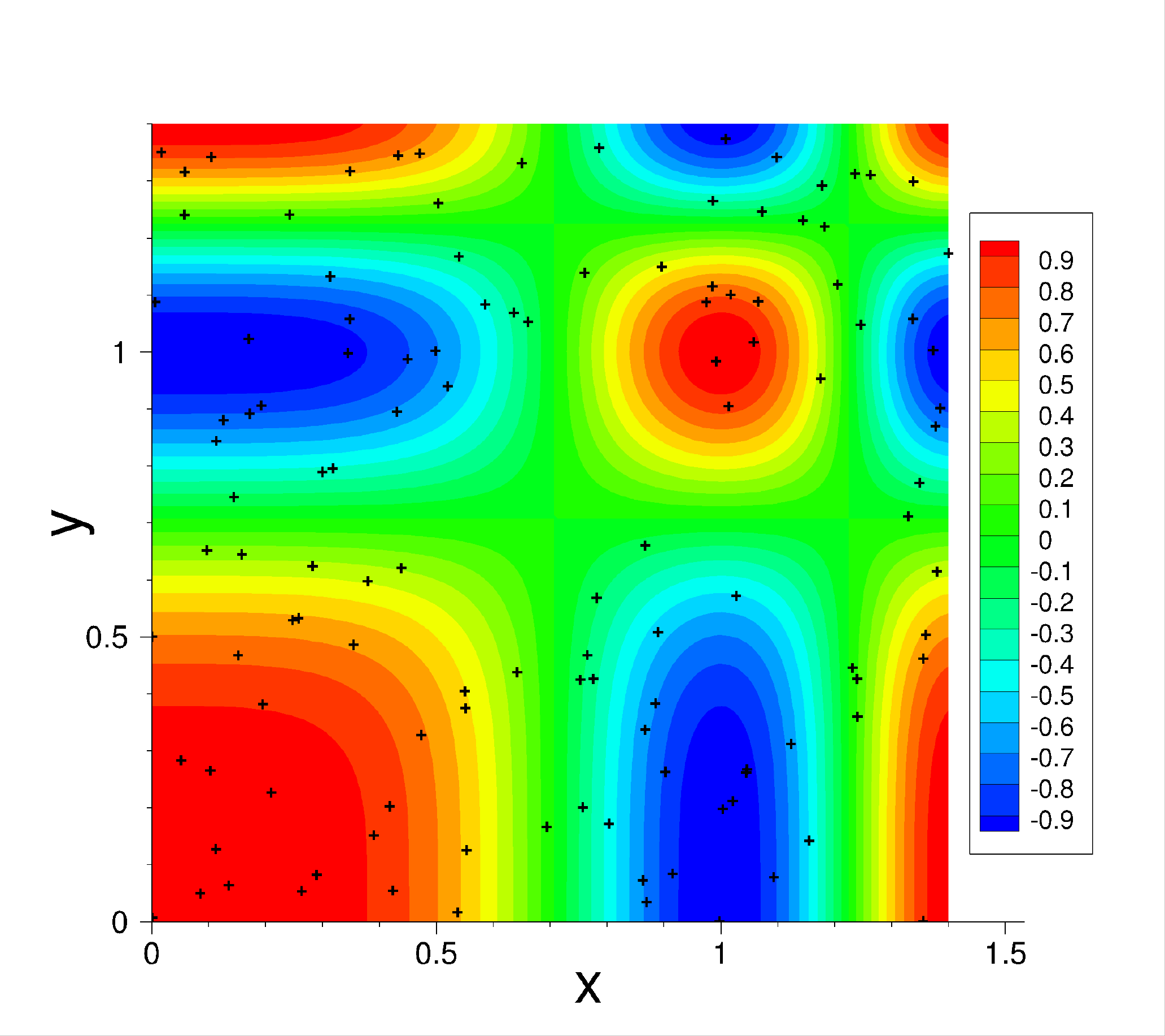}(a)
    \includegraphics[height=1.8in]{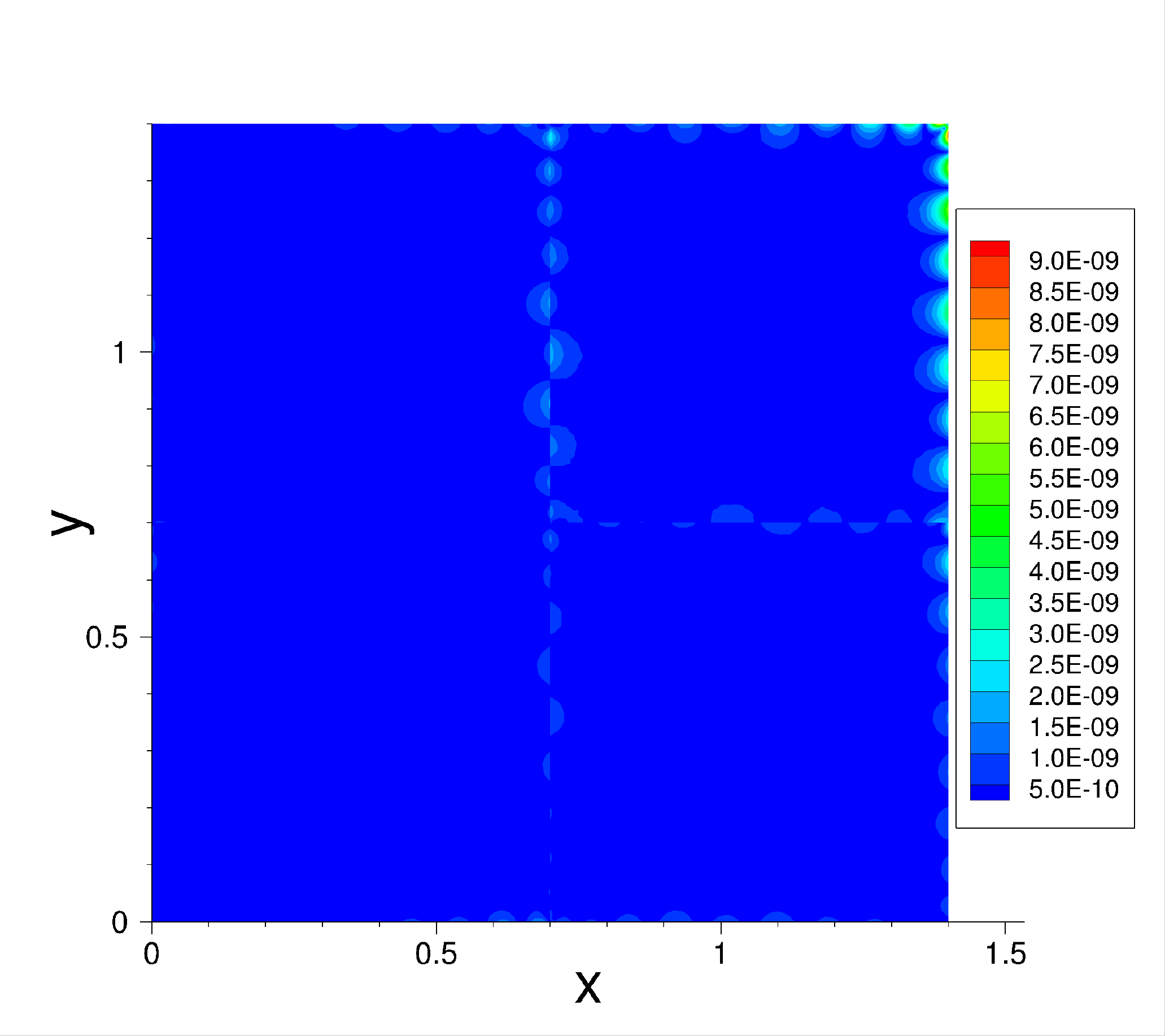}(b)
  }
  \caption{\small Inverse nonlinear Helmholtz problem: distributions of
    (a) the VarPro-F1 solution for $u(x,y)$ and (b) its point-wise absolute error,
    with the random measurement points shown in (a) as ``+'' symbols.
    Four uniform sub-domains ($2$ in each direction), local NN $[2,300,1]$,
    $Q=20\times 20$,
    $Q_s=30$, $R_m=1.5$, $\lambda_{mea}$=1, $\epsilon$=0 (no noise in measurement data).
  }
  \label{fg_5}
\end{figure}

Figure~\ref{fg_5} shows distributions of the $u(x,y)$ solution
and its point-wise absolute error computed by the VarPro-F1 algorithm on $4$ uniform
sub-domains, with
the $120$ random measurement points in total ($Q_s=30$ points per sub-domain)
displayed in Figure~\ref{fg_5}(a).
The figure caption lists the crucial simulation parameters for
this test. VarPro-F1 exhibits a high accuracy, with the maximum
$u$ error on the order of $10^{-8}$.
The relative errors of the computed $\alpha_1$ and $\alpha_2$ are
$3.52\times 10^{-11}$ and $4.76\times 10^{-10}$, respectively, in this test.


\begin{table}[tb]
  \small\centering
  \begin{tabular}{lcc}
    \hline
    $Q$ & computed $\alpha_1$ & computed $\alpha_2$  \\ 
    5$\times$5 & 9.946591149073E+1 & 5.169760481373E+0  \\
    10$\times$10 & 9.999987125506E+1 & 4.999987933629E+0  \\
    15$\times$15 & 9.999999986638E+1 & 5.000000027512E+0  \\
    20$\times$20 & 9.999999982078E+1 & 4.999999813483E+0  \\
    25$\times$25 & 1.000000001774E+2 & 4.999999946859E+0  \\
    30$\times$30 & 1.000000001832E+2 & 4.999999843159E+0  \\
    35$\times$35 & 9.999999989070E+1 & 5.000000059829E+0  \\
    40$\times$40 & 9.999999957958E+1 & 5.000001280912E+0  \\
    \hline
  \end{tabular}
  \caption{\small Inverse nonlinear Helmholtz problem: $\alpha_1$ and $\alpha_2$ versus $Q$
    (number of collocation points) obtained by the NLLSQ algorithm.
    Single sub-domain, NN $[2,500,1]$,  $Q_s=100$,
    $R_m=2.25$, $\epsilon=0$, $\lambda_{mea}$=1.
  }
  \label{tab_12}
\end{table}

\begin{table}[tb]
  \small\centering
  \begin{tabular}{l|ccc|ccc|ccc}
    \hline
     & &  NLLSQ & & & VarPro-F1 & & & VarPro-F2 & \\ \cline{2-10}
    $Q$ & $e_{\alpha_1}$ & $e_{\alpha_2}$ & $l^2$-u &  $e_{\alpha_1}$ & $e_{\alpha_2}$ & $l^2$-u &
       $e_{\alpha_1}$ & $e_{\alpha_2}$ & $l^2$-u \\ \hline
    5$\times$5 & 5.34E-3 & 3.40E-2 & 2.26E-2 & 2.35E-4 & 5.72E-4 & 4.16E-3 & 1.00E+0 & 1.00E+0 & 1.56E+0 \\
    10$\times$10 & 1.29E-6 & 2.41E-6 & 2.21E-4 & 1.29E-6 & 2.42E-6 & 2.20E-4 & 1.00E+0 & 1.00E+0 & 1.66E+1 \\
    15$\times$15 & 1.34E-9 & 5.50E-9 & 6.55E-7 & 8.81E-11 & 1.60E-9 & 5.64E-7 & 5.41E-1 & 3.72E-1 & 4.84E+1 \\
    20$\times$20 & 1.79E-9 & 3.73E-8 & 4.93E-8 & 7.48E-10 & 5.74E-8 & 5.67E-8 & 1.28E-8 & 4.63E-8 & 1.65E-7 \\
    25$\times$25 & 1.77E-9 & 1.06E-8 & 8.52E-9 & 1.11E-9 & 3.30E-9 & 7.12E-9 & 1.43E-8 & 2.85E-7 & 4.30E-8 \\
    30$\times$30 & 1.83E-9 & 3.14E-8 & 1.00E-8 & 5.46E-10 & 4.46E-9 & 8.39E-9 & 5.45E-9 & 2.24E-7 & 6.02E-8 \\
    35$\times$35 & 1.09E-9 & 1.20E-8 & 9.42E-9 & 2.58E-9 & 1.79E-7 & 1.18E-8 & 1.23E-8 & 3.86E-7 & 7.50E-8 \\
    40$\times$40 & 4.20E-9 & 2.56E-7 & 1.54E-8 & 5.14E-9 & 2.73E-7 & 1.58E-8 & 1.22E-8 & 1.85E-7 & 9.12E-8 \\
    \hline
  \end{tabular}
  \caption{\small Inverse nonlinear Helmholtz problem: relative errors of
    $\alpha_1$, $\alpha_2$ and $u$ versus~$Q$
    obtained by the NLLSQ, VarPro-F1 and VarPro-F2 algorithms.
    Single sub-domain, NN $[2, 500, 1]$,
    $Q_s=100$; $R_m=2.25$ with NLLSQ and VarPro-F1, and $R_m=2.5$ with
    VarPro-F2; $\lambda_{mea}$=1, $\epsilon$=0.
    $e_{\alpha_1}$ and $e_{\alpha_2}$ denote the relative errors
    of $\alpha_1$ and $\alpha_2$.
  }
  \label{tab_13}
\end{table}

The convergence of the simulation results with respect to
the number of collocation points ($Q$) is illustrated by
the Tables~\ref{tab_12} and~\ref{tab_13}.
Table~\ref{tab_12} lists the computed $\alpha_1$ and $\alpha_2$
by the NLLSQ algorithm corresponding to a range of $Q$
values. Table~\ref{tab_13} shows the relative $\alpha_1$ and $\alpha_2$
errors and the $l^2$ norm of the relative $u$ error
corresponding to different $Q$ obtained by the NLLSQ, VarPro-F1
and VarPro-F2 algorithms. The crucial simulation parameter values
are provided in the table captions.
A general exponential convergence in the errors with respect to
$Q$ can be observed. One can also observe that the convergence of the
VarPro-F2 algorithm appears to be less regular.
The VarPro-F2 results are inaccurate with a small $Q$ ($Q=15\times 15$ or less),
and its errors abruptly drop to $10^{-7}\sim 10^{-8}$ as the collocation
points reach $Q=20\times 20$ and beyond.

\begin{table}[tb]
  \small\centering
  \begin{tabular}{l|ccc|ccc|ccc}
    \hline
     & &  NLLSQ & & & VarPro-F1 & & & VarPro-F2 & \\ \cline{2-10}
    $M$ & $e_{\alpha_1}$ & $e_{\alpha_2}$ & $l^2$-u &  $e_{\alpha_1}$ & $e_{\alpha_2}$ & $l^2$-u &
       $e_{\alpha_1}$ & $e_{\alpha_2}$ & $l^2$-u \\ \hline
    100 & 5.58E-1 & 7.37E-1 & 3.21E-1 & 5.58E-1 & 7.37E-1 & 3.21E-1 & 1.08E+1 & 3.61E+0 & 8.80E-1 \\
    200 & 4.39E-3 & 2.60E-3 & 7.66E-3 & 4.39E-3 & 2.60E-3 & 7.66E-3 & 3.70E-3 & 4.05E-2 & 1.04E-2 \\
    300 & 6.52E-5 & 7.27E-5 & 5.15E-5 & 6.51E-5 & 7.07E-5 & 5.14E-5 & 5.28E-5 & 2.58E-4 & 8.91E-5 \\
    400 & 1.08E-7 & 2.17E-6 & 6.22E-7 & 1.06E-7 & 1.80E-6 & 6.33E-7 & 3.62E-7 & 1.18E-5 & 1.12E-6 \\
    500 & 1.83E-9 & 3.14E-8 & 1.00E-8 & 5.46E-10 & 4.46E-9 & 8.39E-9 & 5.45E-9 & 2.24E-7 & 6.02E-8 \\
    600 & 5.16E-10 & 1.47E-8 & 2.14E-9 & 1.45E-10 & 5.33E-10 & 2.29E-9 & 4.47E-10 & 3.30E-9 & 1.77E-8 \\
    \hline
  \end{tabular}
  \caption{\small Inverse nonlinear Helmholtz problem: relative errors of
    $\alpha_1$, $\alpha_2$ and $u$ versus~$M$
    (number of training parameters)
    obtained by the NLLSQ, VarPro-F1 and VarPro-F2 algorithms.
    Single sub-domain, NN $[2, M, 1]$,
    $Q_s=100$, $Q=30\times 30$;
    $R_m=2.25$ with NLLSQ and VarPro-F1, and $R_m=2.5$ with VarPro-F2;
    $\epsilon$=0, $\lambda_{mea}$=1.
  }
  \label{tab_14}
\end{table}

Table~\ref{tab_14} illustrates the convergence of the
$\alpha_1$, $\alpha_2$ and $u$ errors, obtained by the NLLSQ,
VarPro-F1 and VarPro-F2 algorithms, with respect to the training parameters ($M$).
The table caption lists values of the main simulation parameters.
The relative errors of $\alpha_1$, $\alpha_2$ and $u$ decrease exponentially
with increasing $M$.

\begin{table}[tb]
  \small\centering
  \begin{tabular}{l|cc|cc|cc}
    \hline
    $Q_s$ & $\alpha_1$ & $\alpha_2$ & $e_{\alpha_1}$ & $e_{\alpha_2}$ & $l^{\infty}$-u & $l^2$-u  \\ \hline
    5 & 1.000000005409E+2 & 4.999994839886E+0 & 5.41E-9 & 1.03E-6 & 7.79E-8 & 4.13E-8 \\
    10 & 9.999999992698E+1 & 4.999999378668E+0 & 7.30E-10 & 1.24E-7 & 6.03E-8 & 8.65E-9 \\
    20 & 1.000000001181E+2 & 5.000000028574E+0 & 1.18E-9 & 5.71E-9 & 9.67E-8 & 8.94E-9 \\
    30 & 9.999999993263E+1 & 5.000000391330E+0 & 6.74E-10 & 7.83E-8 & 7.97E-8 & 8.18E-9 \\
    50 & 9.999999980057E+1 & 5.000000659553E+0 & 1.99E-9 & 1.32E-7 & 7.56E-8 & 1.05E-8 \\
    100 & 1.000000001832E+2 & 4.999999843159E+0 & 1.83E-9 & 3.14E-8 & 9.37E-8 & 1.00E-8 \\
    \hline
  \end{tabular}
  \caption{\small Inverse nonlinear Helmholtz problem: the values/relative-errors of
    $\alpha_1$ and $\alpha_2$, and the $u$
    relative errors, versus~the number of random measurement points ($Q_s$),
    computed by the NLLSQ algorithm.
    Single sub-domain, NN $[2, 500, 1]$,
    $Q=30\times 30$, $R_m=2.25$, $\lambda_{mea}$=1, $\epsilon$=0.
  }
  \label{tab_15}
\end{table}

Table~\ref{tab_15} shows the computed $\alpha_1$ and $\alpha_2$ values,
their relative errors, and the $u$ relative errors ($l^{\infty}$ and $l^2$ norms)
obtained by the NLLSQ algorithm corresponding to a range of $Q_s$ (number of
random measurement points).
The effect of $Q_s$ on the errors appears to
be not significant, unless $Q_s$ is very small. This is
similar to what has been observed with linear forward PDEs (see e.g.~Section~\ref{poisson}).


\begin{table}[tb]
  \small\centering
  \begin{tabular}{l|ccc|ccc|ccc}
    \hline
     & &  NLLSQ & & & VarPro-F1 & & & VarPro-F2 & \\ \cline{2-10}
    $\epsilon$ & $e_{\alpha_1}$ & $e_{\alpha_2}$ & $l^2$-u &  $e_{\alpha_1}$ & $e_{\alpha_2}$ & $l^2$-u &
       $e_{\alpha_1}$ & $e_{\alpha_2}$ & $l^2$-u \\ \hline
    0.0 & 1.99E-9 & 1.32E-7 & 1.05E-8 & 1.96E-11 & 1.76E-8 & 7.85E-9 & 7.18E-9 & 4.77E-7 & 5.95E-8 \\
    0.001 & 4.31E-4 & 1.12E-4 & 2.46E-4 & 4.31E-4 & 1.04E-4 & 2.46E-4 & 4.31E-4 & 1.45E-4 & 2.46E-4 \\
    0.002 & 8.62E-4 & 2.19E-4 & 4.92E-4 & 8.62E-4 & 2.40E-4 & 4.92E-4 & 8.62E-4 & 2.60E-4 & 4.92E-4 \\
    0.005 & 2.16E-3 & 5.21E-4 & 1.23E-3 & 2.16E-3 & 5.67E-4 & 1.23E-3 & 2.16E-3 & 6.00E-4 & 1.23E-3 \\
    0.01 & 4.32E-3 & 9.64E-4 & 2.46E-3 & 4.32E-3 & 9.54E-4 & 2.46E-3 & 4.32E-3 & 1.26E-3 & 2.46E-3 \\
    0.02 & 8.67E-3 & 1.41E-3 & 4.92E-3 & 8.67E-3 & 1.61E-3 & 4.92E-3 & 8.67E-3 & 1.84E-3 & 4.92E-3 \\
    0.05 & 2.19E-2 & 3.72E-4 & 1.23E-2 & 2.19E-2 & 6.72E-4 & 1.23E-2 & 2.19E-2 & 2.29E-3 & 1.23E-2 \\
    0.1 & 4.44E-2 & 8.72E-3 & 2.45E-2 & 4.44E-2 & 8.77E-3 & 2.45E-2 & 4.44E-2 & 7.37E-3 & 2.45E-2 \\
    \hline
  \end{tabular}
  \caption{\small Inverse nonlinear Helmholtz problem: relative errors of $\alpha_1$, $\alpha_2$ and $u$
    versus the noise level ($\epsilon$) in the measurement data,
    obtained by the NLLSQ, VarPro-F1 and VarPro-F2 algorithms.
    Single sub-domain, NN $[2, 500, 1]$,
    $Q=30\times 30$, $Q_s=50$, $\lambda_{mea}$=1;
    $R_m=2.25$ with NLLSQ and VarPro-F1, and $R_m=2.5$ with VarPro-F2.
  }
  \label{tab_16}
\end{table}

No noise is considered in the measurement data in the foregoing tests.
Table~\ref{tab_16} illustrates the effect of the noise level ($\epsilon$)
on the accuracy of the computed $\alpha_1$, $\alpha_2$ and $u$ by
the NLLSQ, VarPro-F1 and VarPro-F2 algorithms.
The main parameters for these simulations are listed in
the table caption.
The accuracy of these algorithms appears quite robust to the noise.
For example, with $1\%$ noise ($\epsilon=0.01$) in the measurement data
the relative errors for $\alpha_1$ and $\alpha_2$ obtained by
the three methods are on the order of $0.1\%$, and
with $10\%$ noise ($\epsilon=0.1$) in the measurement data
the relative errors for $\alpha_1$
and $\alpha_2$ are on the order of $1\sim 4\%$.

\subsection{Parametric Viscous Burgers Equation}
\label{burger}

Consider the spatial-temporal domain, $(x,t)\in\Omega=[0,2]\times[0,1.5]$, and the inverse
problem with the parametric Burgers' equation,
\begin{subequations}
  \begin{align}\small
    &
    \frac{\partial u}{\partial t} + \alpha_1 u\frac{\partial u}{\partial x}
    = \alpha_2\frac{\partial^2u}{\partial x^2} + f(x,t), \label{eq_58a} \\
    &
    u(0,t) = g_1(t), \quad
    u(2,t) = g_2(t), \quad
    u(x,0) = h(x), \label{eq_58c} \\
    &
    u(\xi_i,\eta_i) = S(\xi_i,\eta_i), \quad
    (\xi_i,\eta_i)\in\mbb Y\subset\Omega, \
    1\leqslant i\leqslant NQ_s, \label{eq_58d}
  \end{align}
\end{subequations}
where $f$ is a prescribed source term, $g_1$ and $g_2$ are prescribed Dirichlet
boundary data, $h$ is the initial distribution,
the constants $\alpha_i$ ($i=1,2$) and the field $u(x,t)$ are the unknowns to be solved for,
$\mbb Y$ denotes the set of random measurement points, $N$ is the number of sub-domains, and
$Q_s$ is the number of measurement points per sub-domain. We employ the following
manufactured solution  in the tests,
\begin{equation}\label{eq_59}\small
  \left\{\small
  \begin{split}
    &
    \alpha_1^{ex} = 0.1, \quad \alpha_2^{ex}=0.01, \\
    &
    u_{ex}(x,t) = \left(1+\frac{x}{20} \right)\left(1+\frac{t}{20} \right)\left[
      \frac32\cos\left(\pi x+\frac{7\pi}{20} \right)
      +\frac{27}{20}\cos\left(2\pi x-\frac{3\pi}{5} \right)
      \right] \left[
      \frac32\cos\left(\pi t+\frac{7\pi}{20} \right) \right. \\
      & \qquad\qquad\ \
      \left.
      +\frac{27}{20}\cos\left(2\pi t-\frac{3\pi}{5} \right)
      \right].
  \end{split}
  \right.
\end{equation}
The source term $f$ and the boundary/initial data are chosen such that
the expressions in~\eqref{eq_59} satisfy the equations~\eqref{eq_58a}--\eqref{eq_58c}.
The measurement data $S(\xi_i,\eta_i)$ is assumed to be given by~\eqref{eq_53},
in which $u_{ex}$ is given by~\eqref{eq_59}.
In the following the $u$ errors are computed on
a $101\times 101$ uniform grid points in each sub-domain,
and we adopt the same notations (e.g.~$Q$, $Q_s$, $M$, $R_m$ and $\epsilon$)
as in previous sub-sections.

\begin{figure}
  \small
  \centerline{
    \includegraphics[width=2.2in]{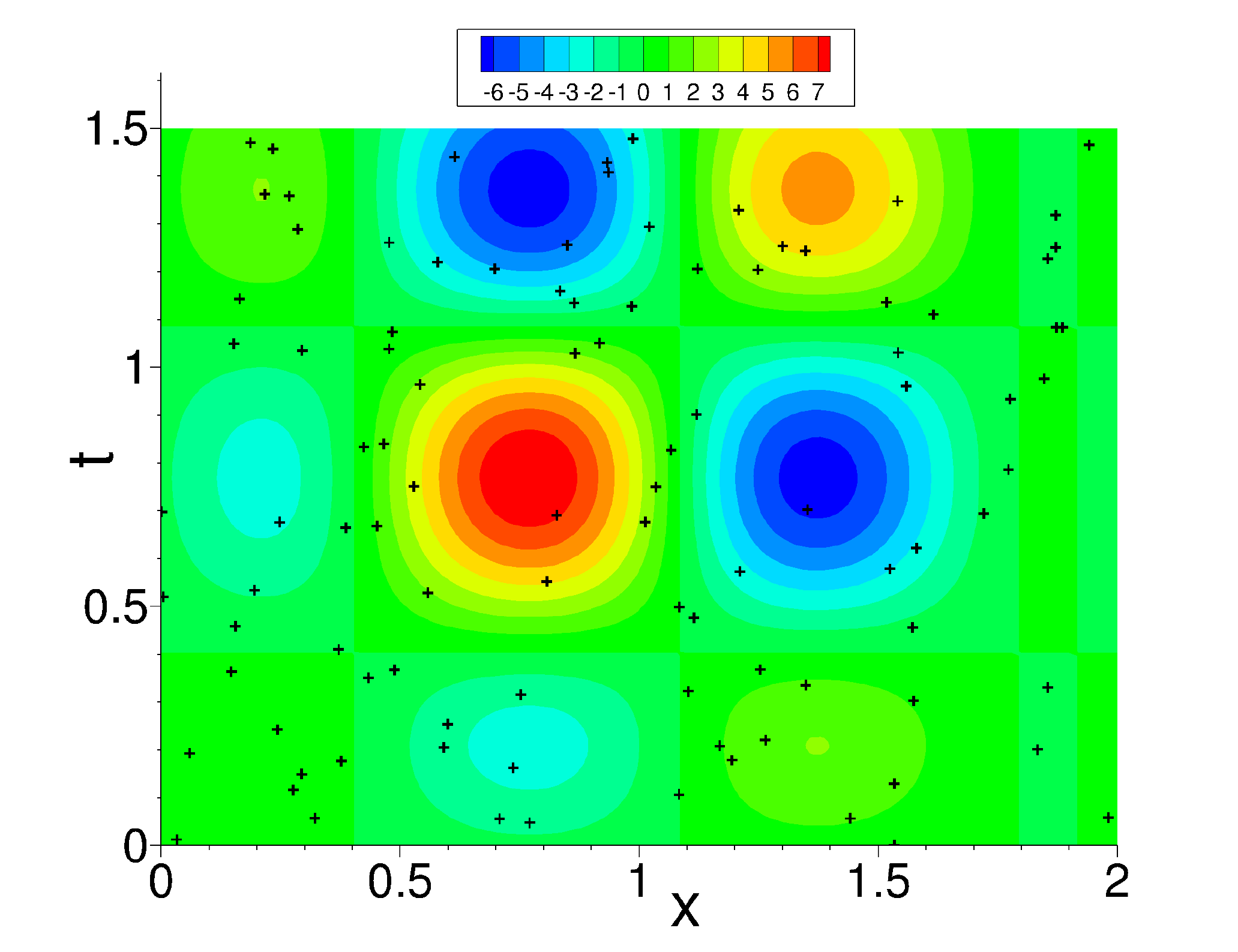}(a)
    \includegraphics[width=2.2in]{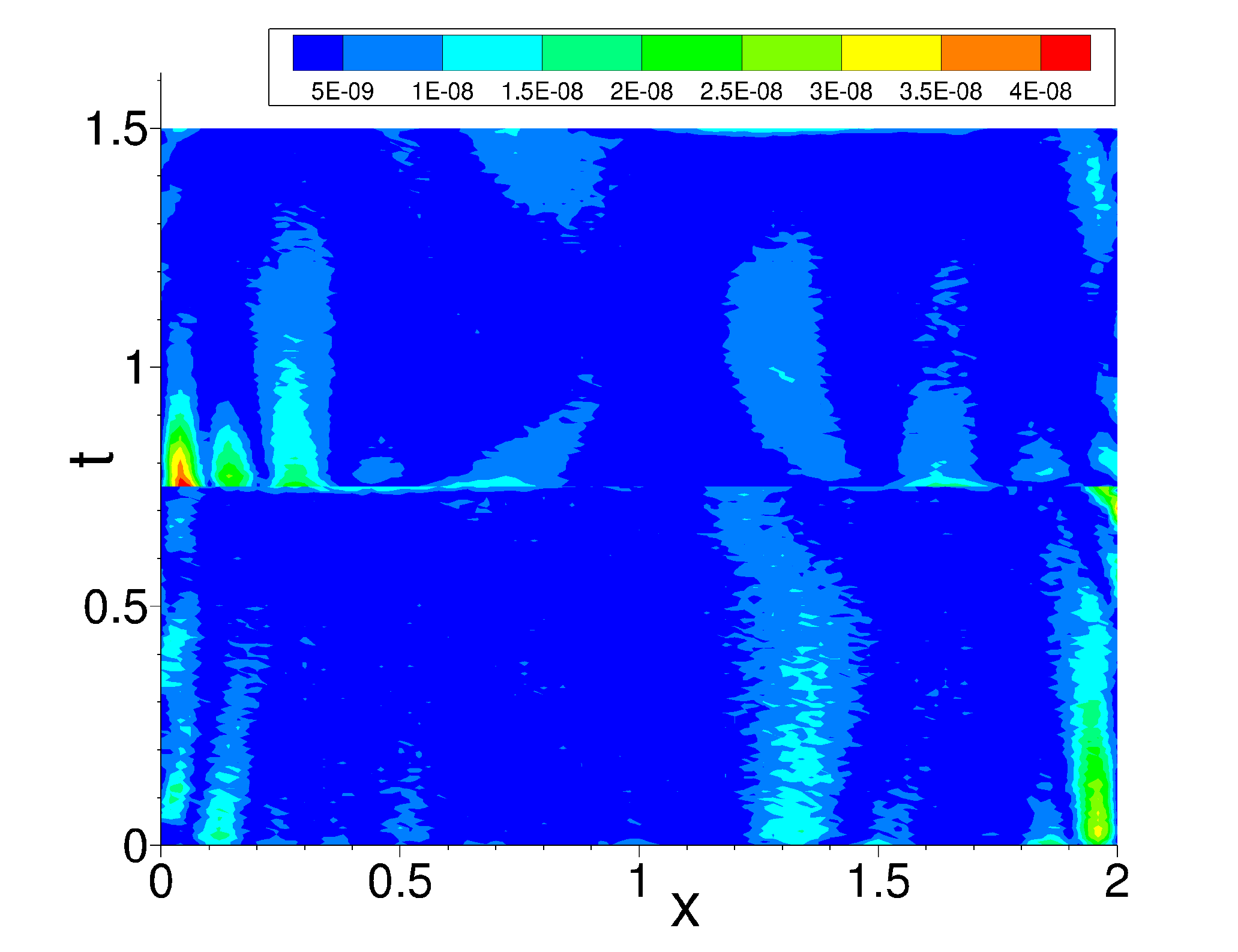}(b)
  }
  \caption{\small Inverse Burgers' problem: distributions of (a) the NLLSQ solution and (b) its
    point-wise absolute error, with the random measurement
    points shown by the ``+'' symbols in (a). Two uniform sub-domains (along $t$),
    local NN $[2,300,1]$, 
    $Q=25\times 25$, $Q_s=50$,
    $R_m=1.5$, $\epsilon$=0 (no noise in measurement data), $\lambda_{mea}$=1.
  }
  \label{fg_6}
\end{figure}

Figure~\ref{fg_6} illustrates the $u(x,t)$ solution and its point-wise
absolute error computed by the NLLSQ algorithm with two uniform sub-domains along $t$,
and the $100$ random measurement points ($50$ points per sub-domain)
in the domain are shown in Figure~\ref{fg_6}(a).
The figure caption provides the main parameter values in this simulation.
The results signify a high accuracy for the computed $u$ solution,
with the maximum error on the order of $10^{-8}$.
The relative errors of the computed $\alpha_1$ and $\alpha_2$ are
$1.30\times 10^{-9}$ and $1.48\times 10^{-8}$, respectively, for this simulation.

\begin{table}[tb]
  \small\centering
  \begin{tabular}{lcc}
    \hline
    $Q$ & $\alpha_1$ & $\alpha_2$  \\ 
    5$\times$5 & 9.999660775275E-2 & 9.994514983290E-3  \\
    10$\times$10 & 9.999998874379E-2 & 9.999992607237E-3  \\
    15$\times$15 & 1.000000000074E-1 & 1.000000000018E-2  \\
    20$\times$20 & 1.000000000060E-1 & 1.000000000487E-2  \\
    25$\times$25 & 9.999999999967E-2 & 9.999999999698E-3  \\
    30$\times$30 & 1.000000000049E-1 & 9.999999998052E-3  \\
    \hline
  \end{tabular}
  \caption{\small Inverse Burgers' problem: the computed $\alpha_1$ and $\alpha_2$ versus $Q$
    (number of collocation points) obtained with the NLLSQ algorithm.
    Single sub-domain, NN $[2,400,1]$, $Q_s=100$,
    $R_m=1.9$, $\lambda_{mea}$=1, $\epsilon$=0.
  }
  \label{tab_17}
\end{table}

\begin{table}[tb]
  \small\centering
  \begin{tabular}{l|ccc|ccc|ccc}
    \hline
     & &  NLLSQ & & & VarPro-F1 & & & VarPro-F2 &  \\ \cline{2-10}
    $Q$ & $e_{\alpha_1}$ & $e_{\alpha_2}$ & $l^2$-u &  $e_{\alpha_1}$ & $e_{\alpha_2}$ & $l^2$-u
    & $e_{\alpha_1}$ & $e_{\alpha_2}$ & $l^2$-u  \\ \hline
    5$\times$5 & 3.39E-5 & 5.49E-4 & 8.77E-4 & 3.75E-5 & 5.96E-4 & 8.64E-4 & 3.04E-1 & 5.71E-1 & 6.24E+0 \\
    10$\times$10 & 1.13E-7 & 7.39E-7 & 8.89E-6 & 1.09E-7 & 6.84E-7 & 8.57E-6 & 9.79E-1 & 1.66E+1 & 7.06E+1 \\
    15$\times$15 & 7.40E-11 & 1.76E-11 & 5.41E-8 & 1.32E-11 & 7.76E-12 & 1.99E-8 & 1.54E-3 & 1.13E-2 & 7.21E-1 \\
    20$\times$20 & 5.97E-11 & 4.87E-10 & 2.03E-9 & 8.42E-11 & 4.02E-10 & 1.42E-9 & 1.53E-10 & 1.06E-9 & 3.73E-9 \\
    25$\times$25 & 3.32E-12 & 3.02E-11 & 6.93E-10 & 5.54E-11 & 1.38E-10 & 4.51E-10 & 1.40E-10 & 3.04E-10 & 1.42E-9 \\
    30$\times$30 & 4.86E-11 & 1.95E-10 & 3.01E-10 & 1.52E-11 & 1.72E-10 & 2.36E-10 & 2.59E-10 & 1.16E-9 & 8.94E-10 \\
    \hline
  \end{tabular}
  \caption{\small Inverse Burgers' problem: relative errors of $\alpha_1$, $\alpha_2$ and $u$ versus $Q$
    obtained by the NLLSQ, VarPro-F1 and VarPro-F2 algorithms.
    Single sub-domain, NN $[2, 400, 1]$,
    $Q_s=100$; $R_m=1.9$ with NLLSQ and VarPro-F1, and $R_m=2.0$ with VarPro-F2;
    $\epsilon$=0, $\lambda_{mea}$=1.
  }
  \label{tab_18}
\end{table}

Tables~\ref{tab_17} and~\ref{tab_18} illustrate the convergence behavior of
the NLLSQ, VarPro-F1 and VarPro-F2 algorithms with respect to
the number of collocation points ($Q$).
Table~\ref{tab_17} shows the computed $\alpha_1$ and $\alpha_2$ values
by the NLLSQ algorithm for several $Q$.
Table~\ref{tab_18} shows the relative errors of $\alpha_1$ and $\alpha_2$ and
the $l^2$ norm of relative $u$ error corresponding to different $Q$.
We refer the reader to the table captions for the simulation parameter
values. One can observe the familiar exponential convergence
with respect to $Q$ (before stagnation when $Q$ reaches a certain level).

\begin{table}[tb]
  \small\centering
  \begin{tabular}{l|ccc|ccc|ccc}
    \hline
     & &  NLLSQ & & & VarPro-F1 & & & VarPro-F2 & \\ \cline{2-10}
    $M$ & $e_{\alpha_1}$ & $e_{\alpha_2}$ & $l^2$-u &  $e_{\alpha_1}$ & $e_{\alpha_2}$ & $l^2$-u  \\ \hline
    50 & 1.07E-1 & 4.07E+0 & 4.22E-1 & 1.07E-1 & 4.07E+0 & 4.22E-1 & 1.01E-1 & 3.62E+0 & 4.11E-1 \\
    100 & 4.84E-3 & 8.21E-2 & 5.15E-2 & 4.84E-3 & 8.21E-2 & 5.15E-2 & 6.72E-3 & 5.10E-3 & 5.52E-2 \\
    200 & 5.27E-6 & 1.86E-4 & 1.69E-5 & 5.27E-6 & 1.86E-4 & 1.69E-5 & 6.39E-6 & 1.85E-4 & 2.01E-5 \\
    300 & 7.97E-9 & 6.35E-8 & 3.34E-8 & 7.98E-9 & 6.40E-8 & 3.34E-8 & 3.11E-8 & 6.45E-8 & 5.19E-8 \\
    400 & 2.64E-11 & 1.78E-10 & 3.03E-10 & 1.02E-11 & 7.89E-11 & 2.34E-10 & 2.14E-10 & 1.31E-9 & 8.44E-10 \\
    500 & 5.23E-12 & 1.64E-12 & 2.66E-11 & 4.85E-12 & 6.55E-12 & 2.39E-11 & 3.49E-11 & 2.09E-10 & 1.05E-10 \\
    \hline
  \end{tabular}
  \caption{\small Inverse Burgers' problem: relative errors
    of $\alpha_1$, $\alpha_2$ and $u$ versus $M$ (training parameters)
    obtained by the NLLSQ, VarPro-F1 and VarPro-F2 algorithms.
    Single sub-domain, NN $[2, M, 1]$ ($M$ varied),
    $Q_s$=150, $Q=30\times 30$, $\epsilon$=0, $\lambda_{mea}$=1;
    $R_m=1.9$ with NLLSQ and VarPro-F1, and $R_m=2.0$ with VarPro-F2.
  }
  \label{tab_19}
\end{table}

The exponential convergence of the simulation results with respect to
the number of training parameters ($M$) for the
NLLSQ, VarPro-F1 and VarPro-F2 algorithms is illustrated by Table~\ref{tab_19}.
This table shows the relative errors of the computed $\alpha_1$, $\alpha_2$
and $u$ obtained by the three algorithms.
One should again refer to the caption for the main settings and simulation parameters.

\begin{table}[tb]
  \small\centering
  \begin{tabular}{l|ccc|ccc|ccc}
    \hline
     & &  NLLSQ & & & VarPro-F1 & & & VarPro-F2 &  \\ \cline{2-10}
    $\epsilon$ & $e_{\alpha_1}$ & $e_{\alpha_2}$ & $l^2$-u &  $e_{\alpha_1}$ & $e_{\alpha_2}$ & $l^2$-u  \\ \hline
    0.0 & 4.86E-11 & 1.95E-10 & 3.01E-10 & 1.52E-11 & 1.72E-10 & 2.36E-10 & 2.59E-10 & 1.16E-9 & 8.94E-10 \\
    0.001 & 2.58E-5 & 2.20E-4 & 9.16E-5 & 2.64E-5 & 2.17E-4 & 9.40E-5 & 2.56E-5 & 2.34E-4 & 9.03E-5 \\
    0.005 & 1.25E-4 & 1.14E-3 & 4.66E-4 & 1.30E-4 & 1.14E-3 & 4.78E-4 & 1.29E-4 & 1.17E-3 & 4.56E-4 \\
    0.01 & 2.57E-4 & 2.25E-3 & 9.39E-4 & 2.57E-4 & 2.23E-3 & 9.35E-4 & 2.57E-4 & 2.36E-3 & 9.07E-4 \\
    0.05 & 1.36E-3 & 1.23E-2 & 4.69E-3 & 1.35E-3 & 1.22E-2 & 4.76E-3 & 1.31E-3 & 1.26E-2 & 4.55E-3 \\
    0.1 & 2.67E-3 & 2.69E-2 & 9.44E-3 & 2.70E-3 & 2.70E-2 & 9.54E-3 & 2.71E-3 & 2.73E-2 & 9.20E-3 \\
    0.5 & 1.81E-2 & 2.31E-1 & 5.06E-2 & 1.84E-2 & 2.37E-1 & 5.19E-2 & 1.77E-2 & 2.28E-1 & 4.98E-2 \\
    1.0 & 6.92E-2 & 9.82E-1 & 1.30E-1 & 7.71E-2 & 1.05E+0 & 1.37E-1 & 6.75E-2 & 9.78E-1 & 1.30E-1 \\
    \hline
  \end{tabular}
  \caption{\small Inverse Burgers' problem: relative errors
    of $\alpha_1$, $\alpha_2$ and $u$
    versus~$\epsilon$ (noise level)
    obtained with the NLLSQ, VarPro-F1 and VarPro-F2 algorithms.
    Single sub-domain, NN $[2, 400, 1]$,
    $Q=30\times 30$, $Q_s$=100, $\lambda_{mea}$=1;
    $R_m=1.9$ with NLLSQ and VarPro-F1, and $R_m=2.0$ with VarPro-F2.
  }
  \label{tab_20}
\end{table}

Table~\ref{tab_20} illustrates the effect of the noisy measurement data ($\epsilon$)
 on the simulation accuracy of the NLLSQ, VarPro-F1 and
VarPro-F2 algorithms for the inverse Burgers problem.
It is observed that the accuracy of these algorithms
is quite robust to the noise.
For example, with $1\%$ noise in the measurement data the relative
errors of these methods are around $0.026\%$ for the computed $\alpha_1$
and around $0.2\%$ for the computed $\alpha_2$;
with $10\%$ noise in the measurement the relative errors are around $0.27\%$
for $\alpha_1$ and around $2.7\%$ for $\alpha_2$.

\subsection{Parametric Sine-Gordan Equation}
\label{sg_sec}

Consider the inverse parametric Sine-Gordan equation on
the domain $(x,t)\in\Omega=[0,1]\times[0,1]$,
\begin{subequations}\label{eq_60}
  \begin{align}\small
    &
    \frac{\partial^2u}{\partial t^2} - \alpha_1\frac{\partial^2u}{\partial x^2}
    + \alpha_2 u + \alpha_3\sin(u) = f(x,t), \label{eq_60a} \\
    &
    u(0,t) = g_1(t), \quad u(1,t) = g_2(t), \quad
    u(x,0) = h_1(x), \quad \left.\frac{\partial u}{\partial t}\right|_{(x,0)} = h_2(x),
    \label{eq_60c} \\
    &
    u(\xi_i,\eta_i) = S(\xi_i,\eta_i), \quad
    (\xi_i,\eta_i)\in\mbb Y\subset\Omega, \ 1\leqslant i\leqslant NQ_s,
  \end{align}
\end{subequations}
where $f$ is a prescribed source term,
$g_i$ ($i=1,2$) and $h_i$ ($i=1,2$) are prescribed boundary and initial conditions,
$\mbb Y$ is the set of random measurement points, and the constants $\alpha_i$ ($i=1,2,3$) and
the field $u(x,t)$ are the unknowns to be determined.
We employ the following manufactured analytic solution in
the tests,
\begin{equation}\label{eq_61}
  \small
  \left\{
  \begin{split}
    &
    \alpha_1^{ex} = \alpha_2^{ex} = \alpha_3^{ex} = 1, \\
    &
    u_{ex}(x,t) = \left[
      \frac52\cos\left(\pi x - \frac{2\pi}{5} \right)
      + \frac32\cos\left(2\pi x + \frac{3\pi}{10} \right)
      \right] \left[
      \frac52\cos\left(\pi t - \frac{2\pi}{5} \right)
      + \frac32\cos\left(2\pi t + \frac{3\pi}{10} \right)
      \right].
  \end{split}
  \right.
\end{equation}
Accordingly, $f$, $g_i$ ($i=1,2$), and $h_i$ ($i=1,2$) are chosen such that
the expressions in~\eqref{eq_61} satisfy~\eqref{eq_60a}--\eqref{eq_60c}.
The measurement data are given by equation~\eqref{eq_53},
in which $u_{ex}$ is given in~\eqref{eq_61}.
The $u$ errors are computed on a uniform $101\times 101$ grid in each sub-domain.
The notations here follow those of previous sub-sections.

\begin{figure}
  \small
  \centerline{
    \includegraphics[width=1.8in]{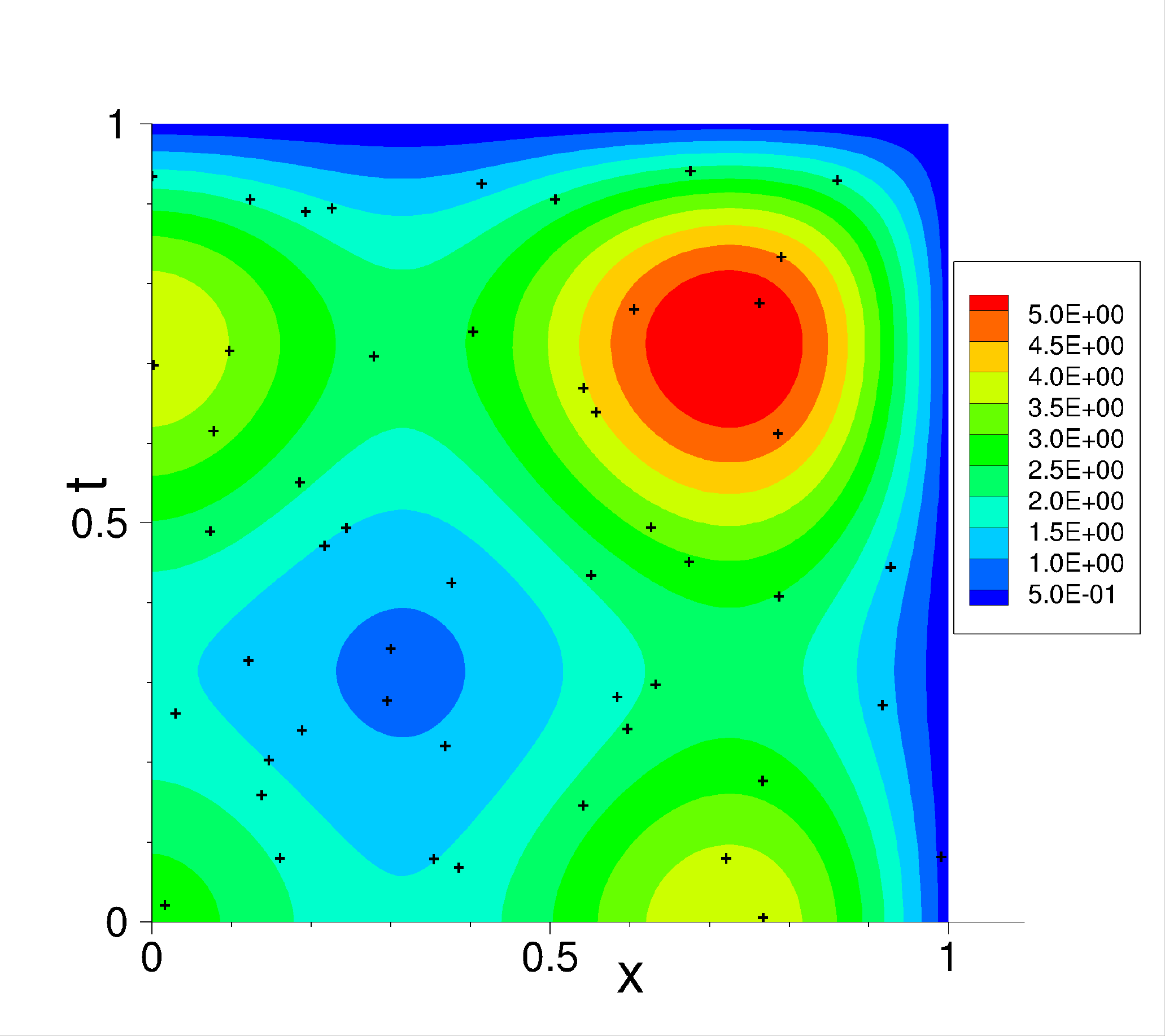}(a)
    \includegraphics[width=1.8in]{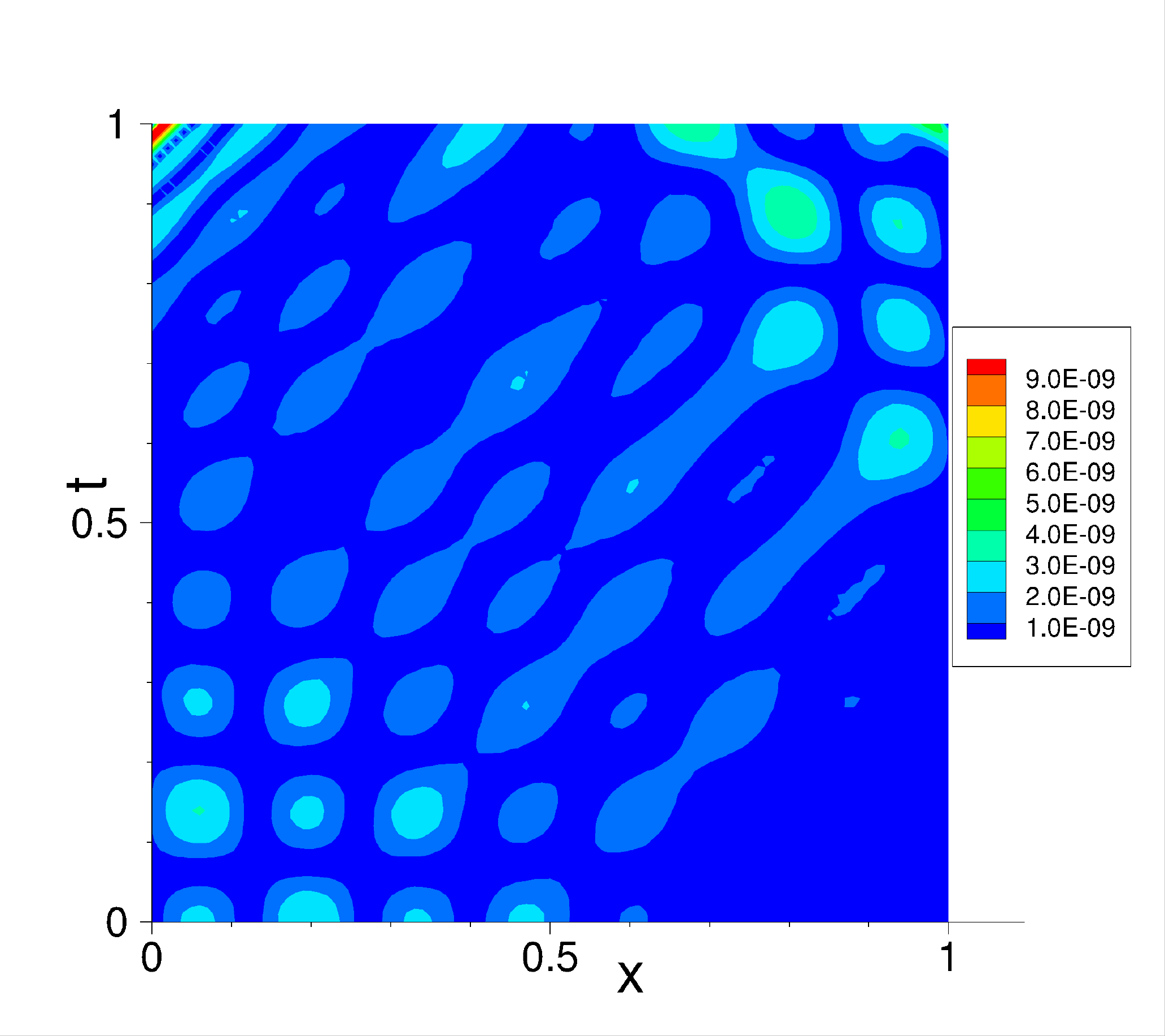}(b)
  }
  \caption{\small Inverse Sine-Gordan problem: distributions of (a) the VarPro-F2 solution for
    $u(x,t)$ and (b) its
    point-wise absolute error, with the measurement
    points shown as ``+'' symbols in (a).
    Single sub-domain, NN $[2,300,1]$,  $Q=25\times 25$,
    $Q_s=50$, $R_m=1.3$, $\lambda_{mea}$=1, $\epsilon=0$ (no noise in measurement data).
  }
  \label{fg_8}
\end{figure}

Figure~\ref{fg_8} shows distributions of the $u(x,t)$ solution and
its point-wise absolute error in $\Omega$
obtained by the VarPro-F2 algorithm, with $50$ random measurement points (no noise).
The other parameter values are provided in
the figure caption.
We can observe a high accuracy in the solution, with the maximum error
on the order of $10^{-8}$ in the domain.
In this simulation the relative errors for the
computed $\alpha_1$, $\alpha_2$ and $\alpha_3$ are
$2.07\times 10^{-10}$, $7.54\times 10^{-9}$ and $2.39\times 10^{-8}$,
respectively.


\begin{table}[tb]
  \centering\small
  \begin{tabular}{l c c c}
    \hline
    $M$ & $\alpha_1$ & $\alpha_2$ & $\alpha_3$ \\
    50 & -3.525085809204E+1 & -6.376776028198E+0 & 6.877670126056E+1 \\
    100 & 1.006414746681E+0 & 9.536423027578E-1 & 1.058573631607E+0  \\
    200 & 1.000000694463E+0 & 9.999830213887E-1 & 1.000056517879E+0 \\
    300 & 9.999999995066E-1 & 1.000000013620E+0 & 9.999999510014E-1 \\
    400 & 1.000000000001E+0 & 9.999999999962E-1 & 1.000000000096E+0 \\
    \hline
  \end{tabular}
  \caption{\small Inverse Sine-Gordan problem: $\alpha_i$ ($i=1,2,3$) versus $M$ (number of
    training parameters) obtained by the NLLSQ algorithm. Single sub-domain,
    NN $[2,M,1]$,
    $Q=25\times 25$, $Q_s=100$,   $R_m=1.5$, $\lambda_{mea}$=1, $\epsilon$=0.
  }
  \label{tab_21}
\end{table}

\begin{table}[tb]
  \centering\small
  \begin{tabular}{l|ccc|ccc|ccc}
    \hline
     & &  NLLSQ & & & VarPro-F1 & & & VarPro-F2 & \\ \cline{2-10}
    $M$ & $e_{\alpha_1}$ & $e_{\alpha_2}$ & $e_{\alpha_3}$ &  $e_{\alpha_1}$ & $e_{\alpha_2}$ & $e_{\alpha_3}$ &
       $e_{\alpha_1}$ & $e_{\alpha_2}$ & $e_{\alpha_3}$ \\ \hline
    50 & 3.63E+1 & 7.38E+0 & 6.78E+1 & 3.30E+1 & 6.51E+1 & 7.86E+0 & 4.55E+1 & 1.32E+2 & 4.14E+2  \\
    100 & 6.41E-3 & 4.64E-2 & 5.86E-2 & 1.76E-4 & 8.88E-3 & 1.14E-2 & 1.76E-4 & 8.88E-3 & 1.14E-2  \\
    200 & 6.94E-7 & 1.70E-5 & 5.65E-5 & 1.37E-7 & 4.34E-6 & 1.77E-5 & 1.37E-7 & 4.35E-6 & 1.77E-5  \\
    300 & 4.93E-10 & 1.36E-8 & 4.90E-8 & 2.92E-13 & 1.90E-10 & 9.64E-10 & 1.31E-10 & 1.99E-9 & 7.19E-9  \\
    400 & 1.43E-12 & 3.81E-12 & 9.58E-11 & 2.44E-12 & 3.67E-11 & 7.65E-11 & 3.31E-13 & 4.30E-12 & 7.85E-11  \\
    500 & 3.30E-13 & 7.76E-12 & 5.40E-11 & 2.84E-13 & 1.04E-11 & 4.73E-11 & 9.83E-13 & 5.90E-12 & 3.99E-11  \\
    \hline
  \end{tabular}
  \caption{\small Inverse Sine-Gordan problem: $\alpha_1$, $\alpha_2$ and $\alpha_3$ relative errors versus $M$ 
    obtained by the NLLSQ, VarPro-F1 and VarPro-F2 algorithms.
    Single sub-domain, NN $[2, M, 1]$,
    $Q_s=100$, $Q=25\times 25$, $\lambda_{mea}$=1, $\epsilon$=0;
    $R_m=1.5$ with NLLSQ, $R_m=1.3$ with VarPro-F1 and VarPro-F2.
  }
  \label{tab_22}
\end{table}

\begin{table}[tb]
  \centering\small
  \begin{tabular}{l|ll|ll|ll}
    \hline
     &   NLLSQ & &  VarPro-F1 & &  VarPro-F2 & \\ \cline{2-7}
    $M$ & $l^{\infty}$-u & $l^2$-u &  $l^{\infty}$-u & $l^2$-u & $l^{\infty}$-u & $l^2$-u  \\ \hline
    50 & 1.45E+0 & 4.90E-1 & 1.30E+0 & 4.10E-1 & 1.48E+0 & 4.85E-1 \\
    100 & 2.66E-2 & 5.24E-3 & 1.07E-2 & 2.23E-3 & 1.07E-2 & 2.23E-3 \\
    200 & 5.36E-6 & 6.69E-7 & 1.30E-6 & 2.90E-7 & 1.30E-6 & 2.90E-7 \\
    300 & 6.01E-9 & 4.75E-10 & 7.86E-10 & 9.21E-11 & 5.21E-9 & 3.44E-10 \\
    400 & 6.29E-11 & 3.77E-12 & 6.43E-11 & 4.25E-12 & 5.94E-10 & 8.53E-11 \\
    500 & 2.67E-11 & 1.13E-12 & 2.67E-11 & 1.10E-12 & 3.05E-10 & 2.08E-11 \\
    \hline
  \end{tabular}
  \caption{\small Inverse Sine-Gordan problem:  $u$ relative errors versus $M$
    obtained by the NLLSQ, VarPro-F1 and VarPro-F2 algorithms.
    Simulation settings and parameters follow those of Table~\ref{tab_22}.
  }
  \label{tab_23}
\end{table}

The convergence of the simulation results obtained by the NLLSQ, VarPro-F1
and VarPro-F2 algorithms is demonstrated by the data in Tables~\ref{tab_21} to \ref{tab_23}.
In these tests the number of training parameters ($M$) is varied systematically (no noise
in measurement),
while the other simulation parameters are fixed and their values
are provided in the table captions.
Table~\ref{tab_21} lists the computed $\alpha_i$ ($i=1,2,3$) values by the NLLSQ algorithm
corresponding to a set of $M$.
Table~\ref{tab_22} lists the relative errors of $\alpha_1$, $\alpha_2$ and $\alpha_3$
computed by NLLSQ, VarPro-F1 and VarPro-F2 corresponding to different $M$.
Table~\ref{tab_23} shows the $l^{\infty}$ and $l^2$ norms of the relative error for $u(x,t)$
in this set of simulations.
It is evident that the errors decrease exponentially with increasing number of
training parameters with these algorithms.

\begin{table}[tb]
  \centering\small
  \begin{tabular}{l|ccc|ccc|ccc}
    \hline
     & &  NLLSQ & & & VarPro-F1 & & & VarPro-F2 & \\ \cline{2-10}
    $\epsilon$ & $e_{\alpha_1}$ & $e_{\alpha_2}$ & $e_{\alpha_3}$ &  $e_{\alpha_1}$ & $e_{\alpha_2}$ & $e_{\alpha_3}$ &
       $e_{\alpha_1}$ & $e_{\alpha_2}$ & $e_{\alpha_3}$ \\ \hline
    0.0 & 1.93E-12 & 4.90E-11 & 1.35E-10 & 1.50E-12 & 4.42E-12 & 2.23E-11 & 3.61E-11 & 9.72E-10 & 3.02E-9  \\
    0.001 & 6.90E-4 & 2.66E-3 & 3.00E-3 & 6.88E-4 & 2.64E-3 & 3.09E-3 & 6.86E-4 & 2.63E-3 & 3.08E-3  \\
    0.005 & 3.44E-3 & 1.31E-2 & 1.45E-2 & 3.44E-3 & 1.32E-2 & 1.54E-2 & 3.43E-3 & 1.32E-2 & 1.54E-2  \\
    0.01 & 6.88E-3 & 2.63E-2 & 2.93E-2 & 6.86E-3 & 2.63E-2 & 3.08E-2 & 6.84E-3 & 2.62E-2 & 3.04E-2  \\
    0.05 & 3.38E-2 & 1.27E-1 & 1.39E-1 & 3.38E-2 & 1.30E-1 & 1.53E-1 & 3.39E-2 & 1.32E-1 & 1.60E-1  \\
    0.1 & 6.65E-2 & 2.49E-1 & 2.76E-1 & 6.65E-2 & 2.55E-1 & 3.05E-1 & 6.65E-2 & 2.57E-1 & 3.12E-1  \\
    0.5 & 2.67E-1 & 8.07E-1 & 5.90E-1 & 2.65E-1 & 7.69E-1 & 4.39E-1 & 2.66E-1 & 7.95E-1 & 5.41E-1  \\
    1.0 & 4.09E-1 & 1.01E+0 & 2.18E-1 & 4.12E-1 & 1.07E+0 & 5.43E-1 & 4.15E-1 & 1.10E+0 & 6.27E-1  \\
    \hline
  \end{tabular}
  \caption{\small Inverse Sine-Gordan problem: $\alpha_1$, $\alpha_2$ and $\alpha_3$ relative errors
    versus the noise level ($\epsilon$)
    obtained by the NLLSQ, VarPro-F1 and VarPro-F2 algorithms.
    Single sub-domain, NN: $[2, 400, 1]$,
    $Q=30\times 30$, $Q_s=50$, $\lambda_{mea}$=1;
    $R_m=1.5$ with NLLSQ, $R_m=1.3$ with VarPro-F1 and VarPro-F2;
  }
  \label{tab_24}
\end{table}

\begin{table}[tb]
  \centering\small
  \begin{tabular}{l|ll|ll|ll}
    \hline
     &   NLLSQ & &  VarPro-F1 & &  VarPro-F2 & \\ \cline{2-7}
    $\epsilon$ & $l^{\infty}$-u & $l^2$-u &  $l^{\infty}$-u & $l^2$-u & $l^{\infty}$-u & $l^2$-u  \\ \hline
    0.0 & 3.44E-11 & 4.16E-12 & 7.03E-11 & 5.01E-12 & 7.73E-10 & 1.57E-10 \\
    0.001 & 8.76E-4 & 3.71E-4 & 8.49E-4 & 3.70E-4 & 8.51E-4 & 3.69E-4 \\
    0.005 & 4.39E-3 & 1.86E-3 & 4.26E-3 & 1.85E-3 & 4.25E-3 & 1.85E-3 \\
    0.01 & 8.77E-3 & 3.71E-3 & 8.50E-3 & 3.70E-3 & 8.51E-3 & 3.70E-3 \\
    0.05 & 4.35E-2 & 1.87E-2 & 4.25E-2 & 1.87E-2 & 4.24E-2 & 1.86E-2 \\
    0.1 & 8.65E-2 & 3.78E-2 & 8.46E-2 & 3.77E-2 & 8.43E-2 & 3.76E-2 \\
    0.5 & 4.10E-1 & 1.95E-1 & 4.08E-1 & 1.95E-1 & 4.07E-1 & 1.95E-1 \\
    1.0 & 8.90E-1 & 3.83E-1 & 8.97E-1 & 3.84E-1 & 9.07E-1 & 3.86E-1 \\
    \hline
  \end{tabular}
  \caption{\small Inverse Sine-Gordan problem:  $u$ relative errors versus $\epsilon$
    obtained by the NLLSQ, VarPro-F1 and VarPro-F2 algorithms.
    Simulation settings and parameters follow those of Table~\ref{tab_24}.
  }
  \label{tab_024}
\end{table}

The effect of noise in the measurement data on the simulation accuracy
is illustrated by Tables~\ref{tab_24} and~\ref{tab_024} for
the NLLSQ, VarPro-F1 and VarPro-F2 algorithms.
The relative errors of $\alpha_1$, $\alpha_2$, $\alpha_3$, and
$u(x,t)$ corresponding to a range of noise levels are
provided in these two tables.
The other crucial simulation parameters are provided in the caption
of Table~\ref{tab_24}.
The accuracy in the computation results deteriorates as the measurement data becomes more noisy.
With $1\%$ measurement noise ($\epsilon=0.01$) the relative errors
of the computed $\alpha_i$ ($i=1,2,3$) are around $0.7\sim 3\%$, and the relative error of
$u$ ($l^2$ norm) is around $0.4\%$ with the three algorithms.
With $5\%$ measurement noise ($\epsilon=0.05$) the relative errors
of the computed $\alpha_i$ are around $3\sim 15\%$ and the relative error of $u$ ($l^2$ norm)
is less than $2\%$.

\subsection{Helmholtz Equation with Inverse Variable Coefficient }
\label{varhelm}

In the last example, we use our method to study
a problem involving an inverse coefficient field.
Consider the domain $(x,y)\in\Omega=[0,1.5]\times[0,1.5]$
and the inverse problem on $\Omega$,
\begin{subequations}\label{eq_62}
  \begin{align}\small
    &
    \frac{\partial^2u}{\partial x^2} + \frac{\partial^2u}{\partial y^2} - \gamma(x,y)u
    = f(x,y), \label{eq_62a}\\
    &
    u(a_1,y) = g_1(y), \quad
    u(b_1,y) = g_2(y), \quad
    u(x,a_2) = g_3(x), \quad
    u(x,b_2) = g_4(x), \label{eq_62b} \\
    &
    u(\xi_i,\eta_i) = S(\xi_i,\eta_i), \quad
    (\xi_i,\eta_i)\in\mbb Y\subset\Omega,\ 1\leqslant i\leqslant NQ_s,
  \end{align}
\end{subequations}
where $f$ is a prescribed source term, $g_i$ ($1\leqslant i\leqslant 4$) denote
the prescribed boundary data, $\mbb Y$ is the set of measurement points, $S(\xi_i,\eta_i)$
denotes the measurement data at the random point $(\xi_i,\eta_i)$,
and $\gamma(x,y)$ and $u(x,y)$ are two field functions to be determined.
We employ the following manufactured solutions,
\begin{equation}\label{eq_63}
  \small
  \left\{
  \begin{split}
    &
    \gamma_{ex}(x,y) = 100\left[1+\frac14\sin(2\pi x) + \frac14\sin(2\pi y) \right], \\
    &
    u_{ex}(x,y) = \left[
    \frac52\sin\left(\pi x - \frac{2\pi}{5} \right) + \frac32\cos\left(2\pi x + \frac{3\pi}{10} \right)
    \right] \left[
    \frac52\sin\left(\pi y - \frac{2\pi}{5} \right) + \frac32\cos\left(2\pi y + \frac{3\pi}{10} \right)
    \right].
  \end{split}
  \right.
\end{equation}
$f$ and $g_i$ ($1\leqslant i\leqslant 4$) are chosen accordingly such that
the expressions in~\eqref{eq_63} satisfy the equations~\eqref{eq_62a}--\eqref{eq_62b}.
The measurement data $S(\xi_i,\eta_i)$ are given by equation~\eqref{eq_53},
in which the $u_{ex}$ is given in~\eqref{eq_63}.
The relative errors for $u(x,y)$ are defined in~\eqref{eq_a55}, and the relative errors for
$\gamma(x,y)$ are defined analogously.
The  $\gamma(x,y)$ and $u(x,y)$ errors reported below
are computed on a uniform $Q_{eval}=101\times 101$ grid in each sub-domain.
The notations here follow those of previous subsections.

We employ the algorithm modification as outlined in Remark~\ref{rem_05} for solving this problem.
Compared with previous subsections,
the main change here lies in that the local neural network on each sub-domain
contains two nodes in the output layer, one representing $u(x,y)$ and
the other $\gamma(x,y)$.
The output-layer coefficients contributing to $\gamma(x,y)$ play the role
of the inverse parameters.
We also find it preferable to regularize the the output-layer coefficients
that contribute to $\gamma(x,y)$ (or $u(x,y)$) for this problem.
For regularization we employ the extra terms for the underlying loss
function, $\frac{\lambda_1^2}{2}\|\bm\alpha \|^2 + \frac{\lambda_2^2}{2}\|\bm\beta \|^2$,
where $\bm\alpha$ and $\bm\beta$ denote the vectors of output-layer coefficients
for $\gamma(x,y)$ and $u(x,y)$, respectively, and the prescribed non-negative
constants $\lambda_1$ and $\lambda_2$ are the corresponding
regularization coefficients.


\begin{figure}
  \small
  \centerline{
    \includegraphics[width=1.5in]{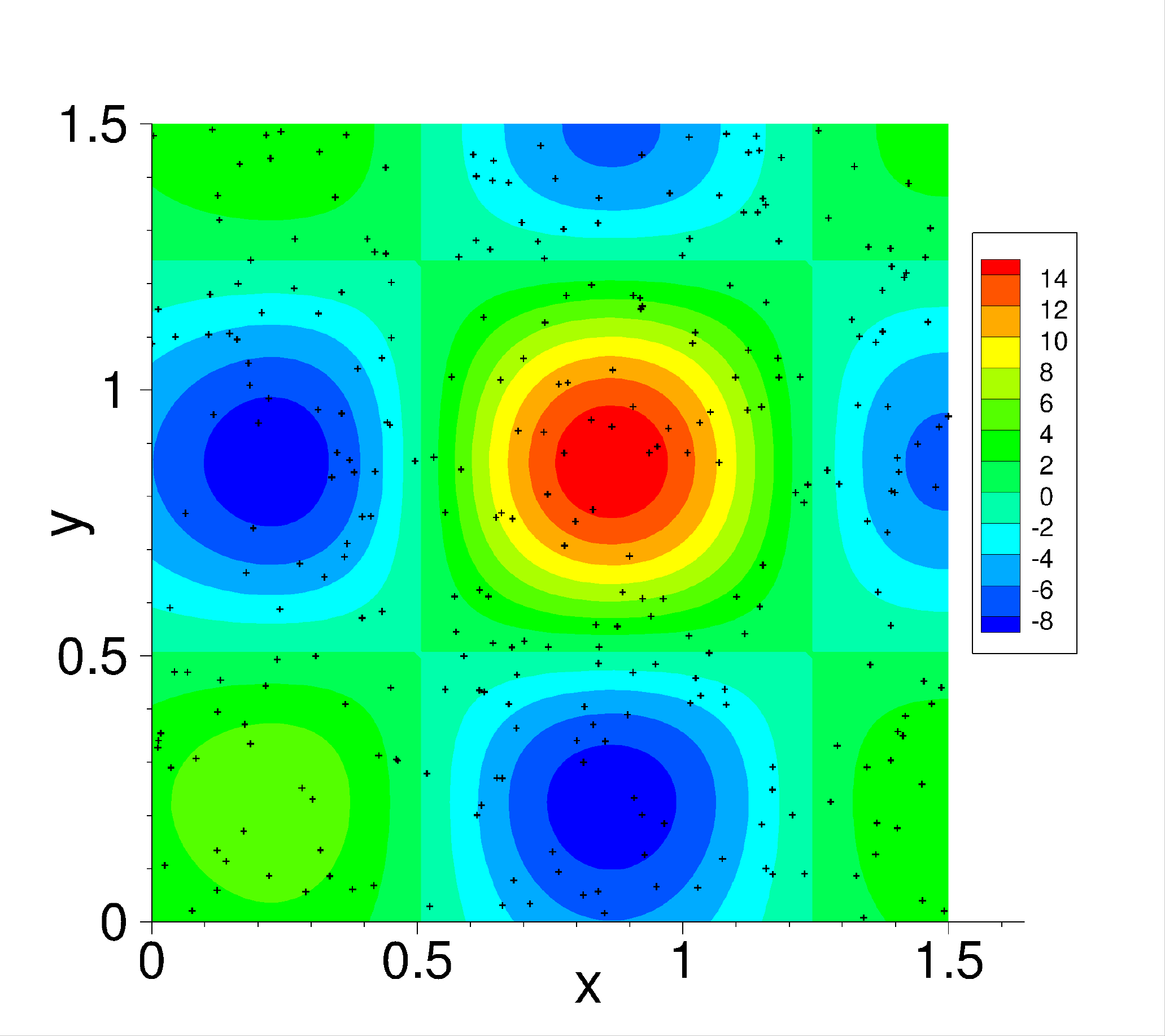}(a)
    \includegraphics[width=1.5in]{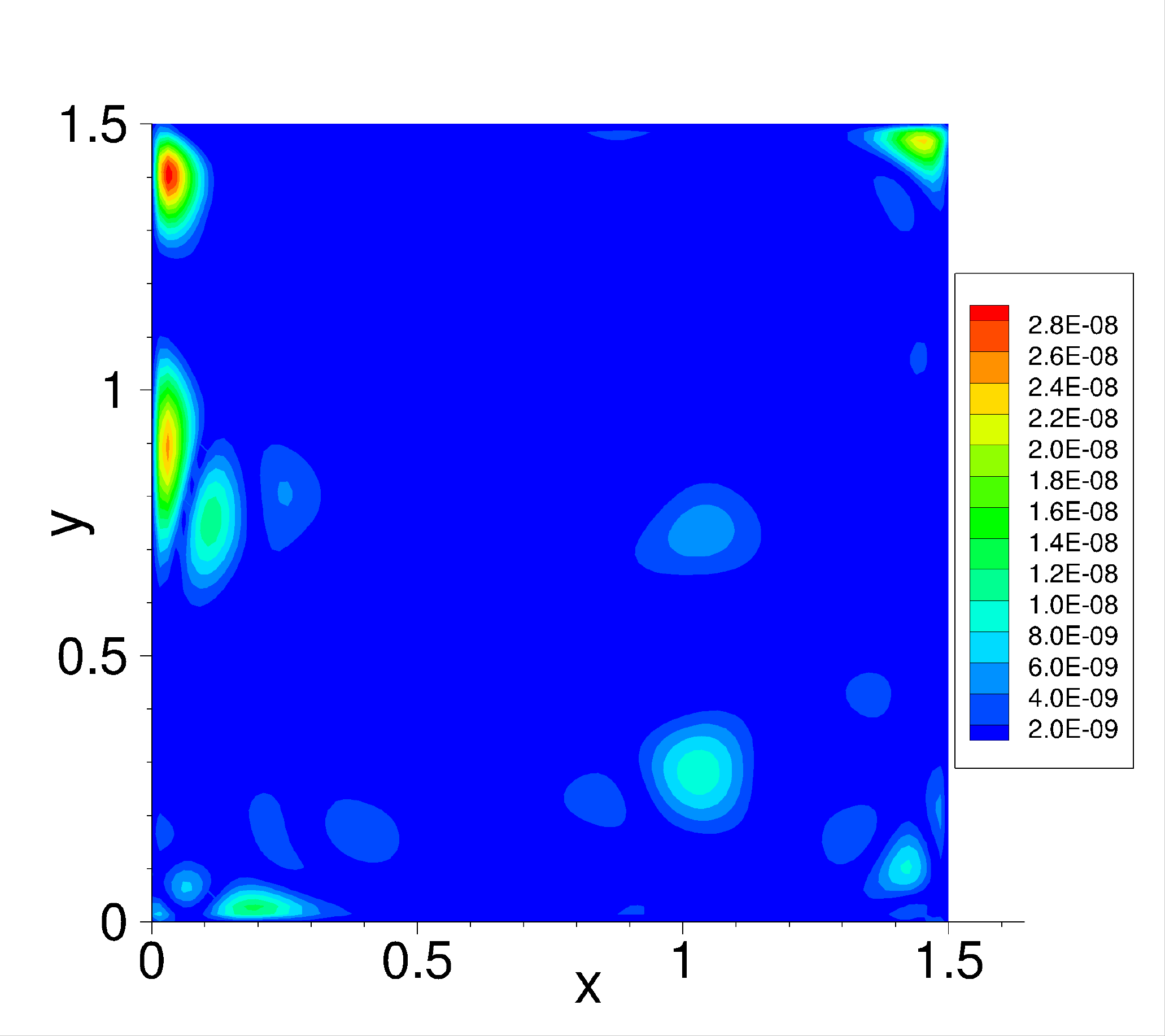}(b)
    \includegraphics[width=1.5in]{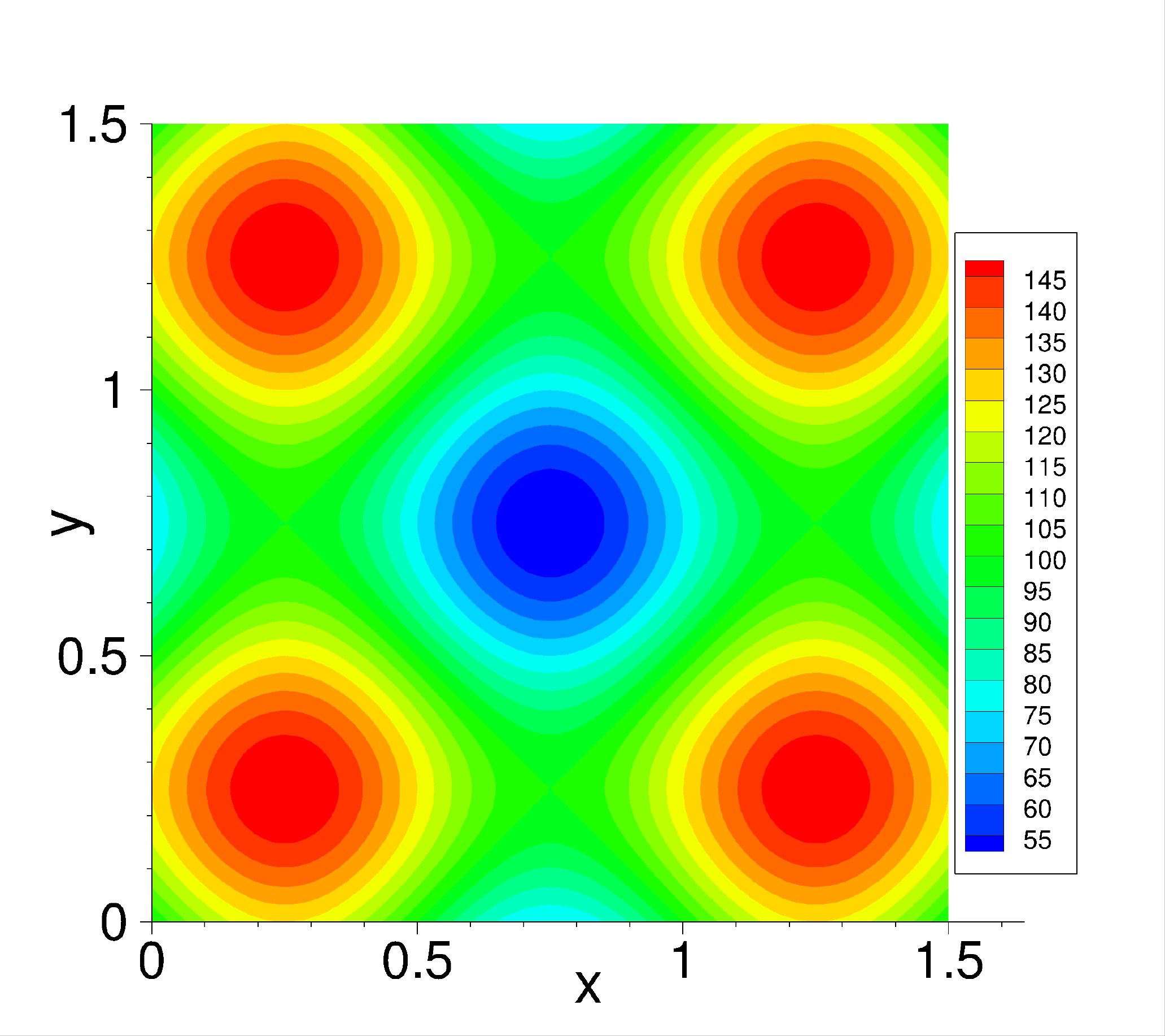}(c)
    \includegraphics[width=1.5in]{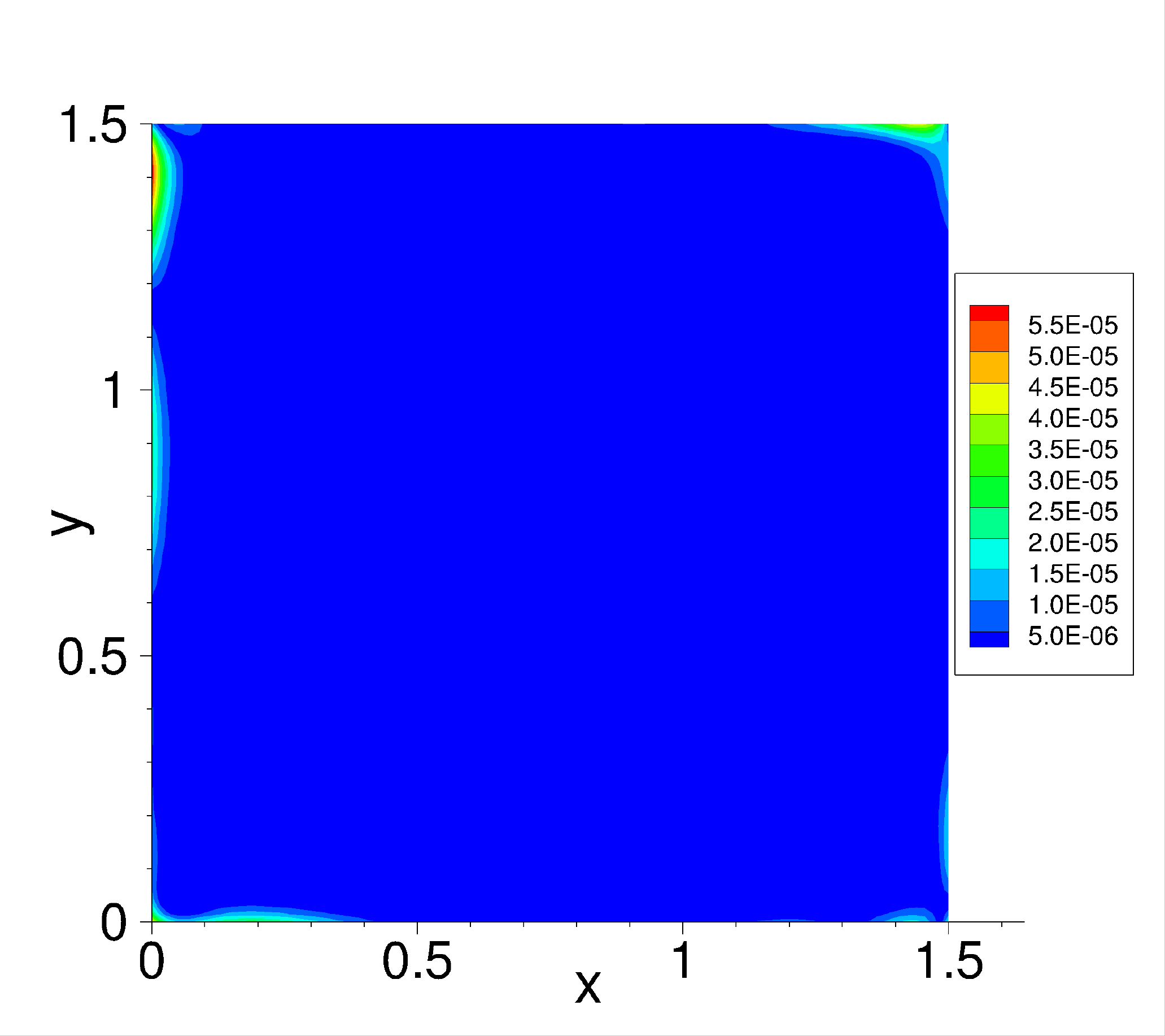}(d)
  }
  \caption{\small Inverse variable-coefficient Helmholtz problem:
    distributions of (a) the NLLSQ solution for $u(x,y)$ and (b) its point-wise absolute error,
    (c) the NLLSQ solution for $\gamma(x,y)$ and (d) its point-wise absolute error,
    with the $Q_s=300$ random measurement points shown in (a) as ``+'' symbols. 
    Single sub-domain, NN $[2,400,1]$, $Q=30\times 30$,
    $R_m=1.5$, $\lambda_{mea}=1$, $\epsilon=0$ (no noise), $\lambda_1=\lambda_2=0$ (no regularization).
  }
  \label{fg_9}
\end{figure}

Figure~\ref{fg_9} shows distributions of the solutions for $u(x,y)$
and $\gamma(x,y)$, and their point-wise absolute errors, obtained by
the NLLSQ algorithm. The random measurement points are also shown in
Figure~\ref{fg_9}(a). The figure caption lists the main parameter
values for this simulation. We observe a fairly high accuracy,
with the maximum $u$ error on the order of $10^{-8}$
and the maximum $\gamma$ error on the order of $10^{-5}$ in the domain.


\begin{table}[tb]
  \centering\small
  \begin{tabular}{l|ll| c c c c}
    \hline\hline
     & & & $l^{\infty}$-$\gamma$ & $l^2$-$\gamma$ & $l^{\infty}$-u & $l^2$-u \\ \hline
    collocation & $Q=$ & 5$\times$5 & 9.84E-2 & 3.15E-2 & 7.94E-3 & 4.38E-4 \\
    point test & & 10$\times$10 & 4.05E-3 & 2.90E-4 & 1.63E-4 & 1.02E-5 \\
    & & 15$\times$15 & 1.15E-4 & 4.94E-6 & 3.38E-6 & 2.99E-7 \\
    & & 20$\times$20 & 4.99E-6 & 4.84E-7 & 2.65E-7 & 2.70E-8 \\
    & & 25$\times$25 & 4.42E-6 & 3.57E-7 & 3.34E-7 & 2.90E-8 \\
    & & 30$\times$30 & 1.56E-6 & 1.46E-7 & 1.19E-7 & 1.08E-8 \\
    & & 35$\times$35 & 1.61E-6 & 1.73E-7 & 1.26E-7 & 1.15E-8 \\
    \hline\hline
    training & $M=$ & 50 & 6.79E+0 & 1.67E+0 & 1.63E+0 & 5.67E-1 \\
    parameter & & 100 & 8.64E-2 & 1.10E-2 & 8.44E-3 & 2.00E-3 \\
    test & & 200 & 1.98E-4 & 2.08E-5 & 6.46E-6 & 6.42E-7 \\
    & & 300 & 2.07E-6 & 1.26E-7 & 3.13E-8 & 2.46E-9 \\
    & & 400 & 5.08E-7 & 2.95E-8 & 5.76E-9 & 5.29E-10 \\
    & & 500 & 1.61E-7 & 1.18E-8 & 1.95E-9 & 1.91E-10 \\
    \hline\hline
    measurement & $Q_s=$ & 10 & 1.36E-3 & 3.11E-4 & 2.49E-4 & 6.60E-5 \\
    point test & & 30 & 1.93E-4 & 3.69E-5 & 2.03E-5 & 4.24E-6 \\
    & & 50 & 2.69E-5 & 4.08E-6 & 2.03E-6 & 3.83E-7 \\
    & & 100 & 3.18E-6 & 2.95E-7 & 2.25E-7 & 1.78E-8 \\
    & & 200 & 5.30E-6 & 2.53E-7 & 8.21E-8 & 5.48E-9 \\
    & & 300 & 1.45E-6 & 9.93E-8 & 2.45E-8 & 2.00E-9 \\
    & & 400 & 2.65E-6 & 8.80E-8 & 1.40E-8 & 1.79E-9 \\
    \hline\hline
    noise level & $\epsilon=$ & 0.0 & 1.61E-6 & 1.73E-7 & 1.26E-7 & 1.15E-8 \\
    test & & 0.0005 & 5.15E-2 & 5.54E-3 & 3.32E-3 & 4.54E-4 \\
    & & 0.001 & 7.57E-2 & 7.31E-3 & 5.46E-3 & 6.63E-4 \\
    & & 0.005 & 2.69E-1 & 2.92E-2 & 2.28E-2 & 3.12E-3 \\
    & & 0.01 & 4.26E-1 & 4.96E-2 & 3.94E-2 & 5.98E-3 \\
    & & 0.05 & 1.50E+0 & 2.00E-1 & 1.44E-1 & 2.71E-2 \\
    & & 0.1 & 1.85E+0 & 2.97E-1 & 2.26E-1 & 4.97E-2 \\
    \hline\hline
  \end{tabular}
  \caption{\small Inverse variable-coefficient Helmholtz problem:
    relative $l^{\infty}$ and $l^2$ errors of $\gamma(x,y)$ and $u(x,y)$ in several tests with
    the NLLSQ algorithm.
    Single sub-domain, NN $[2,M,1]$.
    In collocation point test, $M=400$, $Q_s=100$, $\epsilon=0$, $Q$ is varied.
    In training parameter test, $Q=30\times 30$, $Q_s=300$, $\epsilon=0$, $M$ is varied.
    In measurement point test, $Q=25\times 25$, $M=300$, $\epsilon=0$, $Q_s$ is varied.
    In noise level test, $Q=35\times 35$, $M=400$, $Q_s=100$, $\epsilon$ is varied.
    $R_m=1.5$ and $\lambda_{mea}=1$ in all tests.
    No regularization ($\lambda_1=\lambda_2=0$).
  }
  \label{tab_25}
\end{table}

Table~\ref{tab_25}
lists the relative errors for $\gamma$ and $u$ ($l^{\infty}$ and $l^2$ norms)
computed by the NLLSQ algorithm in several sets of tests, with respect to
$Q$ (number of collocation points), $M$ (number of training parameters),
$Q_s$ (number of random measurement points), and $\epsilon$ (noise level).
The settings and the simulation parameter values are provided in
the table caption for each set of tests.
One can observe an approximately exponential decrease in the $\gamma$ and $u$ errors
with respect to $Q$, $M$ and $Q_s$ (before saturation).
One can also observe the deterioration in the simulation accuracy
with increasing noise level in the measurement data.
Note that no regularization is employed in these simulations.
The noise appears to affect the $\gamma$ results more significantly than $u$.
For example, with $1\%$ noise ($\epsilon=0.01$) in the measurement data,
the maximum ($l^{\infty}$) relative error for $\gamma(x,y)$ is around $43\%$ and
the $l^{2}$ relative error is around $5\%$, while for $u(x,y)$ these
errors are around $4\%$ and $0.6\%$ respectively.


\begin{table}[tb]
  \centering\small
  \begin{tabular}{l|rl|cccc}
    \hline\hline
    & noise & level  & $l^{\infty}$-$\gamma$ & $l^2$-$\gamma$ & $l^{\infty}$-u & $l^2$-u \\ \hline
    NLLSQ & $\epsilon=$ & 0.0 & 1.39E-4 & 7.11E-6 & 2.57E-6 & 2.24E-7 \\
    & & 0.0005 & 1.47E-2 & 1.42E-3 & 7.48E-4 & 2.17E-4 \\
    & & 0.001 & 2.92E-2 & 2.85E-3 & 1.50E-3 & 4.34E-4 \\
    & & 0.005 & 1.27E-1 & 1.36E-2 & 7.33E-3 & 2.17E-3 \\
    & & 0.01 & 2.14E-1 & 2.58E-2 & 1.48E-2 & 4.32E-3 \\
    & & 0.05 & 7.38E-1 & 1.02E-1 & 7.35E-2 & 2.14E-2 \\
    & & 0.1 & 8.43E-1 & 1.84E-1 & 1.48E-1 & 4.26E-2 \\
    \hline\hline
    VarPro-F1 & $\epsilon=$ & 0.0 & 3.16E-4 & 3.02E-5 & 7.75E-6 & 1.11E-6 \\
    & & 0.0005 & 1.65E-2 & 1.93E-3 & 7.23E-4 & 2.18E-4 \\
    & & 0.001 & 3.21E-2 & 3.81E-3 & 1.44E-3 & 4.35E-4 \\
    & & 0.005 & 1.30E-1 & 1.78E-2 & 7.29E-3 & 2.16E-3 \\
    & & 0.01 & 2.31E-1 & 3.34E-2 & 1.48E-2 & 4.31E-3 \\
    & & 0.05 & 9.30E-1 & 1.16E-1 & 7.19E-2 & 2.14E-2 \\
    & & 0.1 & 1.02E+0 & 2.08E-1 & 1.42E-1 & 4.29E-2 \\
    \hline\hline
    VarPro-F2 & $\epsilon=$ & 0.0 & 5.71E-4 & 4.93E-5 & 1.28E-5 & 1.29E-6 \\
    & & 0.0005 & 2.18E-2 & 2.46E-3 & 7.19E-4 & 2.23E-4 \\
    & & 0.001 & 4.33E-2 & 4.90E-3 & 1.44E-3 & 4.46E-4 \\
    & & 0.005 & 1.88E-1 & 2.30E-2 & 6.97E-3 & 2.21E-3 \\
    & & 0.01 & 3.37E-1 & 4.42E-2 & 1.37E-2 & 4.39E-3 \\
    & & 0.05 & 1.60E+0 & 1.93E-1 & 6.91E-2 & 2.18E-2 \\
    & & 0.1 & 2.66E+0 & 3.11E-1 & 1.38E-1 & 4.31E-2 \\
    \hline\hline
  \end{tabular}
  \caption{\small Inverse variable-coefficient Helmholtz problem:
    relative $l^{\infty}$ and $l^2$ errors of $\gamma(x,y)$ and $u(x,y)$ versus 
    $\epsilon$ (noise level) by the NLLSQ, VarPro-F1 and VarPro-F2 algorithms.
    Single sub-domain, NN [2,400,1], $Q=30\times 30$, $Q_s=300$, $\lambda_{mea}=1$;
    $R_m=1.5$ with NLLSQ, $R_m=2.5$ with VarPro-F1, $R_m=3.0$ with VarPro-F2;
    Regularization coefficients: $(\lambda_1,\lambda_2)=(1E-8,1E-8)$ with NLLSQ,
    $(\lambda_1,\lambda_2)=(1E-7,0)$ with VarPro-F1, and
    $(\lambda_1,\lambda_2)=(1E-6,0)$ with VarPro-F2.
  }
  \label{tab_26}
\end{table}

Table~\ref{tab_26} illustrates the effect of noise in the measurement data
on the accuracy of the NLLSQ, VarPro-F1 and VarPro-F2 algorithms.
The relative errors for $\gamma$ and $u$ corresponding to different noise levels
have been shown. In these simulations the output-layer coefficients
for $\gamma(x,y)$ (and also for $u(x,y)$ with NLLSQ) have been regularized,
with the regularization coefficients and the other simulation parameter values
given in the table caption.
The regularization generally improves the accuracy in the presence of noise.



\section{Concluding Remarks}
\label{sec:summary}



In this paper we have presented a method for solving inverse parametric PDE problems
based on randomized neural networks.
This method extends the local extreme learning machine (locELM) technique
to inverse PDEs.
The field solution is represented by a set of
local random-weight neural networks (randomly assigned but fixed hidden-layer coefficients,
trainable output-layer coefficients), one for each sub-domain.
The local neural networks are coupled through the $C^k$ (with $k$ related to
the PDE order) continuity conditions
on the shared sub-domain boundaries.
The inverse parameters of the PDE and the trainable parameters of the local neural networks
are the unknowns to be determined in the system.

Three algorithms are developed for training
the neural network  to solve the inverse problem.
The first algorithm (NLLSQ) computes the inverse parameters and the trainable
network parameters all together by the nonlinear least squares method and
is an extension of the nonlinear least squares
method with perturbations (NLLSQ-perturb) of~\cite{DongL2021}
(developed for forward PDEs) to inverse PDE problems.
The second and the third algorithms are based on the variable projection idea.
The second algorithm (VarPro-F1) employs variable projection
to eliminate the inverse parameters from the problem and attain
a reduced problem about the trainable network parameters only.
Then the reduced problem is solved first by the NLLSQ-perturb algorithm for
the trainable network parameters, and the inverse parameters are computed afterwards by
the linear least squares method.
The third algorithm (VarPro-F2) provides a reciprocal formulation with
variable projection. It eliminates the trainable network parameters (or equivalently
the field solution) from the problem first to arrive at a reduced problem
about the inverse parameters only. Then the inverse parameters are computed first
by solving the reduced problem with the NLLSQ-perturb algorithm,
and afterwards the trainable network parameters are computed based on the inverse
parameters already obtained. The VarPro-F2 algorithm is suitable
for parametric PDEs that are linear with respect to the field solution.
For PDEs that are nonlinear with respect to the field solution,
this algorithm needs to be combined with a Newton iteration.


The presented method is numerically tested using several inverse parametric
PDE problems. It is also compared with the PINN method (see Appendix C).
For smooth solutions and noise-free data, the errors for the inverse parameters and
the field solution computed by the NLLSQ, VarPro-F1 and
VarPro-F2 algorithms decrease exponentially with respect to
the number of collocation points and the number of training parameters
in the neural network. When these parameters become large the errors
can reach a level close to the machine accuracy.
These characteristics are in some sense analogous to those observed for
the forward PDE problems in~\cite{DongL2021,DongY2022rm}.
For noisy measurement data, these algorithms can
produce computation results with good accuracy, indicating robustness
of the method.
We observe that, in the presence of noise,
by scaling the measurement residual by a factor $\lambda_{mea}$
($0<\lambda_{mea}<1$) one can in general improve the simulation
accuracy of the current method markedly, while this scaling
only slightly degrades the accuracy for the noise-free data.
The comparison with PINN shows that the current method has an advantage
in terms of both accuracy and the network training time. In particular,
for the noise-free data the current method exhibits an accuracy significantly
higher than PINN.

These test results suggest that the method developed herein
is an effective and promising technique
for computing inverse PDEs. The exponential convergence exhibited by the method
is especially interesting, suggesting a high accuracy of this technique.
We anticipate that this technique will be a useful and meaningful addition
to the arsenal for tackling
this class of problems and be instrumental in computational science and engineering
applications.


\section*{Acknowledgment}
This work was partially supported by
NSF (DMS-2012415). 

\section*{Appendix A.~Nonlinear Least Squares Algorithm
  with Perturbations (NLLSQ-perturb)}

We summarize the nonlinear least squares algorithm
with perturbations (NLLSQ-perturb) below in
Algorithm~\ref{alg_1}, which is adapted from
the one developed in~\cite{DongL2021} with certain modifications.

\begin{algorithm}[thb]\small
  \DontPrintSemicolon
  \SetKwInOut{Input}{input}\SetKwInOut{Output}{output}

  \Input{max perturbation magnitude $\delta>0$;
    initial guess vector  $\bm\theta_0$; routine for computing
    residual; routine for computing Jacobian matrix;
    perturbation  flag $\eta$ (integer, $0$ or $1$);
    tolerance $\varepsilon>0$; max-nllsq-iterations (positive integer);
    max-sub-iterations (non-negative integer).
  }
  \Output{solution vector $\bm \theta$, associated cost $c$}
  \BlankLine
  invoke the scipy.optimize.least\_squares routine, with the inputs (initial guess
  $\bm \theta_0$, routines for residual/jacobian-matrix calculations,
  and max-nllsq-iterations)\;
  set $\bm\theta\leftarrow$ returned solution, and set $c\leftarrow$ returned cost\;
  \If{c is below $\varepsilon$}{return $\bm \theta$ and $c$\;}
  \BlankLine
  \For{$i\leftarrow 1$ \KwTo max-sub-iterations}{
    generate a uniform random number $\xi$ on the interval $[0,1]$\;
    set $\delta_1 \leftarrow \xi\delta$\;
    generate a uniform random vector $\Delta\bm\theta$ of the
    same shape as $\bm\theta$ on the interval [$-\delta_1$, $\delta_1$]\;
    \BlankLine
    \eIf{$\eta$ is $0$}{
      $\bm\vartheta_0 \leftarrow \Delta\bm\theta$\;
    }
    {$\bm\vartheta_0 \leftarrow \bm\theta + \Delta\bm\theta$\; }
    \BlankLine
    invoke the scipy.optimize.least\_squares routine, with the inputs (initial guess
    $\bm\vartheta_0$, routines for residual/jacobian-matrix calculations,
    and max-nllsq-iterations)\;
    \If{the returned cost is less than $c$}{
      set $\bm\theta\leftarrow$ the returned solution, and
      set $c\leftarrow$ the returned cost\;
      \If{c is below $\varepsilon$}{
        break\;
      }
    } 
  }
  return $\bm\theta$ and $c$\;
  \caption{NLLSQ-perturb (nonlinear least squares with perturbations) algorithm}
  \label{alg_1}
\end{algorithm}

In this algorithm, $\delta$ controls the maximum magnitude of each component
of the random perturbation vector $\Delta\bm\theta$.
The vector $\bm\theta_0$ provides the initial guess to the solution
of the nonlinear least squares problem. 
If the returned solution
from the scipy least\_squares() routine corresponding to $\bm\theta_0$ is not acceptable
(i.e.~the returned cost exceeding the tolerance $\varepsilon$), then
a sub-iteration procedure is triggered in which new initial guesses ($\bm\vartheta_0$)
are generated by perturbing either the origin or the best approximation
to the solution obtained so far with a random vector. The scipy
least\_squares() routine is invoked with the new initial guesses until
an acceptable solution is obtained or until the maximum number of sub-iterations
is reached. The integer flag $\eta$ controls around which point the perturbation is
performed. If $\eta=0$ the new initial guess is generated by perturbing
the origin. Otherwise, the current best approximation to the solution
is perturbed to generate a new initial guess.
The parameter ``max-nllsq-iterations'' controls the maximum number of
iterations (e.g.~the maximum number of residual function evaluations)
in the scipy least\_squares() routine.
The parameter ``max-sub-iterations'' controls the maximum number of 
sub-iterations for the initial guess perturbation.
One can turn off the perturbation in the NLLSQ-perturb algorithm
by setting max-sub-iterations to zero.
Note that the scipy least\_squares() function requires
two routines in the input, one for computing
the residual and the other for computing the Jacobian
matrix for an arbitrary given approximation to the solution.

\section*{Appendix B.~Matrices in the VarPro-F2 Algorithm}

The matrices in the expressions~\eqref{eq_43} are given by,
\begin{equation}\label{eq_ab1}\footnotesize
  \begin{split}
    &
    \mbs b^{\text{pde}}=\begin{bmatrix}\vdots\\ f(\mbs x_p^e)\\ \vdots  \end{bmatrix}_{NQ\times 1},\
    \mbs b^{\text{mea}}=\begin{bmatrix}\vdots\\ S(\bm\xi_p^e)\\ \vdots  \end{bmatrix}_{NQ_s\times 1},\
    \mbs b^{\text{bc1}}=\begin{bmatrix}\vdots\\ g(a_1,y_{p(1,j)}^{e(1,l)})\\ \vdots
    \end{bmatrix}_{N_2Q_2\times 1}, \\
    &
    \mbs b^{\text{bc2}}=\begin{bmatrix}\vdots\\ g(b_1,y_{p(Q_1,j)}^{e(N_1,l)})\\ \vdots
    \end{bmatrix}_{N_2Q_2\times 1}, \
    \mbs b^{\text{bc3}}=\begin{bmatrix}\vdots\\ g(x_{p(i,1)}^{e(m,1)},a_2)\\ \vdots
    \end{bmatrix}_{N_1Q_1\times 1}, \
    \mbs b^{\text{bc4}}=\begin{bmatrix}\vdots\\ g(x_{p(i,Q_2)}^{e(m,N_2)},a_2)\\ \vdots
    \end{bmatrix}_{N_1Q_1\times 1}, \\
    &
    \mbs H^{\text{pde}}=\begin{bmatrix}h_{ij}^{\text{pde}}  \end{bmatrix}_{NQ\times NM},\
    \mbs H^{\text{mea}}=\begin{bmatrix}h_{ij}^{\text{mea}}  \end{bmatrix}_{NQ_s\times NM},\
    \mbs H^{\text{bc1}}=\begin{bmatrix}h_{ij}^{\text{bc1}}  \end{bmatrix}_{N_2Q_2\times NM},\
    \mbs H^{\text{bc2}}=\begin{bmatrix}h_{ij}^{\text{bc2}}  \end{bmatrix}_{N_2Q_2\times NM},\\
    &
    \mbs H^{\text{bc3}}=\begin{bmatrix}h_{ij}^{\text{bc3}}  \end{bmatrix}_{N_1Q_1\times NM},\
    \mbs H^{\text{bc4}}=\begin{bmatrix}h_{ij}^{\text{bc4}}  \end{bmatrix}_{N_1Q_1\times NM},\
    \mbs H^{\text{ck1}}=\begin{bmatrix}h_{ij}^{\text{ck1}}  \end{bmatrix}_{(N-N_2)Q_2\times NM},\\
    &
    \mbs H^{\text{ck2}}=\begin{bmatrix}h_{ij}^{\text{ck2}}  \end{bmatrix}_{(N-N_2)Q_2\times NM},\
    \mbs H^{\text{ck3}}=\begin{bmatrix}h_{ij}^{\text{ck3}}  \end{bmatrix}_{(N-N_1)Q_1\times NM},\
    \mbs H^{\text{ck4}}=\begin{bmatrix}h_{ij}^{\text{ck4}}  \end{bmatrix}_{(N-N_1)Q_1\times NM}.
  \end{split}
\end{equation}
In the matrices $\mbs H^{\text{pde}}$ and $\mbs H^{\text{mea}}$ the only non-zero terms are,
\begin{equation}\label{eq_ab2}\footnotesize
  \left\{
  \begin{split}
  h_{ij}^{\text{pde}}=&\alpha_1\mathcal{L}_1\phi_{eq}(\mbs x_p^e) + \dots +
  \alpha_n\mathcal{L}_n\phi_{eq}(\mbs x_p^e)
  + \mathcal{F}\phi_{eq}(\mbs x_p^e), \\
  & i=(e-1)Q+p, \ j=(e-1)M+q,
  \quad \text{for}\ 1\leqslant (e,p,q)\leqslant (N,Q,M); \\
  h_{ij}^{\text{mea}}=& \mathcal{M}\phi_{eq}(\bm\xi_p^e), \ \
  i=(e-1)Q_s+p, \ j=(e-1)M+q,
  \quad \text{for}\ 1\leqslant (e,p,q)\leqslant (N,Q_s,M).
  \end{split}
  \right.
\end{equation}
In the matrices $\mbs H^{\text{bc1}}$, $\mbs H^{\text{bc2}}$, $\mbs H^{\text{bc3}}$
and $\mbs H^{\text{bc4}}$ the only non-zero terms are,
\begin{equation}\label{eq_ab3}\footnotesize
  \left\{
  \begin{split}
    h_{ij}^{\text{bc1}}=& \mathcal{B}\phi_{eq}(a_1,y_p^e), \ 
    e = e(1,l), \ p = p(1,k), \
    i=(l-1)Q_2+k, \ j=(e-1)M+q, 
    \ \text{for}\ 1\leqslant (l,k)\leqslant (N_2,Q_2); \\
    h_{ij}^{\text{bc2}}=& \mathcal{B}\phi_{eq}(b_1,y_p^e), \ 
    e = e(N_1,l), \ p = p(Q_1,k), \
    i=(l-1)Q_2+k, \ j=(e-1)M+q, 
    \ \text{for}\ 1\leqslant (l,k)\leqslant (N_2,Q_2); \\
    h_{ij}^{\text{bc3}}=& \mathcal{B}\phi_{eq}(x_p^e,a_2), \ 
    e = e(m,1), \ p = p(k,1), \
    i=(m-1)Q_1+k, \ j=(e-1)M+q, 
    \ \text{for}\ 1\leqslant (m,k)\leqslant (N_1,Q_1); \\
    h_{ij}^{\text{bc4}}=& \mathcal{B}\phi_{eq}(x_p^e,b_2), \ 
    e = e(m,N_2), \ p = p(k,Q_2), \
    i=(m-1)Q_1+k, \ j=(e-1)M+q, 
    \ \text{for}\ 1\leqslant (m,k)\leqslant (N_1,Q_1),
  \end{split}
  \right.
\end{equation}
where the functions $e(\cdot,\cdot)$ and $p(\cdot,\cdot)$
are given by~\eqref{eq_7} and~\eqref{eq_10}.
In the matrices $\mbs H^{\text{ck1}}$, $\mbs H^{\text{ck2}}$, $\mbs H^{\text{ck3}}$
and $\mbs H^{\text{ck4}}$ the only non-zero terms are,
\begin{equation}\label{eq_ab4}\footnotesize
  \left\{
  \begin{split}
    h_{ij_1}^{\text{ck1}}=& \phi_{e_1q}(X_m,y_{p_1}^{e_1}), \
    h_{ij_2}^{\text{ck1}}= -\phi_{e_2q}(X_m,y_{p_2}^{e_2}), \
    e_1 = e(m,l), \ p_1 = p(Q_1,k), \
    e_2 = e(m+1,l), \\& p_2 = p(1,k), \
    i=(m-1)N_2Q_2+(l-1)Q_2+k, \ j_1=(e_1-1)M+q, \
    j_2=(e_2-1)M+q, \\&
    \text{for}\ 1\leqslant (m,l,k,q)\leqslant (N_1-1,N_2,Q_2,M); \\
    h_{ij_1}^{\text{ck2}}=& \left.\frac{\partial\phi_{e_1q}}{\partial x}\right|_{(X_m,y_{p_1}^{e_1})}, \
    h_{ij_2}^{\text{ck2}}= -\left.\frac{\partial\phi_{e_2q}}{\partial x}\right|_{(X_m,y_{p_2}^{e_2})}, \
    e_1 = e(m,l), \ p_1 = p(Q_1,k), \
    e_2 = e(m+1,l), \\& p_2 = p(1,k), \
    i=(m-1)N_2Q_2+(l-1)Q_2+k, \ j_1=(e_1-1)M+q, \
    j_2=(e_2-1)M+q, \\&
    \text{for}\ 1\leqslant (m,l,k,q)\leqslant (N_1-1,N_2,Q_2,M); \\
    h_{ij_1}^{\text{ck3}}=& \phi_{e_1q}(x_{p_1}^{e_1},Y_l), \
    h_{ij_2}^{\text{ck3}}= -\phi_{e_2q}(x_{p_2}^{e_2},Y_l), \
    e_1 = e(m,l), \ p_1 = p(k,Q_2), \
    e_2 = e(m,l+1), \\& p_2 = p(k,1), \
    i=(l-1)N_1Q_1+(m-1)Q_1+k, \ j_1=(e_1-1)M+q, \
    j_2=(e_2-1)M+q, \\&
    \text{for}\ 1\leqslant (m,l,k,q)\leqslant (N_1,N_2-1,Q_1,M); \\
    h_{ij_1}^{\text{ck4}}=& \left.\frac{\partial\phi_{e_1q}}{\partial y}\right|_{(x_{p_1}^{e_1},Y_l)}, \
    h_{ij_2}^{\text{ck4}}= -\left.\frac{\partial\phi_{e_2q}}{\partial y}\right|_{(x_{p_2}^{e_2},Y_l)}, \
    e_1 = e(m,l), \ p_1 = p(k,Q_2), \
    e_2 = e(m,l+1), \\& p_2 = p(k,1), \
    i=(l-1)N_1Q_1+(m-1)Q_1+k, \ j_1=(e_1-1)M+q, \
    j_2=(e_2-1)M+q, \\&
    \text{for}\ 1\leqslant (m,l,k,q)\leqslant (N_1,N_2-1,Q_1,M).
  \end{split}
  \right.
\end{equation}

\section*{Appendix C.~Comparison with PINN}

This appendix provides a comparison of the simulation results obtained by
the current method (NLLSQ algorithm) and the physics-informed neural network (PINN)
method~\cite{RaissiPK2019} for several test 
problems from Section~\ref{sec:tests}.
The PINN method is also implemented in Python based on the Tensorflow
and Keras libraries. The PINN loss function consists of those contributions
from the parametric PDE, the measurement, and the boundary/initial conditions (BC/IC).
Let $\gamma_{bc}\in(0,1)$ denote the penalty coefficient in front of
the BC/IC loss term, and we employ $(1-\gamma_{bc})$ as the penalty coefficient
for the PDE and measurement loss terms. We have varied $\gamma_{bc}$,
the learning rate schedule, and the random initialization for the weights/biases
of PINN systematically. PINN is trained by the Adam optimizer.
The PINN/Adam results reported below are the best we have obtained for these problems using PINN.
We have also tried the L-BFGS optimizer with PINN, and its results for these inverse problems are
quite poor and worse than the Adam results.


\begin{table}[tb]
  \centering\small
  \begin{tabular}{l|cccc|cccc}
    \hline
    & & $\epsilon$=0 & & & & $\epsilon$=0.01 & &  \\ \cline{2-9}
    method & $e_{\alpha}$ & $l^{\infty}$-u & $l^2$-u & time(sec)
    & $e_{\alpha}$ & $l^{\infty}$-u & $l^2$-u & time(sec) \\ \hline
    PINN (Adam) & 6.31E-3 & 1.08E-2 & 3.56E-3 & 134.5 & 5.53E-3 & 1.03E-2 & 3.30E-3 & 130.9 \\ \hline
    current (NLLSQ) & 1.66E-8 & 3.66E-6 & 2.62E-7 & 11.5 & 9.72E-4 & 1.76E-3 & 5.26E-4 & 10.4 \\
    \hline
  \end{tabular}
  \caption{\small Inverse Poisson problem: relative errors of $\alpha$ and $u$ and the network training
    time (seconds) obtained by PINN (Adam) and the current NLLSQ algorithm.
    In both PINN and NLLSQ, Q=30$\times$30, $Q_s$=100, Gaussian activation function.
    In PINN, neural network $[2,30,30,30,1]$; $20,000$ training epochs; $\gamma_{bc}=0.99$;
    learning rate decreasing linearly from
    0.01 to 1.0E-4 in first $10,000$ epochs, and fixed at 1.0E-4 afterwards.
    In NLLSQ, single sub-domain, neural network $[2,500,1]$, $R_m$=3.0, $\lambda_{mea}$=0.1.
  }
  \label{tab_29}
\end{table}

\begin{table}[tb]
  \centering\small
  \begin{tabular}{l|cccc|cccc}
    \hline
    & & $\epsilon$=0 & & & & $\epsilon$=0.01 & &  \\ \cline{2-9}
    method & $e_{c}$ & $l^{\infty}$-u & $l^2$-u & time(sec)
    & $e_{c}$ & $l^{\infty}$-u & $l^2$-u & time(sec) \\ \hline
    PINN (Adam) & 1.18E-5 & 7.18E-3 & 7.63E-4 & 133.5 & 1.47E-4 & 9.15E-3 & 1.67E-3 & 134.9 \\ \hline
    current (NLLSQ) & 2.32E-10 & 8.51E-8 & 4.66E-9 & 29.6 & 2.61E-5 & 2.85E-4 & 1.10E-4 & 39.3 \\ \hline
    \hline
  \end{tabular}
  \caption{\small Inverse advection problem: relative errors of $c$ and $u$ and the network training
    time (seconds) obtained by the PINN (Adam) and the current NLLSQ algorithm.
    In both PINN and NLLSQ, Q=30$\times$30, $Q_s$=100, Gaussian activation function.
    In PINN, neural network $[2,30,30,30,30,1]$; $20,000$ training epochs; $\gamma_{bc}=0.2$;
    learning rate decreasing linearly from
    0.01 to 1.0E-4 in first $10,000$ epochs, and fixed at 1.0E-4 afterwards.
    In NLLSQ, single sub-domain, neural network $[2,400,1]$, $R_m$=2.5, $\lambda_{mea}$=0.1.
  }
  \label{tab_31}
\end{table}

\begin{table}[tb]
  \centering\small
  \begin{tabular}{l|l|ccccc}
    \hline
    noise level & method & $e_{\alpha_1}$ & $e_{\alpha_2}$ & $l^{\infty}$-u & $l^2$-u & training-time(sec) \\ \hline
    $\epsilon=0$ & PINN (Adam) & 7.08E-1 & 2.68E-1 & 1.48E+0 & 5.65E-1 & 3049.2 \\ \cline{2-7}
    & current (NLLSQ) & 5.71E-9 & 3.05E-7 & 5.98E-8 & 1.49E-8 & 10.3 \\ \hline
    $\epsilon=0.01$ & PINN (Adam) & 6.74E-1 & 7.76E-1 & 1.56E+0 & 6.79E-1 & 2742.9 \\ \cline{2-7}
    & current (NLLSQ) & 4.34E-3 & 7.13E-4 & 5.25E-3 & 2.38E-3 & 10.0 \\ \hline
    \hline
  \end{tabular}
  \caption{\small Inverse nonlinear Helmholtz problem: relative errors of $\alpha_1$, $\alpha_2$
    and $u$ and the network training
    time (seconds) obtained by PINN (Adam) and the current NLLSQ algorithm.
    In both PINN and NLLSQ, Q=30$\times$30, $Q_s$=100, Gaussian activation function.
    In PINN, neural network $[2,30,30,30,30,30,30,1]$; $200,000$ training epochs; $\gamma_{bc}=0.99$;
    learning rate decreasing linearly from
    0.01 to 1.0E-4 in first $10,000$ epochs, and fixed at 1.0E-4 afterwards.
    In NLLSQ, single sub-domain, neural network $[2,500,1]$, $R_m$=2.25, $\lambda_{mea}$=0.25.
  }
  \label{tab_33}
\end{table}

\begin{table}[tb]
  \centering\small
  \begin{tabular}{l|l|ccccc}
    \hline
    noise level & method & $e_{\alpha_1}$ & $e_{\alpha_2}$ & $l^{\infty}$-u & $l^2$-u & training-time(sec) \\ \hline
    $\epsilon=0$ & PINN (Adam) & 2.40E-4 & 7.26E-4 & 1.59E-3 & 2.06E-4 & 529.0 \\ \cline{2-7}
    & current (NLLSQ) & 4.31E-10 & 2.33E-9 & 3.31E-9 & 5.86E-10 & 4.2  \\ \hline
    $\epsilon=0.01$ & PINN (Adam) & 1.80E-3 & 2.66E-3 & 5.65E-3 & 7.83E-4 & 540.1  \\ \cline{2-7}
    & current (NLLSQ) & 1.12E-5 & 1.23E-4 & 9.87E-5 & 3.75E-5 & 6.4 \\ 
    \hline
  \end{tabular}
  \caption{\small Inverse Burgers' problem: relative errors of $\alpha_1$, $\alpha_2$
    and $u$ and the network training
    time (seconds) obtained by PINN (Adam) and the current NLLSQ algorithm.
    In both PINN and NLLSQ, Q=30$\times$30, $Q_s$=100, Gaussian activation function.
    In PINN, neural network $[2,30,30,30,30,1]$; $50,000$ training epochs; $\gamma_{bc}=0.9$;
    learning rate decreasing linearly from
    0.01 to 1.0E-4 in first $10,000$ epochs, and fixed at 1.0E-4 afterwards.
    In NLLSQ, single sub-domain, neural network $[2,400,1]$, $R_m$=1.9, $\lambda_{mea}$=0.1.
  }
  \label{tab_34}
\end{table}

\begin{table}[tb]
  \centering\small
  \begin{tabular}{l|l|cccccc}
    \hline
    noise level & method & $e_{\alpha_1}$ & $e_{\alpha_2}$ & $e_{\alpha_3}$ & $l^{\infty}$-u & $l^2$-u & training-time(sec) \\ \hline
    $\epsilon=0$ & PINN (Adam) & 9.21E-3 & 2.30E-1 & 7.33E-1 & 1.86E-2 & 3.32E-3 & 1853.3 \\ \cline{2-8}
    & current (NLLSQ) & 7.65E-10 & 4.49E-9 & 6.60E-9 & 7.97E-10 & 3.50E-10 & 23.6  \\ \hline
    $\epsilon=0.01$ & PINN (Adam) & 1.16E-2 & 1.35E-1 & 3.51E-1 & 1.18E-2 & 2.97E-3 & 1833.2 \\ \cline{2-8}
    & current (NLLSQ) & 5.45E-3 & 2.59E-2 & 5.09E-3 & 5.76E-3 & 2.63E-3 & 30.1 \\
    \hline
  \end{tabular}
  \caption{\small Inverse Sine-Gordan problem: relative errors of $\alpha_i$ ($i=1,2,3$)
    and $u$ and the network training
    time (seconds) obtained by PINN (Adam) and the current NLLSQ algorithm.
    In both PINN and NLLSQ, Q=30$\times$30, $Q_s$=100, Gaussian activation function.
    In PINN, neural network $[2,30,30,30,30,1]$; $200,000$ training epochs; $\gamma_{bc}=0.99$,
    learning rate decreasing linearly from
    0.01 to 1.0E-4 in first $10,000$ epochs, and fixed at 1.0E-4 afterwards.
    In NLLSQ, single sub-domain, neural network $[2,400,1]$, $R_m$=1.5, $\lambda_{mea}$=0.01.
  }
  \label{tab_35}
\end{table}

Tables~\ref{tab_29} through~\ref{tab_35} summarize the errors of the inverse parameters
and the solution field, as well as the network training time, obtained by
the current and the PINN methods for the inverse Poisson,
advection, nonlinear Helmholtz, Burgers', and
the Sine-Gordan problems. The table captions provide the respective parameter values
in these simulations for the two methods.
We observe that the current method produces more accurate results than PINN for both
the inverse parameters and the solution field, and that the network training time
of the current method is markedly smaller than that of PINN.
For the noise-free data, the current method is significantly more accurate
(typically by several orders of magnitude) than PINN.

\section*{Appendix D.~Parameter Values in Algorithm~\ref{alg_1} for Numerical Tests}

\noindent{\bf Section~\ref{poisson} (Parametric Poisson Equation):}\\
\noindent\underline{For NLLSQ:} \\
In Figure~\ref{fg_3}, Tables~\ref{tab_1}, \ref{tab_2}, \ref{tab_3}, \ref{tab_5},
and \ref{tab_6}: (max-nllsq-iterations,max-sub-iterations,$\varepsilon$,$\delta$,$\eta$,$\bm\theta_0$)=(80,2,1E-8,1.0,1,0). \\
In Table~\ref{tab_4}: (max-nllsq-iterations,max-sub-iterations,$\varepsilon$,$\delta$,$\eta$,$\bm\theta_0$)=(80,5,1E-8,1.0,1,0). \\
\noindent\underline{For VarPro-F1:}\\
In Tables~\ref{tab_2} and~\ref{tab_6}: (max-nllsq-iterations,max-sub-iterations,$\varepsilon$,$\delta$,$\eta$,$\bm\theta_0$)=(80,2,1E-8,1.0,1,0). \\
In Table~\ref{tab_3}: (max-nllsq-iterations,max-sub-iterations,$\varepsilon$,$\delta$,$\eta$,$\bm\theta_0$)=(80,5,1E-8,1.0,1,0). \\
In Table~\ref{tab_4}: (max-nllsq-iterations,max-sub-iterations,$\varepsilon$,$\delta$,$\eta$,$\bm\theta_0$)=(100,2,1E-8,4.0,1,0). \\
\noindent\underline{For VarPro-F2:}\\
In Tables~\ref{tab_2}, \ref{tab_3}, and \ref{tab_6}: (max-nllsq-iterations,max-sub-iterations,$\varepsilon$,$\delta$,$\eta$,$\bm\theta_0$)=(80,2,1E-8,1.0,1,0). \\
In Table~\ref{tab_4}: (max-nllsq-iterations,max-sub-iterations,$\varepsilon$,$\delta$,$\eta$,$\bm\theta_0$)=(80,5,1E-8,1.0,1,0). 

\vspace{3pt}
\noindent{\bf Section~\ref{advec} (Parametric Advection Equation):}\\
\noindent\underline{For NLLSQ:}\\
In Figure~\ref{fg_4}, Tables~\ref{tab_7}, \ref{tab_8}, \ref{tab_9}, \ref{tab_10},
\ref{tab_11} and \ref{tab_011}:
(max-nllsq-iterations,max-sub-iterations,$\varepsilon$,$\delta$,$\eta$,$\bm\theta_0$)=(80,10,1E-8,10.0,0,$\bm\vartheta_0$). \\
\noindent\underline{For VarPro-F1:}\\
In Table~\ref{tab_8}: (max-nllsq-iterations,max-sub-iterations,$\varepsilon$,$\delta$,$\eta$,$\bm\theta_0$)=(80,2,1E-8,5.0,0,$\bm\vartheta_0$). \\
In Tables~\ref{tab_9} and \ref{tab_11}: (max-nllsq-iterations,max-sub-iterations,$\varepsilon$,$\delta$,$\eta$,$\bm\theta_0$)=(100,5,1E-8,5.0,0,$\bm\vartheta_0$). \\
\noindent\underline{For VarPro-F2:}\\
In Table~\ref{tab_8}: (max-nllsq-iterations,max-sub-iterations,$\varepsilon$,$\delta$,$\eta$,$\bm\theta_0$)=(80,2,1E-8,5.0,0,$\bm\xi_0$). \\
In Table~\ref{tab_9}: (max-nllsq-iterations,max-sub-iterations,$\varepsilon$,$\delta$,$\eta$,$\bm\theta_0$)=(80,2,1E-8,1.0,1,$\bm\xi_0$). \\
In Table~\ref{tab_11}: (max-nllsq-iterations,max-sub-iterations,$\varepsilon$,$\delta$,$\eta$,$\bm\theta_0$)=(80,5,1E-8,1.0,1,$\bm\xi_0$). \\
In the above, $\bm\xi_0$ is a uniform random vector from [-1,1].
$\bm\vartheta_0$ is a uniform random vector generated by the lines 7 through 14 of
Algorithm~\ref{alg_1} with the $\delta$ as specified above and $\eta=0$.

\vspace{3pt}
\noindent{\bf Section~\ref{nonl_helm} (Parametric Nonlinear Helmholtz Equation):}\\
\noindent\underline{For NLLSQ:}\\
In Tables~\ref{tab_12}, \ref{tab_13}, \ref{tab_14}, \ref{tab_15}, \ref{tab_16}:
(max-nllsq-iterations,max-sub-iterations,$\varepsilon$,$\delta$,$\eta$,$\bm\theta_0$)=(80,2,1E-8,0.5,1,0). \\
\noindent\underline{For VarPro-F1:}\\
In Figure~\ref{fg_5}, Tables~\ref{tab_13}, \ref{tab_14} and \ref{tab_16}:
(max-nllsq-iterations,max-sub-iterations,$\varepsilon$,$\delta$,$\eta$,$\bm\theta_0$)=(80,2,1E-8,0.5,1,0). \\
\noindent\underline{For VarPro-F2:}\\
In Tables~\ref{tab_13}, \ref{tab_14} and \ref{tab_16}:
(max-nllsq-iterations,max-sub-iterations,$\varepsilon$,$\delta$,$\eta$,$\bm\theta_0$)=(80,0,1E-8,0.5,1,0); max-newton-iterations=15.

\vspace{3pt}
\noindent{\bf Section~\ref{burger} (Parametric Viscous Burgers' Equation):}\\
\noindent\underline{For NLLSQ:}\\
In Figure~\ref{fg_6}, Tables~\ref{tab_17}, \ref{tab_18}, \ref{tab_19} and \ref{tab_20}:
(max-nllsq-iterations,max-sub-iterations,$\varepsilon$,$\delta$,$\eta$,$\bm\theta_0$)=(80,2,1E-8,0.5,1,0).\\
\noindent\underline{For VarPro-F1:}\\
In Tables~\ref{tab_18}, \ref{tab_19} and \ref{tab_20}:
(max-nllsq-iterations,max-sub-iterations,$\varepsilon$,$\delta$,$\eta$,$\bm\theta_0$)=(80,2,1E-8,1.0,1,0).\\
\noindent\underline{For VarPro-F2:}\\
In Tables~\ref{tab_18}, \ref{tab_19} and \ref{tab_20}:
(max-nllsq-iterations,max-sub-iterations,$\varepsilon$,$\delta$,$\eta$,$\bm\theta_0$)=(80,2,1E-12,1.0,0,$\bm\xi_0$); max-newton-iterations=15.
Here $\bm\xi_0$ is a uniform random vector from [-1,1].

\vspace{3pt}
\noindent{\bf Section~\ref{sg_sec} (Parametric Sine-Gordan Equation):}\\
\noindent\underline{For NLLSQ:}\\
In Tables~\ref{tab_21}, \ref{tab_22}, \ref{tab_23}, \ref{tab_24} and \ref{tab_024}:
(max-nllsq-iterations,max-sub-iterations,$\varepsilon$,$\delta$,$\eta$,$\bm\theta_0$)=(80,5,1E-8,5.0,0,0).\\
\noindent\underline{For VarPro-F1:}\\
In Tables~\ref{tab_22}, \ref{tab_23}, \ref{tab_24} and \ref{tab_024}:
(max-nllsq-iterations,max-sub-iterations,$\varepsilon$,$\delta$,$\eta$,$\bm\theta_0$)=(80,5,1E-8,5.0,0,0).\\
\noindent\underline{For VarPro-F2:}\\
In Figure~\ref{fg_8}, Tables~Tables~\ref{tab_22}, \ref{tab_23},
\ref{tab_24} and \ref{tab_024}:
(max-nllsq-iterations,max-sub-iterations,$\varepsilon$,$\delta$,$\eta$,$\bm\theta_0$)=(80,5,1E-8,1.0,0,0); max-newton-iterations=15. 

\vspace{3pt}
\noindent{\bf Section~\ref{varhelm} (Helmholtz Equation with Inverse Variable Coefficient):}\\
\noindent\underline{For NLLSQ:}\\
In Figure~\ref{fg_9}, Tables~\ref{tab_25} and \ref{tab_26}:
(max-nllsq-iterations,max-sub-iterations,$\varepsilon$,$\delta$,$\eta$,$\bm\theta_0$)=(80,2,1E-8,1.0,1,0).\\
\noindent\underline{For VarPro-F1:}\\
In Table~\ref{tab_26}: (max-nllsq-iterations,max-sub-iterations,$\varepsilon$,$\delta$,$\eta$,$\bm\theta_0$)=(80,2,1E-8,0.01,1,0).\\
\noindent\underline{For VarPro-F2:}\\
In Table~\ref{tab_26}: (max-nllsq-iterations,max-sub-iterations,$\varepsilon$,$\delta$,$\eta$,$\bm\theta_0$)=(50,2,1E-8,0.5,1,0).

\bibliographystyle{plain}
\bibliography{ml,elm,elm1,mypub,dnn,sem,obc,varpro}

\end{document}